\documentclass[a4paper,landscape,12pt]{article}
\usepackage{amsmath}
\usepackage{amsfonts}
\usepackage{amssymb}
\usepackage[english]{babel}
\usepackage{amsthm}
\usepackage{extarrows}
\usepackage{amscd}
\usepackage{tikz}
\usepackage{mathdots}
\usepackage{wasysym}
\usepackage{mathtools}

\usetikzlibrary{arrows,shapes}

\def\pf{\par\noindent {\bf Proof}~\par\noindent}

\newcommand{\mR}{\mathbb{R}}
\newcommand{\mC}{\mathbb{C}}
\newcommand{\mN}{\mathbb{N}}
\newcommand{\mZ}{\mathbb{Z}}
\newcommand{\mE}{\mathbb{E}}

\newcommand{\mcH}{\mathcal{H}}

\newcommand{\mcC}{\mathcal{C}}
\newcommand{\mcD}{\mathcal{D}}

\newcommand{\ux}{\underline{x}}
\newcommand{\uy}{\underline{y}}
\newcommand{\la}{\lambda}

\newcommand{\uom}{\underline{\omega}}
\newcommand{\uomr}{\vec{e}_r}

\newcommand{\up}{\underline{p}}

\newcommand{\uD}{\underline{D}}
\newcommand{\uE}{\underline{E}}
\newcommand{\uS}{\underline{S}}
\newcommand{\uB}{\underline{B}}
\newcommand{\uM}{\underline{M}}

\newcommand{\bZ}{{\bf Z}}

\newcommand{\p}{\partial}
\newcommand{\dirac}{\underline{\p}}

\newcommand{\eop}{\hfill$\square$}

\newcommand{\onehalf}{\frac{1}{2}}
\newcommand{\onethird}{\frac{1}{3}}

\newcommand{\invr}{\frac{1}{r}}
\newcommand{\invrsq}{\frac{1}{r^2}}
\newcommand{\invrcub}{\frac{1}{r^3}}
\newcommand{\invux}{\frac{1}{\ux}}

\newtheorem{lemma}{Lemma}
\newtheorem{proposition}{Proposition}
\newtheorem{definition}{Definition}
\newtheorem{remark}{Remark}

\newtheorem{example}{Example}

\newcommand{\Eqref}[1]{(\ref{#1})}

\begin{document}

\title{Two Families of Clifford Distributions Revisited}
\author{Fred Brackx}

\date{Clifford Research Group, Foundation Lab,\\Department of Electronics and Information Systems,\\
Faculty of Engineering and Architecture, Ghent University
}

\maketitle


\begin{abstract}
Two specific families of distributions in harmonic and Clifford analysis are further studied through a spherical co-ordinates approach. In particular actions involving spherical co-ordinates, such as the radial derivative and the multiplication and division by the radial distance, are computed, giving rise to two families of so--called signumdistributions, i.e. bounded linear functionals on a space of test functions showing a singularity at the origin.
\\

\noindent
Keywords: distribution, radial derivative, signumdistribution\\
MSC: 46F05, 46F10, 30G35
\end{abstract}


\newpage
\section{Introduction}
\label{intro}


In a series of papers \cite{fb1, fb2, differential, BDSch, multi, fourier,  dancing, overview, convolution, specific} several families of distributions in Euclidean space $\mR^m$ were thoroughly studied in the context of Clifford analysis, a direct and elegant generalization to higher dimension of the theory of holomorphic functions in the complex plane. At the heart of Clifford analysis lies the Dirac operator $\dirac = \sum_{j=1}^m \, e_j\, \p_{x_j}$, ($e_1, e_2, \ldots, e_m$) being an orthonormal basis in $\mR^m$, which may be seen as a Stein--Weiss projection of the gradient operator (see e.g.\ \cite{stein}). Similar to the Cauchy--Riemann operator in the complex plane, the Dirac operator linearizes the Laplace operator in $\mR^m$: $\dirac^2 = -\, \Delta$,  and its null--solutions with values in the Clifford algebra $\mR_{0,m}$  are termed {\em monogenic functions}. For an in--depth study of monogenic function theory we refer to e.g. \cite{dss}.\\

Of particular importance are the distribution families $T_\lambda$ and $U_\lambda$, $\lambda$ being a complex parameter. They are defined as follows, using spherical co-ordinates $\ux = r\, \uom, r = |\ux|, \uom \in S^{m-1}$, $S^{m-1}$ being the unit sphere in $\mR^m$.
\begin{definition}
For all $\lambda \in \mC$ and for all test functions $\varphi(\ux) \in \mcD(\mR^m)$ the distributions $T_\lambda$ and $U_\lambda$ are defined by
$$\langle  \ T_\lambda , \varphi(\ux)  \ \rangle := a_m\, \langle \  {\rm Fp}\, r^{\lambda+m-1}_+  , \Sigma^0[\varphi](r)  \  \rangle_r$$
and
$$\langle  \  U_\lambda , \varphi(\ux) \ \rangle := a_m\, \langle \  {\rm Fp}\, r^{\lambda+m-1}_+  , \Sigma^1[\varphi](r)  \  \rangle_r$$
\noindent
where the so--called spherical means $\Sigma^0$ and $\Sigma^1$ are given by
$$\Sigma^0[\varphi](r) = \frac{1}{a_m}\int_{S^{m-1}}\, \phantom{\uom}\, \varphi(r\uom)\, dS(\uom)$$
and
$$\Sigma^1[\varphi](r) = \frac{1}{a_m}\int_{S^{m-1}}\, \uom\, \varphi(r\uom)\, dS(\uom)$$
with $a_m = \frac{2\pi^{m/2}}{\Gamma(m/2)}$ the area of the unit sphere $S^{m-1}$,
and where {\rm Fp}$\, r_+^\mu$ stands for the {\em finite part} distribution on the one--dimensional $r$--axis.
\end{definition}

Despite the fact that the distributions $T_\lambda$ are spherical in nature and that a formula such as $\p_r\, T_\la = \la\, T_{\la-1}$ seems to be trivial, at the time the above mentioned papers  were written their radial derivative $\p_r\, T_\la$ was not yet studied.
Meanwhile it has become clear that derivation of a distribution with respect to the spherical co-ordinates  is far from trivial. To put it straight: the radial derivative of a distribution cannot be a distribution. Already in his famous and seminal book \cite{zwart} Laurent Schwartz writes on page 51:
{\em
Using co-ordinate systems other than the cartesian ones should be done with the utmost care
}
 [our translation]. As an illustration consider the delta-distribution $\delta(\ux)$: it is pointly supported at the origin, it is rotation invariant: $\delta(A\, \ux) = \delta(\ux), \ \forall \, A \in {\rm SO}(m)$, it is even:  $\delta(-\ux) = \delta(\ux)$ and it is homogeneous of degree $(-m)$: $\delta(a \ux) = \frac{1}{|a|^m}\, \delta(\ux)$. So in a first, naive, approach, one could think of its radial derivative $\p_r\, \delta(\ux)$ as being a distribution which remains pointly supported at the origin, rotation invariant, even and homogeneous of degree $(-m-1)$.  Temporarily leaving aside the even character, on the basis of the other cited characteristics the distribution $\p_r\, \delta(\ux)$ should take the form
$$
\p_r\, \delta(\ux) = c_0\, \p_{x_1} \delta(\ux) + \cdots + c_m\, \p_{x_m} \delta(\ux)
$$
and it becomes clear at once that such an approach to the radial derivation of the delta-distribution is impossible since all distributions appearing in the sum at the right--hand side are odd and not rotation invariant, whereas $\p_r\, \delta(\ux)$ is assumed to be even and rotation invariant. It could be that $\p_r\, \delta(\ux)$ is either the zero distribution or is no longer pointly supported at the origin, but both those possibilities are unacceptable.\\
In \cite{bsv} the problem of defining the radial derivative $\p_r\, \delta(\ux)$ of the delta--distribution was solved by introducing a new concept, similar to but fundamentally different from a distribution, which was termed a {\em signumdistribution}. Roughly speaking a signumdistribution is a continuous linear functional on a space of test functions which are smooth in $\mR^m \setminus \{O \}$ and show a non--removable singularity at the origin. The definition and first properties of signumdistributions are recalled in Section \ref{signumdistributions}. It turns out that the radial derivative of the delta distribution  is indeed a signumdistribution.\\
Actions on general distributions involving spherical co-ordinates  were thoroughly studied in \cite{ffb}. For a distribution $T \in \mcD'(\mR^m)$, its radial derivative $\p_r\, T$ is well--defined, but not uniquely defined, as an equivalence class of signumdistributions. To make the paper self--contained the main results about spherical actions on (signum)distributions are recalled in Sections \ref{signumdistributions},  \ref{cartoperators}, \ref{sigoperators} and \ref{spheroperators}. In Section \ref{uniqueness} we state sufficient conditions under which the actions of a number of operators involving division by $\ux$ or $r$ are uniquely determined. In Section \ref{actions} we will show a.o. that the radial derivative of the distributions $T_\lambda$ and $U_\lambda$ is uniquely defined as a signumdistribution out of two families of signumdistributions $^{s}T_\lambda$ and $^{s}U_\lambda$ which are closely related to the initial distribution families $T_\lambda$ and $U_\lambda$. Also other actions on the distributions $T_\lambda$ and $U_\lambda$ involving spherical co-ordinates are explored. Finally, in Sections \ref{signumcartderiv} and \ref{cartderiv}, special attention is paid to the cartesian derivatives of signumdistributions in general and the signumdistributions $^{s}T_\lambda$ and $^{s}U_\lambda$ in particular.


\newpage
\section{Signumdistributions}
\label{signumdistributions}


We confine ourselves to a sketchy introduction of the new concept of a signumdistribution; for a systematic treatment including the functional analytic justification and numerous examples, we refer to \cite{bsv, ffb}.\\

\noindent
In Euclidean space $\mR^m$ we interpret vectors as Clifford $1$--vectors in the Clifford algebra $\mR_{0,m}$, where the orthonormal basis vectors $(e_1, e_2, \ldots, e_m)$ satisfy the relations $e_j^2 = - 1,\, e_i \wedge e_j = e_ie_j = - e_je_i = - e_j \wedge e_i,\, e_i \cdot e_j = 0, i \neq j=1,\ldots,m$. This allows for the use of the very efficient {\em geometric} or {\em Clifford product} of Clifford vectors:
$$
\ux \,  \uy = \ux \cdot \uy  + \ux \wedge \uy
$$
for which, in particular,
$$
\ux \,  \ux = \ux \cdot \ux  = -\, |\ux|^2
$$
$\ux$ being the Clifford 1--vector $\ux = \sum_{j=1}^m \, e_j\, x_j$, whence also
$$
\uom \,  \uom = \uom \cdot \uom  = -\, |\uom|^2 = -1 \quad , \quad \uom \in S^{m-1}
$$
For more on Clifford algebras we refer to e.g.\ \cite{port}.\\

\noindent
We consider two spaces of test functions: the traditional space $\mcD(\mR^m)$ of compactly supported infinitely differentiable functions $\varphi(\ux)$ and the space $\Omega(\mR^m; \mR^m) = \{ \uom\, \varphi(\ux) : \varphi(\ux) \in \mcD(\mR^m)  \}$. Clearly the test functions in $\Omega(\mR^m; \mR^m)$ are no longer differentiable in the whole of $\mR^m$, since they are not defined at the origin due to the function $\uom = \frac{\ux}{|\ux|}$ which can be seen as the higher--dimensional counterpart of the one--dimensional signum function $sign(t) = \frac{t}{|t|}, t \in \mR$. Obviously there is a one--to--one correspondence between the spaces $\mcD(\mR^m)$ and $\Omega(\mR^m; \mR^m)$. The continuous linear functionals on those spaces of test functions, both equipped with an appropriate topology, are the standard distributions and the signumdistributions respectively.\\

\noindent
Given a standard distribution $T(\ux) \in \mcD'(\mR^m)$,  the signumdistribution  $T^{\vee}(\ux) \in \Omega'(\mR^m; \mR^m)$ is defined in such a way that for all test functions $\uom\, \varphi \in \Omega(\mR^m; \mR^m)$ it holds that
\begin{equation}
\label{unique}
\langle \  T^{\vee}(\ux) ,  \uom\, \varphi(\ux)  \ \rangle = -\, \langle \  T(\ux) ,  \varphi(\ux)  \ \rangle 
\end{equation}
Then $T^{\vee}(\ux)$ is called the signumdistribution associated to $T(\ux)$. In \cite{ffb} it was proven that this associated signumdistribution is {\em unique}.

\noindent
Conversely for a given signumdistribution $^{s}U \in \Omega'(\mR^m; \mR^m)$ we define the associated distribution $^{s}U^{\wedge}$ by
$$
\langle \  ^{s}U^{\wedge}(\ux) ,   \varphi(\ux)  \ \rangle = -\, \langle \  ^{s}U(\ux) ,  \uom\, \varphi(\ux)  \ \rangle \qquad \forall \varphi(\ux) \in \mcD(\mR^m)
$$
Clearly it holds that
$$
T^{\vee \wedge} = T \hspace{8mm} \mbox{and} \hspace{8mm} ^{s}U^{\wedge \vee} =\, ^{s}U
$$


\newpage
\section{Cartesian operators}
\label{cartoperators}


We call an operator acting on distributions {\em cartesian} if it involves partial derivation with respect to the cartesian co-ordinates, and multiplication and division by analytic functions. Although being well--defined, this last operation is not uniquely determined but results into an equivalence class of distributions involving (derivatives of) the delta distribution $\delta(\ux - \underline{y})$, $\underline{y}$ being a zero of the analytic function considered. Two basic cartesian operators on distributions are the multiplication operator $\ux$ and the Dirac operator $\dirac$. Their actions are, quite naturally, well--defined and uniquely determinded, as is the case for their squares: the multiplication operator $-\, \ux^2 =  |\ux|^2$ and the Laplace operator $-\, \dirac^2 = \Delta$.
Division of a distribution by $\ux$ may be approached in two  ways. Either the function $\ux$ is considered as a vector polynomial of the first degree, whence a vector analytic function  with a single zero at the origin, or, as $\invux = -\, \frac{\ux}{|\ux|^2}$, division by $\ux$ is seen as the composition of two operations: first division by $r^2 = |\ux|^2 = x_1^2 + \ldots + x_m^2$, followed by multiplication by $(-\, \ux)$. In \cite{ffb} the following lemma was proven showing that both approaches are equivalent.

\begin{lemma}
\label{divisionbyx}
For a scalar distribution $T$ it holds that
$$
\invux\, T = \underline{S} + \delta(\ux)\, \underline{c}
$$
for any distribution $\underline{S}$ such that $\ux\, \underline{S} = T$.
\end{lemma}
\noindent
Henceforth we use the notation $\left[  \invux\, T \right]$ for the equivalent class of distributions $S$ such that $\ux\, S = T$.\\

\noindent
The two fundamental formulae in monogenic function theory involving the (anti--)commutator of $\ux$ and $\dirac$, viz.
$$
\{ \ux , \dirac\} = \ux \, \dirac + \dirac \, \ux =  -\, 2 \mE - m \hspace{8mm} \mbox{and} \hspace{8mm}  [ \ux , \dirac ] = \ux \, \dirac - \dirac \, \ux = m - 2 \Gamma
$$
give rise to two other well--known cartesian operators: the scalar Euler operator 
$$
\mE = \sum_{j=1}^m \, x_j \, \p_{x_j}
$$ 
and the  bivector angular momentum operator 
$$
\Gamma = -\, \sum_{j<k} \, e_j e_k\, L_{jk} =  -\, \sum_{j<k} \, e_j e_k (x_j \p{x_k} - x_k \p_{x_j})
$$ 
It follows that
$$
\ux\, \dirac = -\, \mE - \Gamma
$$
or more precisely
\begin{equation}
\label{eulergamma}
 \ux  \cdot \dirac  =  -\, \mE \hspace{8mm}  \mbox{and} \hspace{8mm} \ux \wedge \dirac  = -\, \Gamma
\end{equation}
Passing to spherical co-ordinates $\ux = r \uom,\, r = |\ux|,\, \uom = \sum_{j=1}^{m} \, e_j \, \omega_j \, \in S^{m-1}$, the Dirac operator takes the form
$$
\dirac = \dirac_{rad} + \dirac_{ang}
$$
with
$$
\dirac_{rad} = \uom\, \p_r   \hspace{8mm} \mbox{and}  \hspace{8mm}  \dirac_{ang} = \frac{1}{r}\, \p_{\uom}
$$

\noindent
Taking into account that $\p_{\uom}$ is orthogonal to $\uom$, the Euler operator takes the well--known form 
$$
\mE = -\,  \ux  \cdot \dirac   =  -\, r\uom \cdot \dirac_{rad} = -\, r\uom \cdot \uom\, \p_r = r \, \p_r
$$
while the angular momentum operator $\Gamma$ takes the form 
$$
\Gamma = -\,  \ux  \wedge \dirac  =  -\, r\uom \wedge \dirac_{ang} = -\, r\uom \wedge  \frac{1}{r}\, \p_{\uom} = -\, \uom \wedge \p_{\uom} = - \uom \, \p_{\uom} =  -\, \sum_{j<k} \, e_j e_k (\omega_j \p{\omega_k} - \omega_k \p_{\omega_j})
$$
The question now is how to define, if possible, the separate action of the operators $\dirac_{rad}$ and $\dirac_{ang}$ on a distribution. To that end both operators should first be shown to be cartesian, which is achieved by putting
$$
\dirac_{rad} = \uom\, \p_r = -\, \frac{1}{\ux} \, \mE \hspace{8mm} \mbox{and} \hspace{8mm} \dirac_{ang} = \frac{1}{r}\, \p_{\uom} = -\, \frac{1}{\ux} \, \Gamma
$$
leading to the following definition.

\begin{definition}
\label{defdiracparts}
Let $T(\ux) \in \mcD'(\mR^m)$ be a distribution. Then we put
\begin{equation}
\label{equirad}
\dirac_{rad}\, T(\ux) = (\uom\, \p_r)\, T(r \uom) =  -\, \left[\frac{1}{\ux} \, \mE\, T(\ux)\right] 
\end{equation}
and
\begin{equation}
\label{equiang}
\dirac_{ang}\, T(\ux) =  (\frac{1}{r}\, \p_{\uom})\, T(r\uom) = -\, \left[ \frac{1}{\ux} \, \Gamma\, T(\ux)\right]
\end{equation}
\end{definition}

\noindent
It immediately becomes clear that, in this way, the actions of $\dirac_{rad}$ and $\dirac_{ang}$ on the distribution $T(\ux)$ are well-defined but not uniquely determined, due to the division by $\ux$, see Lemma \ref{divisionbyx}. The resulting equivalent classes are denoted by square brackets.
However, if $\underline{S}_1$ and $\underline{S}_2$ are distributions arbitrarily chosen in the equivalent classes (\ref{equirad}) and (\ref{equiang})  respectively, i.e.
$$
\ux \, \underline{S}_1 = -\, \mE\, T(\ux)  \hspace{8mm} \mbox{and} \hspace{8mm} \ux\, \underline{S}_2 = -\, \Gamma T(\ux)
$$
then
\begin{eqnarray*}
\dirac_{rad}\, T(\ux) &=& \underline{S}_1 + \underline{c}_1\, \delta(\ux)\\
\dirac_{ang}\, T(\ux) &=& \underline{S}_2 + \underline{c}_2\, \delta(\ux)
\end{eqnarray*}
and it must hold that
\begin{equation}
\label{entangle}
\underline{S}_1 + \underline{c}_1\, \delta(\ux) + \underline{S}_2 + \underline{c}_2\, \delta(\ux) = \dirac_{rad}\, T(\ux) + \dirac_{ang}\, T(\ux) = \dirac\, T(\ux)
\end{equation}
where the distribution at the utmost right--hand side is, quite naturally, a known distribution once the distribution $T$ is given.
One could say that the differential operators $\dirac_{rad}$ and $\dirac_{ang}$ are {\em entangled} in the sense that the results of their actions on a distribution are subject to (\ref{entangle}) which becomes a condition on the arbitrary constants $\underline{c}_1$ and $\underline{c}_2$. Henceforth we call (\ref{entangle}) the {\em entanglement condition} for the operators $\dirac_{rad}$ and $\dirac_{ang}$.\\

\noindent
Expressing the Laplace operator in spherical co-ordinates we obtain,
in view of
\begin{align}
\dirac_{rad} \, \dirac_{rad} &= -\, \p_{r}^2 \nonumber \\
\dirac_{rad} \, \dirac_{ang} &= -\frac{1}{r^2}\, \uom\, \p_{\uom} + \frac{1}{r}\, \uom\, \p_{\uom}\, \p_r \nonumber \\
\dirac_{ang} \, \dirac_{rad} &= -\, (m-1) \frac{1}{r}\, \p_r - \frac{1}{r}\, \p_r\, \uom\, \p_{\uom} \nonumber \\
\dirac_{ang} \, \dirac_{ang} &= \frac{1}{r^2}\, \p_{\uom}^2 \nonumber
\end{align}
that
$$
\Delta = -\, \dirac^2 = \p_r^2 + (m-1)\, \invr\, \p_r + \invrsq\, \Delta^*
$$
The so--called Laplace--Beltrami operator $\Delta^* = \uom\, \p_{\uom} - \p_{\uom}^2$ is, contrary to its appearance, cartesian, and so is the operator $ \p_{\uom}^2$. The following result indeed holds.
\begin{proposition} (see \cite{ffb})
\label{propangop}
The angular differential operators $\p_{\uom}^2$ and $\Delta^*$ may be written in terms of cartesian derivatives as
$$
\p_{\uom}^2 = \Gamma^2 - (m-1)\, \Gamma
$$
and
$$
\Delta^* = (m-2)\, \Gamma - \Gamma^2
$$
\end{proposition}

\noindent
The actions of the Laplace operator and the Laplace--Beltrami operator on a distribution being uniquely well-defined, the question arises how to define the separate actions on a distribution of the three parts of the Laplace operator expressed in spherical co-ordinates. It turns out that these operators are cartesian,  their actions on a distribution being well--defined, though not uniquely determined, through equivalent classes of distributions.

\begin{proposition}
\label{laplaceparts}
The operators $\p_r^2$, $\invr\, \p_r$ and $\invrsq\, \Delta^*$ are cartesian, and it holds, for a distribution T, that
\begin{eqnarray*}
\p_r^2\, T &=& \left[\,  \invrsq\, \mE\, (\mE - 1)\, T\, \right]\\
\invr\, \p_r\, T &=& \left[\, \invrsq\, \mE\, T\, \right]\\
\invrsq\, \Delta^*\, T &=& \left[\, \invrsq\, ( (m-2)\Gamma - \Gamma^2)\, T\, \right]
\end{eqnarray*}
leading to
\begin{itemize}
\item[(i)] $\p_r^2\, T = S_2 + \delta(\ux)\, c_2 - \sum_{j=1}^m\, c_{1,j}\, \p_{x_j} \delta(\ux)$\\[1mm] for arbitrary constants $c_2$ and $ c_{1,j}, j=1,\ldots,m$ and any distribution $S_2$ such that $\ux\, S_2 = \mE\, \underline{S}_1$ with $\ux\, \underline{S}_1 = -\, \mE\, T$
\item[(ii)] $\invr\, \p_r\, T = S_3 + \frac{1}{m}\, \sum_{j=1}^m\, c_{1,j}\, \p_{x_j} \delta(\ux) + c_3\, \delta(\ux)$\\[1mm] for arbitrarily constant $c_3$ and any distribution $S_3$ such that $\ux\, S_3 = \underline{S}_1$
\item[(iii)] $\frac{1}{r^2}\ \Delta^*\, T = S_4 + c_4\, \delta(\ux) + \sum_{j=1}^m\, c_{5,j}\, \p_{x_j} \delta(\ux)$\\[1mm] for arbitrary constants $c_4$ and $ c_{5,j}, j=1,\ldots,m$ and any distribution $S_4$ such that $r^2\, S_4 = \Delta^*\, T$
\end{itemize}
\end{proposition}

\pf
\begin{itemize}
\item[(i)] A direct computation shows  that $r^2\, \p_r^2 = \mE(\mE-1)$, whence
$$
\p_r^2 = \invrsq\, \mE\, (\mE - 1)
$$
 Further we have
$$
(\uom\, \p_r)\, T = -\, \left[ \frac{1}{\ux}\, \mE\, T  \right] = \underline{S}_1 + \delta(\ux)\, \underline{c}_1 
$$
with $\ux\, \underline{S}_1 = -\, \mE\, T$. It follows that
\begin{align}
\p_r^2\, T &= -\, (\uom\, \p_r)^2\, T   \nonumber\\
&= -\, (\uom\, \p_r)\, (\underline{S}_1 + \delta(\ux)\, \underline{c}_1) \nonumber\\
&= \left[ \frac{1}{\ux}\, \mE\, \underline{S}_1  \right] - \dirac \delta(\ux)\, \underline{c}_1 \nonumber\\
&= S_2 + \delta(\ux)\, c_2 - \dirac \delta(\ux)\, \underline{c}_1 \nonumber
\end{align}
with $\ux\, S_2 = \mE\, \underline{S}_1$.
\item[(ii)] As $r\, \p_r = \mE$, it follows that $\invr\, \p_r = \invrsq\, \mE$ and also 
\begin{align}
\invr\, \p_r\, T &= \frac{1}{\ux}\, (\uom\, \p_r)\, T \nonumber\\
&= \frac{1}{\ux}\, ( \underline{S}_1 + \delta(\ux)\, \underline{c}_1 ) \nonumber\\
&= S_3 +  \frac{1}{\ux}\, \delta(\ux)\, \underline{c}_1 \nonumber\\
&= S_3 + \frac{1}{m}\, \dirac\, \delta(\ux)\,  \underline{c}_1+ \delta(\ux)\, c_3 \nonumber
\end{align}
with $\ux\, S_3 = \underline{S}_1$.
\item[(iii)] The distribution $\Delta^*\, T$ being uniquely well--defined and $r^2 $ being an analytic function with a second order zero at the origin, the result follows immediately.
\end{itemize} \eop

\begin{remark}
{\rm
The operators $\p_r^2$,  $\invr\, \p_r$ and $\frac{1}{r^2}\ \Delta^*$ are {\em entangled} in the sense that, given a distribution $T$ and having chosen appropriately the distributions $\underline{S}_1$, $S_2$, $S_3$ and $S_4$, all arbitrary constants appearing in the expressions of Proposition \ref{laplaceparts} should satisfy the entanglement condition  generated by
$$
\p_r^2\, T + (m-1) \frac{1}{r}\, \p_r\, T + \frac{1}{r^2}\, \Delta^*\, T = \Delta\, T
$$ 
the distribution at the right--hand side being uniquely determined.
}
\end{remark}


\newpage
\section{Signum--pairs and cross--pairs of operators}
\label{sigoperators}


With each well--defined operator $P$ acting between distributions there corresponds an operator $$P^{\vee} = \uom\, P\, (-\uom)$$ acting between signumdistributions according to the following commutative diagram
$$
\begin{array}{ccccccccc}
  &&    & {}_{P}   &                             &&&&\\
  &&  T & \longrightarrow  &  P\,   T &&&& \\[1mm]
  &&\hspace{-8mm}{}_{- \uom}&&&&&&\\[-2mm]
  && \uparrow &   & | &&&&\\[-1.1mm]
  && | &   & \downarrow &&&&\\[-2mm]
  &&&&\hspace{5mm}{}^{\uom}&&&&\\[1mm]
  && T^{\vee} = \uom\, T & \longrightarrow & P^{\vee}\, T^{\vee} = \uom\, P\, T &&&&\\[-1mm]
  &&               & {}_{P^{\vee}}  &           &&&&
 \end{array}
$$
in this way giving rise to a pair of operators $P$ and $P^{\vee}$ which we call a {\em signum--pair of operators}. If the action result of the operator $P$ is uniquely determined then the action result of $P^{\vee}$ is uniquely determined too, in which case we use the notation $(P,P^{\vee})$ for this signum--pair of operators. If, on the contrary, the action result of $P$ is an equivalence class of distributions, then the action result of $P^{\vee}$ will be an equivalence class of signumdistributions, in which case we use the notation $[P, P^{\vee}]$. In Table 1 a number of signum--pairs of operators are listed.

\begin{remark}
{\rm
When the operators $P$ and $P^{\vee}$ form a signum--pair of operators it then follows that $P = (-\uom)\, P^{\vee}\, \uom$, and it becomes tempting to consider the operator $P$ as the signum--partner to the operator $P^{\vee}$, and to define the action of the operator $P$ on signumdistributions through the action of $P^{\vee}$ on distributions. However this is justifiable only if the operator $P^{\vee}$ is a ''legal`` operator, such as a cartesian operator,  between distributions.
}
\end{remark}

The above commutative diagram induces two more operators: the operator $Q$ mapping a distribution to a signumdistribution, and the corresponding operator $Q^c = (-\uom)\, Q\, (-\uom) =  \uom\, Q\, \uom$ mapping a signumdistribution to a distribution according to the following commutative diagram
$$
\begin{array}{ccccccccc}
  &&    & {}_{P}   &                             &&&&\\
  &&  ^{s}U^{\wedge} = -\uom\,  ^{s}U & \longrightarrow  &  Q^c \, ^{s}U &&&& \\[1mm]
  &&\hspace{-8mm}{}_{- \uom} &\hspace{2mm}{} \hspace{11mm}{}_{Q^c}&\hspace{6mm}{}_{-\uom}&&&&\\[-2mm]
  && \uparrow & \diagdown \hspace{-1.3mm} \nearrow & \uparrow &&&&\\[-1.1mm]
  && | & \diagup \hspace{-1.3mm} \searrow & | &&&&\\[-2mm]
  &&&\hspace{2mm}{} \hspace{12mm}{}^{Q}&&&&&\\[1mm]
  &&   ^{s}U & \longrightarrow & Q \, (-\uom)\, ^{s}U &&&&\\[-1mm]
  &&               & {}_{P^{\vee}}  &           &&&&
 \end{array}
$$
Clearly the operators $Q$ and $Q^c$ cannot be cartesian since they map distributions to signumdistributions and vice versa. We call the pair of operators $Q$ and  $Q^c$  a {\em cross--pair of operators} denoted by either $(Q\, , Q^c)$ or $[Q\, , Q^c]$ depending on the nature of their action result similarly as in the case of a signum--pair of operators. In Table 2 a number of cross--pairs of operators are listed.\\

In Section \ref{cartoperators} we saw that the Laplace--Beltrami operator $\Delta^*$ and the square of the angular derivative $\p_{\uom}^2$ are cartesian operators.
Their signum--partners are straightforwardly computed to be
\begin{eqnarray*}
(\p_{\uom}^2)^{\vee} = \uom\, \p_{\uom}^2\, (-\, \uom) =  \p_{\uom}^2
\end{eqnarray*}
and
\begin{eqnarray*}
\Delta^{* \vee} = -\, \Gamma^2 + m\, \Gamma - (m-1)\,{\bf 1}
\end{eqnarray*}
Clearly also the operator $\Delta^{*\, \vee}$ is cartesian. Introducing  the notation ${\bf Z}^* = \Delta^{* \vee}$, the signum--pairs of operators $(\p_{\uom}^2 \, , \, \p_{\uom}^2)$, $( \Delta^* \, , \, {\bf Z}^* )$ and $( {\bf Z}^* \, , \, \Delta^* )$ follow, inducing the definition of the actions of the operators $\p_{\uom}^2$, ${\bf Z}^*$ and $\Delta^{*}$ on a signumdistribution.\\

\noindent
For the signum--partner $\underline{D}$ of the Dirac operator $\dirac$ we obtain the following expressions:
\begin{eqnarray*}
\underline{D} &=& \uom\, \dirac\, (-\uom) =  \uom\, (\uom\, \p_r + \invr\, \p_{\uom})\, (-\uom)\\
                      &=& \uom\, \p_r + \invr\, \uom\, \p_{\uom}\, (-\uom)\\
                      &=& \uom\, \p_r - \invr\,  \p_{\uom} + (m-1)\, \invr\, \uom\\
	             &=& \dirac - 2 \invr\,  \p_{\uom} + (m-1)\, \invr\, \uom\\
                      &=& -\, \dirac + 2\, \uom\, \p_r + (m-1)\, \invr\, \uom
\end{eqnarray*}
giving rise to the signum--pairs of operators $(\dirac\, ,\, \underline{D})$ and $[\invr\, \p_{\uom} \, , \, - \invr\,  \p_{\uom} + (m-1)\, \invr\, \uom]$, which  induce the actions of the operators $\uD$ and $- \invr\,  \p_{\uom} + (m-1)\, \invr\, \uom$ on a signumdistribution.\\
Notice that while the actions of the operator $\dirac$ on distributions and of its signum--partner $\underline{D}$ on signumdistributions are uniquely determined, the action results of the operator $\invr\, \p_{\uom}$ on distributions and of its signum--partner $- \invr\, \p_{\uom} + (m-1)\, \invr\, \uom$ on signumdistributions are equivalent classes. \\
We call the operator $\uD$ the {\em signum--Dirac operator}. At first sight it is not clear if $\uD$ is cartesian. But it is indeed,  and we have the following result.

\begin{proposition}
The operators $\uD$ and $- \invr\,  \p_{\uom} + (m-1)\, \invr\, \uom$ are cartesian operators.
\end{proposition}
\pf
In view of Definition \ref{defdiracparts} it holds that
$$
- \invr\, \p_{\uom} + (m-1)\, \invr\, \uom = \invux\, (\Gamma - (m-1){\bf 1})
$$
and also
\begin{eqnarray*}
\uD &=&  \invux\, ( -\, \mE + \Gamma - (m-1){\bf 1})\\
       &=&  \dirac + \invux\, (2\, \Gamma  - (m-1){\bf 1})\\
       &=& -\, \dirac + \invux\, ( -2\, \mE - (m-1){\bf 1})
\end{eqnarray*}
\eop

\noindent
This leads to the signum--pair of operators $[\uD\, , \, \dirac]$, which induces the action of the Dirac operator $\dirac$ on a signumdistribution. However, due to division by $\ux$, the latter action results into an equivalence class of signumdistributions.\\

\noindent
It is interesting to note that, in the same way as the Dirac operator factorizes the Laplace operator: $\dirac^2 = -\, \Delta$, the signum--Dirac operator $\uD$ factorizes the {\em signum--Laplace operator}, i.e. the signum--partner of the Laplace operator:
$$
\uD^2 = \left(\uom\, \dirac\, (-\uom)\right)^2 = \uom\, \dirac^2 \, (-\, \uom) = -\, \uom\, \Delta\, (-\, \uom) = - \Delta^{\vee}
$$
Introducing the notation ${\bf Z} = \Delta^{\vee}$, it follows that $(\Delta\, , \, {\bf Z})$ is a signum--pair of operators, with
\begin{eqnarray*}
{\bf Z} &=& -\, \uD^2\\
&=& \p_r^2 + (m-1)\, \invr\, \p_r + \invrsq\, {\bf Z}^*
\end{eqnarray*}
Clearly also the operator $\bZ$ is cartesian; the signum--pair of operators $[\bZ\, , \, \Delta]$ follows, inducing the action of the Laplace operator $\Delta$ on a signumdistribution.\\

\noindent
From Section \ref{cartoperators} we know that the operators $\p_r^2$, $\invr\, \p_r$ and $\invrsq\, \Delta^*$, which are the constituents of the Laplace operator $\Delta$, are cartesian operators. Their signum--partners are easily seen to be
\begin{eqnarray*}
(\p_r^2)^{\vee} &=& \p_r^2\\
(\invr\, \p_r)^{\vee} &=& \invr\, \p_r\\
(\invrsq\, \Delta^*)^{\vee} &=& \invrsq\, \bZ^*
\end{eqnarray*}
and the signum--pairs of operators $[\p_r^2\, , \,  \p_r^2]$, $[\invr\, \p_r \, , \, \invr\, \p_r]$, $[\invrsq\, \Delta^*\, , \, \invrsq\, \bZ^*]$ and $[\invrsq\, \bZ^*\, , \, \invrsq\, \Delta^*]$ follow, inducing the action of the operators $\p_r^2$, $\invr\, \p_r$, $\invrsq\, \bZ^*$ and $\invrsq\, \Delta^*$ on a signumdistribution.


\newpage
\section{Spherical operators}
\label{spheroperators}


We say that an operator involving the spherical co-ordinates is {\em spherical} when it is not cartesian. Clearly the multiplication operators $r$ and $\uom$ are spherical operators, as are the differential operators $\p_r$ and $\p_{\uom}$.
The concepts of signumdistribution, signum--pair of operators and cross--pair of operators allow for a definition of the action of spherical operators on (signum)distributions.

\begin{definition}
\label{multiplicomega}
The product of a distribution $T$ by the function $\uom$ is the signumdistribution $T^{\vee}$ associated to $T$, and it holds for all test functions $\uom\, \varphi \in \Omega(\mR^m; \mR^m)$ that
$$
\langle \ \uom\, T \ , \ \uom\, \varphi \ \rangle = \langle \ T^{\vee} \ , \ \uom\, \varphi \ \rangle = -\, \langle \  T \  ,  \varphi \ \rangle
$$
Similarly the product of the signumdistribution $^{s}U$ by the function $(-\, \uom)$ is its associated distribution $ ^{s}U^{\wedge}$.
\end{definition}

\begin{definition}
\label{multiplicomegaj}
The product of a scalar distribution $T^{scal}$ by the function $\omega_j, j=1,\ldots,m$ is the signumdistribution $\omega_j\ T^{scal}$ given by the uniquely determined expression
$$
\omega_j\, T^{scal} = \{ \uom\, T^{scal}  \}_j
$$
Similarly the product of the scalar signumdistribution $^{s}U^{scal}$ by the function $\omega_j$ is the distribution $\omega_j\, ^{s}U^{scal}$ given by
$$
\omega_j\, ^{s}U^{scal} = \{ \uom\, ^{s}U^{scal}  \}_j
$$
\end{definition}

\begin{remark}
{\rm
For a general Clifford algebra valued distribution or signumdistribution, the action of the multiplicative operator $\omega_j$ is defined through linearity with respect to  the scalar components.
}
\end{remark}

\begin{definition}
\label{multiplicr}
The product of a scalar distribution $T$ by the function $r$ is the signumdistribution $r\,T = (-\, \ux\ T)^{\vee}$ given for all test functions $\uom\, \varphi \in \Omega(\mR^m; \mR^m)$ by
$$
\langle \ r\, T \ , \ \uom\, \varphi \ \rangle =  \langle \  \ux\, T \  , \ \varphi \ \rangle =  \langle \  T \  , \  \ux\, \varphi \ \rangle
$$
according to (the boldface part of) the commutative diagram
$$
\begin{array}{ccccccccc}
  &&    & {}_{-\ux}   &                             &  &&&\\
  && {\bf T} & \longrightarrow  & {\bf -\, \ux\, T} & &&& \\[1mm]
  &&\hspace{-8mm}{}_{- \uom} &\hspace{2mm}{} \hspace{11mm}{}_{- r}&\hspace{6mm}{}_{-\uom}&&&&\\[-2mm]
  && \uparrow & \diagdown \hspace{-1.3mm} \nearrow & \uparrow &&&&\\[-1.1mm]
  && \downarrow & \diagup \hspace{-1.3mm} \searrow & \downarrow &&&&\\[-2mm]
  &&\hspace{-4.5mm}{}^{\uom}&\hspace{2mm}{} \hspace{10mm}{}^{r}&\hspace{5mm}{}^{\uom}&&&&\\[1mm]
  && T^{\vee} = \uom\, T & \longrightarrow & {\bf r\, T} &&&&\\[-1mm]
  &&               & {}_{-\ux}  &           &&&&
 \end{array}
$$
involving the signum--pair of operators $(\ux \, , \, \ux)$, which induces the product of the signum--distribution $\uom\, T$ by $\ux$ to be $\ux (\uom\, T) = -\, r\, T$.
\end{definition}

\begin{definition}
\label{radderiv}
The derivative with respect to the radial distance $r$ of a scalar distribution $T$ is the equivalent class of signumdistributions 
$$
\left[ \p_r\,T \right] = \left[ -\, \uom\, \p_r\, T \right]^{\vee} = \left[\frac{1}{\ux} \, \mE\, T \right]^{\vee} = \left( \uS + \underline{c}\, \delta(\ux) \right)^{\vee} =  \uom\, \uS + \uom\, \delta(\ux)\, \underline{c}
$$ 
for any vector distribution $\uS$ satisfying $\ux \, \uS =  \mE\, T$, according to (the boldface part of) the commutative diagram
$$
\begin{array}{ccccccccc}
  &&    & {}_{-\uom \p_r}   &                             &&&&\\
  && {\bf T} & \longrightarrow  & {\bf  \left[\frac{1}{\ux} \, \mE\, T \right] }&&&& \\[1mm]
  &&\hspace{-8mm}{}_{- \uom} &\hspace{2mm}{} \hspace{11mm}{}_{-\p_r}&\hspace{6mm}{}_{-\uom}&&&&\\[-2mm]
  && \uparrow & \diagdown \hspace{-1.3mm} \nearrow & \uparrow &&&&\\[-1.1mm]
  && \downarrow & \diagup \hspace{-1.3mm} \searrow & \downarrow &&&&\\[-2mm]
  &&\hspace{-4.5mm}{}^{\uom}&\hspace{2mm}{} \hspace{10mm}{}^{\p_r}&\hspace{5mm}{}^{\uom}&&&&\\[1mm]
  && T^{\vee} = \uom\, T & \longrightarrow & {\bf \left[ \p_r T \right] } &&&&\\[-1mm]
  &&               & {}_{-\uom \p_r}  &           &&&&
 \end{array}
$$
involving the signum--pair of operators $[\uom\, \p_r \, , \, \uom\, \p_r]$, which induces the action of $(\uom\, \p_r)$ on the signumdistribution $\uom\,  T$ to be be $(\uom\, \p_r)\, (\uom\, T)  = [-\, \p_r\, T]$.
\end{definition}

\begin{definition}
\label{angderiv}
The angular $\p_{\uom}$--derivative of a distribution $T$ is the unique signumdistribution $\p_{\uom}\,T = (\Gamma\, T)^{\vee}$ given for all test functions $\uom\, \varphi \in \Omega(\mR^m; \mR^m)$ by
$$
\langle \  \uom\, \varphi\ , \  \p_{\uom}\, T \ \rangle =  \langle \   \varphi  \  , \  \uom\, \p_{\uom} \, T  \ \rangle = \langle \  \varphi  \  ,  \  -\, \Gamma \, T \ \rangle
$$
according to (the boldface part of) the commutative diagram
$$
\begin{array}{ccccccccc}
  &&    & {}_{- \uom \p_{\uom}}   &                             &  &&&\\
  && {\bf T} & \longrightarrow  & {\bf \Gamma\, T} & &&& \\[1mm]
  &&\hspace{-8mm}{}_{- \uom} &\hspace{2mm}{} \hspace{13mm}{}_{\uom \p_{\uom} \uom}&\hspace{6mm}{}_{-\uom}&&&&\\[-2mm]
  && \uparrow & \diagdown \hspace{-1.3mm} \nearrow & \uparrow &&&&\\[-1.1mm]
  && \downarrow & \diagup \hspace{-1.3mm} \searrow & \downarrow &&&&\\[-2mm]
  &&\hspace{-4.5mm}{}^{\uom} & \hspace{2mm}{} \hspace{10mm}{}^{\p_{\uom}}&\hspace{5mm}{}^{\uom}&&&&\\[1mm]
  && T^{\vee} = \uom\, T & \longrightarrow & {\bf \p_{\uom}\, T}  &&&&\\[-1mm]
  &&               & {}_{- \p_{\uom} \uom}  &           &&&&
 \end{array}
$$
involving  the signumpair of operators $(\,  \uom \p_{\uom}  \, , \, \p_{\uom} \uom \,)$ or $(\Gamma \, , -\, \p_{\uom}\, \uom) = (\Gamma,  (m-1){\bf 1} -\, \Gamma)$, which induces the action of the operator $\p_{\uom}\, \uom$ on the signumdistribution $\uom\, T$ to be $\p_{\uom}\, \uom\, (\uom\, T) = -\, \p_{\uom}\, T$.
\end{definition}

\begin{definition}
\label{angderivj}
The angular $\p_{\omega_j}$--derivative of a scalar distribution $T^{scal}$ is the unique signumdistribution given by
$$
\p_{\omega_j} T^{scal} = \{ \p_{\uom} \, T^{scal} \}_j \quad , \quad j=1,\ldots,m
$$
\end{definition}

\begin{remark}
{\rm\hfill\\
(i) An alternative expression for $\p_{\omega_j} T^{scal}$ is:
$$
\p_{\omega_j} T^{scal} =  r\, \p_{x_j}\, T^{scal} - \omega_j\, \mE\, T^{scal}
$$
Indeed, it follows from
$$
-\, \Gamma\, T^{scal} = \ux\, \dirac\, T^{scal} + \mE\, T^{scal}
$$
that 
$$
-\, \p_{\uom}\, T^{scal} = -\, \uom\, \ux\, \dirac\, T^{scal} + \uom\, \mE\, T^{scal} = -\, r \dirac\, T^{scal} + \uom\, \mE\, T^{scal}
$$
whence the desired formula for each of the components.\\
(ii) For a general Clifford algebra valued distribution, the action of the spherical derivative operator $\p_{\omega_j}$ is defined through linearity with respect to  the scalar components.
}
\end{remark}

\begin{definition}
\label{T/r}
The quotient of a scalar distribution $T$ by the radial distance $r$ is the equivalence class of signumdistributions
$$
\left[ \frac{1}{r}\, T \right] = \uom\, \left[ \frac{1}{\ux}\, T  \right] = \uom\, (\underline{S} + \delta(\ux)\, \underline{c} ) = \uom\, \underline{S} + \uom\, \delta(\ux)\, \underline{c} = \underline{S}^{\vee} + \delta(\ux)^{\vee}\, \underline{c}
$$
for any vector-valued distribution $\underline{S}$ for which $\ux\, \underline{S} = T$,
according to (the boldface part of) the commutative diagram
$$
\begin{array}{ccccccccc}
  &&    & {}_{\frac{1}{\ux}}   &                             &  &&&\\
  && {\bf T} & \longrightarrow  & {\bf \left[ \frac{1}{\ux}\, T \right]} & &&& \\[1mm]
  &&\hspace{-8mm}{}_{- \uom} &\hspace{2mm}{} \hspace{11mm}{}_{- \frac{1}{r}}&\hspace{6mm}{}_{-\uom}&&&&\\[-2mm]
  && \uparrow & \diagdown \hspace{-1.3mm} \nearrow & \uparrow &&&&\\[-1.1mm]
  && \downarrow & \diagup \hspace{-1.3mm} \searrow & \downarrow &&&&\\[-2mm]
  &&\hspace{-4.5mm}{}^{\uom}&\hspace{2mm}{} \hspace{10mm}{}^{\frac{1}{r}}&\hspace{5mm}{}^{\uom}&&&&\\[1mm]
  && T^{\vee} = \uom\, T & \longrightarrow & {\bf \left[ \frac{1}{r}\, T \right]} &&&&\\[-1mm]
  &&               & {}_{\frac{1}{\ux}}  &           &&&&
 \end{array}
$$
involving the signum--pair of operators $[\invux\, , \, \invux]$, which induces the quotient of the signumdistribution $\uom\, T$ by $\ux$ to be $\invux\, (\uom\, T) = [\invr\, T]$.
\end{definition}

\begin{remark}
{\rm\hfill\\
Once the action of the spherical derivative operators $\p_{\omega_j}, j=1,\ldots,m$ is defined, see Definition \ref{angderivj}, we are able to define the action of the corresponding cartesian operators $-\, \uom\,  \p_{\omega_j}, j=1,\ldots,m$, 
through the commutative diagram
$$
\begin{array}{ccccccccc}
  &&    & {}_{-\, \uom\, \p_{\omega_j}}   &                             &  &&&\\
  &&  T & \longrightarrow  &  -\, \uom\, \p_{\omega_j} T & &&& \\[1mm]
  &&\hspace{-8mm}{}_{- \uom} &\hspace{2mm}{} \hspace{15mm}{}_{ \uom\, \p_{\omega_j\, \uom} }&\hspace{6mm}{}_{-\uom}&&&&\\[-2mm]
  && \uparrow & \diagdown \hspace{-1.3mm} \nearrow & \uparrow &&&&\\[-1.1mm]
  && \downarrow & \diagup \hspace{-1.3mm} \searrow & \downarrow &&&&\\[-2mm]
  &&\hspace{-4.5mm}{}^{\uom} & \hspace{2mm}{} \hspace{13mm}{}^{\p_{\omega_j}}&\hspace{5mm}{}^{\uom}&&&&\\[1mm]
  &&  \uom\, T & \longrightarrow &  \p_{\omega_j}\, T  &&&&\\[-1mm]
  &&               & {}_{-\,  \p_{\omega_j}\, \uom}  &           &&&&
 \end{array}
$$
In particular for a scalar distribution $T^{scal}$ it holds that
$$
-\, \uom\, \p_{\omega_j}\, T^{scal} = -\, \ux\, \p_{x_j}\, T^{scal} + \omega_j\, \uom\, \mE\, T^{scal} \quad , \quad j=1,\ldots,m
$$
}
\end{remark}

\begin{example}
{\rm

As $\delta(\ux)$ is a radial distribution we expect  $\p_{\omega_j}\, \delta(\ux)$, and thus also $-\, \uom\, \p_{\omega_j}\, \delta(\ux) $,  to be zero. And indeed it holds that
\begin{align*}
\p_{\omega_j}\, \delta(\ux) &= r\, \p_{x_j}\, \delta(\ux) - \omega_j\, \mE\, \delta(\ux)\\
& = r\, \p_{x_j}\, \delta(\ux) + m\, \omega_j\, \delta(\ux)\\
&= \{ r\, \dirac\, \delta(\ux) + m\, \uom\, \delta(\ux)  \}_j\\
&= \{ \uom\, \mE\, \delta(\ux)   + m\, \uom\, \delta(\ux) \}_j\\
&= 0
\end{align*}
}
\end{example}


\newpage
\section{Action uniqueness of some operators}
\label{uniqueness}


In the preceding sections we encountered operators acting on (signum)distributions with a uniquely determined result and other ones whose actions are not uniquely determined but lead to equivalence classes of (signum)distributions instead. Nevertheless in \cite{ffb} sufficient conditions were found guaranteeing the uniqueness of the latter operators' actions, involving homogeneous, radial and signum--radial (signum)distributions.

\begin{definition}
\label{defsignumradial}
(i) A distribution $T$ or a signumdistibution $^{s}U$ respectively, is said to be {\em radial} if it is SO$(m)$--invariant and so only depends on $r = |\ux|$: $T(\ux) = T(r)$ or $^{s}U(\ux) =\, ^{s}U(r)$ respectively.\\
(ii) A distribution, signumdistribution respectively, is said to be signum--radial if its associated signumdistribution, distribution respectively, is radial.
\end{definition}
\noindent
Let us state these sufficient conditions.

\begin{itemize}
\item[(i)] If the distribution $T^{rad}$ is radial then the following actions are uniquely determined:
$$
(\uom\, \p_r)\, T^{rad} \qquad \left(-\, \invr\, \p_{\uom} + (m-1)\, \invr\, \uom \right)\, T^{rad} \qquad  \invux\, T^{rad} \qquad \uD\, T^{rad} \qquad \p_r\, T^{rad} \qquad \invr\, T^{rad}
$$
and the actions of the corresponding signum--partner operators on the signum--radial signumdistribution $^{s}U^{srad}$, viz.
$$
(\uom\, \p_r)\, ^{s}U^{srad} \qquad \qquad \left(\invr\, \p_{\uom} \right)\, ^{s}U^{srad} \qquad \qquad \qquad  \invux\, ^{s}U^{srad} \qquad \dirac\, ^{s}U^{srad} \qquad \p_r\, ^{s}U^{srad} \qquad \invr\, ^{s}U^{srad}
$$
are uniquely determined too.

\item[(ii)] If the distribution $T^{(k)}$ is homogeneous with homogeneity degree $k \neq -m+1$ then the following actions are uniquely determined:
$$
(\uom\, \p_r)\, T^{(k)} \qquad \left(\invr\, \p_{\uom} \right)\, T^{(k)} \qquad \left(-\, \invr\, \p_{\uom} + (m-1)\, \invr\, \uom \right)\, T^{(k)} \qquad  \invux\, T^{(k)} \qquad \uD\, T^{(k)} \qquad \p_r\, T^{(k)} \qquad \invr\, T^{(k)}
$$
and the actions of the corresponding signum--partner operators on the homogeneous signumdistribution $^{s}U^{(k)}$, viz.
$$
(\uom\, \p_r)\, ^{s}U^{(k)} \qquad \left(-\, \invr\, \p_{\uom} + (m-1)\, \invr\, \uom \right)\, ^{s}U^{(k)} \qquad \left(\invr\, \p_{\uom} \right)\, ^{s}U^{(k)} \qquad \invux\, ^{s}U^{(k)} \qquad \dirac\, ^{s}U^{(k)} \qquad \p_r\, ^{s}U^{(k)} \qquad \invr\, ^{s}U^{(k)}
$$
are uniquely determined too.

\item[(iii)] If the distribution $T^{(k)}$ is homogeneous with homogeneity degree $k \neq -m+1, -m+2$ then the following actions are uniquely determined:
$$
\p_r^{2}\, T^{(k)} \qquad \invr\, \p_{r} \, T^{(k)} \qquad \invrsq\, \Delta^* \, T^{(k)} \qquad  \invrsq\, \bZ^* \, T^{(k)} \qquad \bZ\, T^{(k)} \qquad \invrsq\, T^{(k)}
$$
and the actions of the corresponding signum--partner operators on the homogeneous signumdistribution $^{s}U^{(k)}$, viz.
$$
\p_r^{2}\, ^{s}U^{(k)} \qquad \invr\, \p_{r} \, ^{s}U^{(k)} \qquad \invrsq\, \bZ^* \, ^{s}U^{(k)} \qquad \invrsq\, \Delta^* \,  ^{s}U^{(k)} \qquad \Delta\, ^{s}U^{(k)} \qquad \invrsq\, ^{s}U^{(k)}
$$
are uniquely determined too.

\item[(iv)] The conclusions contained in (iii) remain valid if the distribution, signumdistribution respectively, under consideration is both radial, signum--radial respectively, and homogeneous with homogeneity degree $k \neq -m+2$.

\item[(v)] The conclusions contained in (iii) remain valid if the distribution, signumdistribution respectively, under consideration is both signum--radial, radial respectively, and homogeneous with homogeneity degree $k \neq -m+1$.

\end{itemize}


\newpage
\section{Cartesian derivatives of signumdistributions}
\label{signumcartderiv}


Recalling the signum--pairs of operators $(\dirac\, , \underline{D})$ and $[\underline{D}\, , \dirac]$, we expect the cartesian derivatives of (signum)distributions to show up in signum--pairs of operators of the form $(\p_{x_j}\, , d_j)$ and $[d_j\, , \p_{x_j}]$, $j=1,\ldots,m$. These operators $d_j, j=1,\ldots,m$ are important in the following sense. When defining a derivative of a distribution, which is merely a sophisticated integration by parts, the differential operator, say $\p_{x_j}$ shifts to the test function. It was shown in \cite{ffb} that when computing  the derivative $\p_{x_j}$ of a signumdistribution the derivative shifts to the test function but then in the form of the operator $d_j$.\\

\begin{definition}
For $j=1,\ldots,m$ we define the operator $d_j$ to be the signum--partner of the cartesian derivative $\p_{x_j}$ and it holds that
\begin{eqnarray*}
d_j &=& \uom\, \p_{x_j}\, (-\uom)\\
&=& -\, \invr\, \uom\, e_j -\, \invr\, \omega_j + \p_{x_j}  = \invr\, e_j\, \uom + \invr\, \omega_j + \p_{x_j}\\
&=& -\, \invrsq\, \ux\, e_j - \invrsq\, x_j + \p_{x_j} = \invrsq\, e_j\, \ux + \invrsq\, x_j + \p_{x_j}\\
&=&  \invux\, e_j + (-\, \frac{x_j}{r^2}) + \p_{x_j} = -\, e_j\, \invux -  (-\, \frac{x_j}{r^2}) + \p_{x_j}\\
\end{eqnarray*}
\end{definition}
\noindent
Clearly the operator $d_j\, (j=1,\ldots,m)$ is cartesian. It is the sum of the scalar part $\p_{x_j}$ and the bivector part $\invrsq\, e_j\, \ux + \invrsq\, x_j  = \invrsq\, (e_je_1x_1 + \ldots + e_je_{j-1}x_{j-1} + e_je_{j+1}x_{j+1}+ \ldots + e_je_mx_m )$. This bivector part is a combination of the operator $\invux$ and its components $(-\, \frac{j}{r^2}, j=1,\ldots,m)$.\\
As the operator $d_j$  is a well--defined, but not uniquely defined, operator acting on distributions, also the signum--pairs of operators $[d_j\, , \p_{x_j}]$, $j=1,\ldots,m$ hold, enabling the definition of the cartesian derivative $\p_{x_j}$ of a signumdistribution as an equivalence class of signumdistributions, as is seen by the following commutative diagram:
$$
\begin{array}{ccccccccc}
  &&    & {}_{d_j}   &                             &&&&\\
  &&  ^{s}U^{\wedge} = -\uom\,  ^{s}U & \longrightarrow  & [  d_j\, (-\uom\, ^{s}U) ] &&&& \\[1mm]
  &&\hspace{-8mm}{}_{- \uom}&\hspace{2mm}{} \hspace{11mm}{}_{}&\hspace{6mm}{}_{-\uom}&&&&\\[-2mm]
  && \uparrow & \phantom{\diagdown} \hspace{-1.3mm} \phantom{\nearrow} & \uparrow &&&&\\[-1.1mm]
  && \downarrow & \phantom{\diagup} \hspace{-1.3mm} \phantom{\searrow} & \downarrow &&&&\\[-2mm]
  &&\hspace{-4.5mm}{}^{\uom}&\hspace{2mm}{} \hspace{12mm}{}^{}&\hspace{5mm}{}^{\uom}&&&&\\[1mm]
  &&   ^{s}U & \longrightarrow & \left[ \p_{x_j}\, ^{s}U \right] &&&&\\[-1mm]
  &&               & {}_{\p_{x_j}}  &           &&&&
 \end{array}
$$

\begin{remark}
{\rm
As one may expect there exists a relationship between the operators $d_j, j=1,\ldots,m$, the Dirac operator $\dirac$ and the operator $\uD$. A straightforward calculation shows that
\begin{eqnarray*}
\dirac &=& \phantom{-}\, \sum_{j=1}^m\, e_j\, d_j + (m-1)\, \invrsq\, \ux =  \phantom{-}\, \sum_{j=1}^m\, e_j\, d_j - (m-1)\, \invux \\
\uD &=& -\, \sum_{j=1}^m\, e_j\, d_j  + 2\, \invrsq\, \ux\, \mE =  -\, \sum_{j=1}^m\, e_j\, d_j  - 2\, \invux\, \mE
\end{eqnarray*}
}
\end{remark}

\begin{table}
\begin{center}
\renewcommand{\arraystretch}{2.4}
\begin{tabular}{|c|c|}
\hline
\multicolumn{2}{|c|}{$(\, \ux \, , \, \ux\, )$}   \\ \hline
\multicolumn{2}{|c|}{$(\, r^2\, , \, r^2\, )$}    \\ \hline 
\multicolumn{2}{|c|}{$(\, \mE\, , \, \mE\, )$}  \\ \hline
$(\, \Gamma\, , \, -\, \p_{\uom}\, \uom\, )$ & $(\, -\, \p_{\uom}\, \uom\, , \Gamma\, ) $ \\ \hline
$(\, \Gamma^2\, , \, \Gamma^2 - 2(m-1)\Gamma + (m-1)^2\, )$ & $(\, \Gamma^2 - 2(m-1)\Gamma + (m-1)^2\, , \Gamma^2\, )$ \\ \hline
$(\, \dirac\, , \, \uD\, )$ & $[\, \uD\, , \, \dirac\, ]$ \\ \hline
\multicolumn{2}{|c|}{$[\, \uom\, \p_r\, , \, \uom\, \p_r\, ]$} \\ \hline
$[\, \invr\, \p_{\uom}\, , \, -\, \invr\, \p_{\uom} + (m-1)\, \invr\, \uom\, ]$ & $[\, -\, \invr\, \p_{\uom} + (m-1)\, \invr\, \uom\, , \, \invr\, \p_{\uom}\, ]$ \\ \hline
\multicolumn{2}{|c|}{$(\, \p_{\uom}^2\, , \, \p_{\uom}^2\, )$} \\ \hline
$(\, \Delta^*\, , \, {\bf Z}^*\, )$ & $(\, {\bf Z}^*\, , \, \Delta^*\, )$ \\ \hline
$(\, \Delta\, , \, {\bf Z}\, )$ & $[\, {\bf Z}\, , \, \Delta\, ]$ \\ \hline
\multicolumn{2}{|c|}{$[\, \p_r^2\, , \, \p_r^2\, ]$} \\ \hline
\multicolumn{2}{|c|}{$[\, \invux\, , \, \invux\, ]$} \\ \hline
\multicolumn{2}{|c|}{$[\, \invr\, \p_r\, , \, \invr\, \p_r\, ]$}  \\ \hline
\multicolumn{2}{|c|}{$[\, \invrsq\, , \, \invrsq\, ]$} \\ \hline
$(\, \p_{x_j}\, , \, d_j\,)$ & $[\, d_j\, , \, \p_{x_j}\, ]$ \\ \hline
\end{tabular}
\caption{Signum--pairs of operators}\label{tablesignumpairs}
\end{center}
\end{table}

\begin{table}
\begin{center}
\renewcommand{\arraystretch}{2.4}
\begin{tabular}{|c|c|}
\hline
\multicolumn{2}{|c|}{$(\, \uom \, , \, \uom\, )$} \\ \hline
$(\, r\, , \, -\, r\, )$  &  $(-\, r\, , \,  r\, )$ \\ \hline 
$[\, \p_r\, , \, -\, \p_r\, ]$ & $[-\, \p_r\, , \,  \p_r\, ]$ \\ \hline
$(\, \p_{\uom}\, , \, \uom\, \p_{\uom}\, \uom\, )$ & $(\, \uom\, \p_{\uom}\, \uom\, , \p_{\uom}\, )$ \\ \hline
$[\, \invr\, , \, -\, \invr\, ]$ & $[-\, \invr\, , \,  \invr\, ]$ \\ \hline
$[\, \invr\, \p_{\uom}\, \uom\, , \, -\, \invr\, \uom\, \p_{\uom}\, ]$ & $[-\, \invr\, \p_{\uom}\, \uom\, , \,  \invr\, \uom\, \p_{\uom}\, ]$ \\ \hline
\end{tabular}
\caption{Cross--pairs of operators}\label{tablecrosspairs}
\end{center}
\end{table}


\newpage
\section{Actions on $T_\la$ and $U_\la$}
\label{actions}


The results of the preceding sections will now be applied to the distributions $T_\la$ and $U_\la$.
Let us recall their definitions: for $\la \in \mC$ it holds that
$$\langle  \ T_\lambda , \varphi(\ux)  \ \rangle := a_m\, \langle \  {\rm Fp}\, r^{\lambda+m-1}_+  , \Sigma^0[\varphi](r)  \  \rangle_r$$
and
$$\langle  \  U_\lambda , \varphi(\ux) \ \rangle := a_m\, \langle \  {\rm Fp}\, r^{\lambda+m-1}_+  , \Sigma^1[\varphi](r)  \  \rangle_r$$\\
An alternative, and handy, notation could be $T_\la = {\rm Fp}\, r^\la$ and $U_\la =  \uom\, {\rm Fp}\, r^\la$.\\

\noindent
Let us comment on these definitions; for more details we refer to the series of papers cited at the beginning of Section \ref{intro}.\\

\noindent
A priori the distributions $T_\la$ and $U_\la$ are not defined at the simple poles of the one--dimensional distribution Fp$ \; r^{\la+m-1}_+$ on the $r$--axis, viz., $\lambda = -m-n+1$, $n \in \mN$. To cope with these singularities, the distributions Fp$ \; r^{-n}_+$, $n \in \mN$ are interpreted as the so--called {\em monomial pseudofunctions}  (see \cite{fb1} and \cite{fb2}).\\

\noindent
The distributions $T_\lambda$ are standard scalar  distributions well--known in harmonic analysis. They are radial and homogeneous of degree $\la$. As meromorphic functions of $\lambda \in \mC$ they show genuine simple poles at $\lambda = -m, -m-2, -m-4, \ldots$.  This is due to the fact that the singular points $\lambda=-m-2\ell-1$, $\ell=0,1,2,\ldots$ are removable, since the spherical mean $\Sigma^{(0)}[\phi]$ has its odd order derivatives vanishing at $r=0$. So we can define
$$
\langle T_{-m-2\ell-1},\phi \rangle = \lim_{\mu \rightarrow -2\ell-2} a_m \langle {\rm Fp} \; r^\mu_+, \Sigma^{(0)}[\phi] \rangle
$$
but, remarkably, this limit is precisely $a_m \langle {\rm Fp} \; r^{-2\ell-2}_+,\Sigma^{(0)}[\phi] \rangle$, with
Fp$ \; r^{-2\ell-2}_+$ the monomial pseudofunction.
The most important distribution is this family is $T_{-m+2} = \frac{1}{r^{m-2}}$, which is, up to a constant, the fundamental solution of the Laplace operator $\Delta$. Also note the special cases  $T_0 = 1$, $T_{2\ell} = r^{2\ell} =  (-1)^\ell\, \ux^{2\ell}$, and $T_{2\ell+1} = r^{2\ell+1}, \ell=0,1,2,\ldots$.\\

\noindent
The distributions $U_\lambda$ are typical Clifford analysis constructs. They are homogeneous of degree $\la$. As vector--valued meromorphic functions of $\lambda \in \mC$ they show genuine simple poles at $\lambda = -m-1, -m-3, -m-5, \ldots$. This is due to the fact that the singular points $\lambda=-m-2\ell$, $\ell=0,1,2,\ldots$ are removable, since the spherical mean $\Sigma^{(1)}[\phi]$ has its even order derivatives vanishing at $r=0$. So we can define
$$
\langle U_{-m-2\ell},\phi \rangle = \lim_{\mu \rightarrow -2\ell-1} a_m \langle {\rm Fp} \; r^\mu_+, \Sigma^{(1)}[\phi] \rangle
$$
but this limit is precisely $a_m \langle {\rm Fp} \; r^{-2\ell-1}_+,\Sigma^{(1)}[\phi] \rangle$, with
Fp$ \; r^{-2\ell-1}_+$ the monomial pseudofunction.
The most important distribution in this family is 
$U_{-m+1} = \frac{\uom}{r^{m-1}} = \frac{\ux}{r^m}$ which is, up to a constant, the fundamental solution of the Dirac operator $\dirac$. Also note the special cases $U_0 = \uom$, $U_{2\ell} = \uom\, r^{2\ell}$ and $U_{2\ell+1} = \uom\, r^{2\ell+1} =  (-1)^\ell\, \ux^{2\ell+1}, \ell=0,1,2,\ldots$.\\

\noindent
It is important to note that although the distributions $T_\la$ and $U_\la$ are also defined in their respective singularities, these exceptional values do {\em not} turn these distributions into entire functions of the parameter $\lambda \in \mC$.\\

\noindent
When restricted to the half--plane $\mR e\, \lambda > -\, m$ the distributions $T_\lambda$ and $U_\lambda$ are regular, i.e. they are locally integrable functions.
From \cite{ffb, bsv} we know that a locally integrable function can be seen as a signumdistribution as well. This inspires the definition of the following two families of signumdistributions:

$$\langle  \ ^{s}T_\lambda , \uom\, \varphi(\ux)  \ \rangle := \phantom{-}\, a_m\, \langle \  {\rm Fp}\, r^{\lambda+m-1}_+  , \Sigma^1[\varphi](r)  \  \rangle_r$$

$$\langle  \ ^{s}U_\lambda , \uom\, \varphi(\ux) \ \rangle := -\, a_m\, \langle \  {\rm Fp}\, r^{\lambda+m-1}_+  , \Sigma^0[\varphi](r)  \  \rangle_r$$

\noindent
It is clear that
$$T_\lambda^{\, \vee} = \phantom{-}\, ^{s}U_\lambda \quad , \quad ^{s}U_\lambda^{\wedge} = \phantom{-}\, T_\lambda$$
and
$$U_\lambda^{\, \vee} = -\, ^{s}T_\lambda \quad , \quad ^{s}T_\lambda^{\wedge} \, = -\, U_\lambda$$

\noindent
In this way $^{s}T_\lambda$ inherits the simple poles of $U_\lambda$, viz., $\lambda = -m-1, -m-3, -m-5, \ldots$, while $^{s}U_\lambda$ inherits the simple poles of 
$T_\lambda$, viz., $\lambda = -m, -m-2, -m-4, \ldots$.\\
As the distributions $T_\la$ are radial, by Definition \ref{defsignumradial} the signumdistributions  $^{s}U_\la$ are  signum--radial. As the signumdistributions $^{s}T_\la$ are radial, by the same definition the distributions $U_\la$ are signum--radial.\\

In the following subsections we will systematically compute the actions on $T_\la$, $U_\la$, $^{s}T_\la$ and $^{s}U_\la$ of all operators introduced in the precedings sections, paying attention to the uniqueness of the expressions obtained.


\newpage
\subsection{The operators $\ux$, $r$ and $r^2$}


The multiplication operator $\ux$ is a cartesian operator whose actions are uniquely determined, quite naturally.  It holds for all $\la \in \mC$ that
\begin{equation}
\label{xT}
\boxed{\ux\, T_\lambda = U_{\lambda+1}} \quad , \quad \boxed{\ux\, U_\lambda = -\, T_{\lambda+1}} 
\end{equation}

\noindent
Based on the commutative diagram
$$
\begin{array}{ccccccccc}
  &&    & {}_{-\ux}   &                             &  &&&\\
  &&  T_\lambda & \longrightarrow  &  -\, U_{\lambda+1} & &&& \\[1mm]
  &&\hspace{-8mm}{}_{- \uom} &\hspace{2mm}{} \hspace{11mm}{}_{- r}&\hspace{6mm}{}_{-\uom}&&&&\\[-2mm]
  && \uparrow & \diagdown \hspace{-1.3mm} \nearrow & \uparrow &&&&\\[-1.1mm]
  && \downarrow & \diagup \hspace{-1.3mm} \searrow & \downarrow &&&&\\[-2mm]
  &&\hspace{-4.5mm}{}^{\uom}&\hspace{2mm}{} \hspace{10mm}{}^{r}&\hspace{5mm}{}^{\uom}&&&&\\[1mm]
  && ^{s}U_\lambda & \longrightarrow &  ^{s}T_{\lambda+1} &&&&\\[-1mm]
  &&               & {}_{-\ux}  &           &&&&
 \end{array}
$$
we find the additional formulae
$$
\boxed{r\, T_\lambda =  \, ^{s}T_{\lambda+1}} \quad {\rm and} \quad \boxed{r\, ^{s}U_\lambda =  U_{\lambda+1}} \quad , \quad \la \in \mC
$$
and
$$
\boxed{\ux\, ^{s}U_\lambda = -\, ^{s}T_{\lambda+1}} \quad , \quad \la \in \mC
$$
In a similar way, based on the commutative diagram
$$
\begin{array}{ccccccccc}
  &&    & {}_{\ux}   &                             &  &&&\\
  &&  -\, U_\lambda & \longrightarrow  &  T_{\lambda+1} & &&& \\[1mm]
  &&\hspace{-8mm}{}_{- \uom} &\hspace{2mm}{} \hspace{11mm}{}_{r}&\hspace{6mm}{}_{-\uom}&&&&\\[-2mm]
  && \uparrow & \diagdown \hspace{-1.3mm} \nearrow & \uparrow &&&&\\[-1.1mm]
  && \downarrow & \diagup \hspace{-1.3mm} \searrow & \downarrow &&&&\\[-2mm]
  &&\hspace{-4.5mm}{}^{\uom}&\hspace{2mm}{} \hspace{10mm}{}^{- r}&\hspace{5mm}{}^{\uom}&&&&\\[1mm]
  && ^{s}T_\lambda & \longrightarrow &  ^{s}U_{\lambda+1} &&&&\\[-1mm]
  &&               & {}_{\ux}  &           &&&&
 \end{array}
$$
we obtain the additional formulae
$$
\boxed{r\, ^{s}T_\lambda =  T_{\lambda+1}} \quad {\rm and} \quad \boxed{r\, U_\lambda =  \, ^{s}U_{\lambda+1}} \quad , \quad \la \in \mC
$$
and also
$$
\boxed{\ux\, ^{s}T_\lambda = \, ^{s}U_{\lambda+1}} \quad , \quad \la \in \mC
$$
Iteration of the multiplication operator $\ux$ results into
$$
\boxed{r^2\ T_\la = T_{\la+2}} \quad , \quad \boxed{r^2\, U_\la = U_{\la+2}} \quad , \quad \la \in \mC
$$
and
$$
\boxed{r^2\ ^{s}T_\la =\, ^{s}T_{\la+2} }\quad , \quad \boxed{r^2\, ^{s}U_\la =\, ^{s}U_{\la+2}}\quad , \quad \la \in \mC
$$\\

The natural powers of $\ux$ are cartesian operators too. We find, for all $\la \in \mC$,
\begin{align*}
(-1)^{\ell}\, \ux^{2\ell}\, T_\la &= r^{2\ell}\, T_\la = T_{\la + 2\ell}\\[1mm]
(-1)^{\ell}\, \ux^{2\ell+1}\, T_\la &= (\uom\, r^{2\ell+1})\, T_\la = U_{\la + 2\ell+1}\\
\end{align*}
and
\begin{align*}
(-1)^{\ell}\, \ux^{2\ell}\, U_\la &= r^{2\ell}\, U_\la = U_{\la + 2\ell}\\[1mm]
(-1)^{\ell+1}\, \ux^{2\ell+1}\, U_\la &= -\, (\uom\, r^{2\ell+1})\, U_\la = T_{\la + 2\ell+1}\\
\end{align*}
Through the appropriate commuative diagrams we find, for all $\la \in \mC$,
\begin{align*}
 r^{2\ell}\, ^{s}T_\la &=\,  ^{s}T_{\la + 2\ell}\\[1mm]
r^{2\ell}\, ^{s}U_\la &=\,  ^{s}U_{\la + 2\ell}
\end{align*}
and
\begin{align*}
 r^{2\ell+1}\, T_\la &=\,  ^{s}T_{\la + 2\ell+1}\\[1mm]
r^{2\ell+1}\, U_\la &=\,  ^{s}U_{\la + 2\ell+1}
\end{align*}


\newpage
\subsection{The operators $\dirac$ and $\Delta$}


As does the multiplication operator $\ux$, also the Dirac operator $\dirac$ intertwines the $T_\la$ and $U_\la$ distribution families. It clearly is a cartesian operator whose action is uniquely determined.  It holds that
\begin{equation}
\label{diracTlambda}
\boxed{\dirac\, T_\lambda =  \lambda\, U_{\lambda - 1}} \quad , \quad \lambda \neq -m, -m-2, -m-4, \ldots
\end{equation}
and
\begin{equation}
\label{diracUlambda}
\boxed{\dirac\, U_\lambda =  - (\lambda+m-1)\, T_{\lambda - 1}} \quad , \quad \lambda \neq -m+1, -m-1, -m-3, \ldots
\end{equation}
while for $\ell = 0, 1, 2, \ldots$
\begin{equation}
\label{diracTspecial}
\boxed{
\dirac\, T_{-m-2\ell} = -\, (m+2\ell)\, U_{-m-2\ell-1}  + (-1)^{\ell+1}\, \frac{1}{C(m,\ell)}\, a_m\, \dirac^{2\ell+1}\, \delta(\ux)
}
\end{equation}
and
\begin{equation}
\label{diracUspecial}
\boxed{
\dirac\, U_{-m-2\ell+1} = (2\ell)\, T_{-m-2\ell} + (-1)^{\ell-1}\, \frac{m+2\ell}{C(m,\ell)}\, a_m\, \dirac^{2\ell}\, \delta(\ux)
}
\end{equation}
with
$$
C(m,\ell) = 2^{2\ell+1}\, \ell!\, \frac{\Gamma(\frac{m}{2}+\ell+1)}{\Gamma(\frac{m}{2})} =  2^{\ell}\, \ell!\, m(m+2)(m+4)\cdots(m+2\ell)
$$
In particular, for $\ell=0$, it holds that
\begin{equation}
\label{diracTm}
\boxed{
\dirac\, T_{-m} = (-m)\, U_{-m-1} - \frac{1}{m}\, a_m\, \dirac\delta(\ux)
}
\end{equation}
\begin{equation}
\label{diracUm11}
\boxed{
\dirac\, U_{-m+1} = -\, a_m\, \delta(\ux)
}
\end{equation}
Formula (\ref{diracUm11}) expresses the well--known fact that $-\, \frac{1}{m}\, U_{-m+1}$ is indeed the fundamental solution of the Dirac operator $\dirac$.\\
Note also that
$$
\dirac\, T_0 = 0 \quad {\rm and} \quad \dirac\, \ln{r} = U_{-1}
$$\\

\noindent
Through the signum--pair  of operators $(\dirac , \uD)$ the corresponding formulae for the signumdistributions $^{s}T_\la$ and $^{s}U_\la$ are readily obtained:

\begin{equation*}
\label{diracT}
\boxed{\uD\, ^{s}U_\lambda =  -\, \lambda\, ^{s}T_{\lambda - 1}} \quad , \quad \lambda \neq -m, -m-2, -m-4, \ldots
\end{equation*}
and
\begin{equation*}
\label{diracU}
\boxed{\uD\, ^{s}T_\lambda =   (\lambda+m-1)\, ^{s}U_{\lambda - 1}} \quad , \quad \lambda \neq -m+1, -m-1, -m-3, \ldots
\end{equation*}
while for $\ell = 0, 1, 2, \ldots$
\begin{equation*}
\boxed{
\uD\, ^{s}U_{-m-2\ell} =  (m+2\ell)\, ^{s}T_{-m-2\ell-1}  + (-1)^{\ell+1}\, \frac{1}{C(m,\ell)}\, a_m\, \uom\, \dirac^{2\ell+1}\, \delta(\ux)
}
\end{equation*}
and
\begin{equation*}
\boxed{
\uD\, ^{s}T_{-m-2\ell+1} = -\, (2\ell)\, ^{s}U_{-m-2\ell} + (-1)^{\ell}\, \frac{m+2\ell}{C(m,\ell)}\, a_m\, \uom\, \dirac^{2\ell}\, \delta(\ux)
}
\end{equation*}
and in particular, for $\ell=0$,
\begin{equation*}
\boxed{
\uD\, ^{s}U_{-m} = m\, ^{s}T_{-m-1} - \frac{1}{m}\, a_m\, \uom\, \dirac\delta(\ux)
}
\end{equation*}
\begin{equation*}
\boxed{
\uD\, ^{s}T_{-m+1} =  a_m\, \uom\, \delta(\ux)
}
\end{equation*}
\\

Iteration of the Dirac operator $\dirac$ results into formulae for the action of the Laplace operator on distributions. It holds that

\begin{eqnarray}
\label{lapT} \Delta\, T_\lambda &=&  \lambda\,(\lambda+m-2)\, T_{\lambda-2} \quad , \quad \lambda \neq -m+2, -m, -m-2, \ldots\\
\label{lapU} \Delta\, U_\lambda &=&  (\lambda-1)(\lambda+m-1)\, U_{\lambda-2} \quad , \quad \lambda \neq -m+1, -m-1, \ldots
\end{eqnarray}
and
\begin{eqnarray}
\label{lapTspecial}
\Delta\, T_{-m-2\ell} &=&   (m+2\ell)(2\ell+2)\, T_{-m-2\ell-2} + (-1)^\ell\, a_m\, \frac{(m+4\ell+2)(m+2\ell+2)}{C(m,\ell+1)}\, \dirac^{2\ell+2}\, \delta(\ux)\\
\label{lapUspecial}
\Delta\, U_{-m-2\ell+1} &=& (m+2\ell)(2\ell)\, U_{-m-2\ell-1} + (-1)^\ell\, a_m\, \frac{m+4\ell}{C(m,\ell)}\, \dirac^{2\ell+1}\, \delta(\ux)
\end{eqnarray}
and, in particular, for $\ell=0$
\begin{eqnarray}
\label{lapTm}
\Delta\, T_{-m} &=&   2m\, T_{-m-2} + \frac{m+2}{2m}\, a_m\, \dirac^2\, \delta(\ux) \\
\label{lapUm1}
\Delta\, U_{-m+1} &=& a_m\, \dirac \delta(\ux)
\end{eqnarray}
and also
\begin{eqnarray}
\label{lapTm2}
\Delta\, T_{-m+2} &=& -\, (m-2)\, a_m\, \delta(\ux)
\end{eqnarray}
this last formula expressing the fact that $-\, \frac{1}{m-2}\, \frac{1}{a_m}\, T_{-m+2}$ is, indeed, the fundamental solution of the Laplace operator.\\

As an aside notice that continuing the iteration by the Dirac operator $\dirac$ leads to the fundamental solutions of the natural powers of $\dirac$. We find
$$
\dirac^{2\ell}\, E_{2\ell} = \dirac^{2\ell} \left( \frac{1}{2^{\ell-1}\, (\ell-1)!\, (m-2)(m-4)\cdots(m-2\ell)}\, \frac{1}{a_m}\, T_{-m+2\ell}  \right) = \delta(\ux)
$$
and
$$
\dirac^{2\ell+1}\, E_{2\ell+1} = \dirac^{2\ell+1} \left( -\, \frac{1}{2^{\ell}\, (\ell)!\, (m-2)(m-4)\cdots(m-2\ell)}\, \frac{1}{a_m}\, U_{-m+2\ell+1}  \right) = \delta(\ux)
$$
If the dimension $m$ is odd then the above formulae are valid for all natural values of $\ell$. However if the dimension $m$ is even, then these expressions are only valid for $\ell < m/2$; this already becomes clear from the fundamental solution $E_m$ of the operator
$\dirac^m$ which is logarithmic in nature:
$$
\dirac^m\, E_m = \dirac^m \left(  -\, \frac{1}{2^{m-1}\, \pi^{m/2}\, \Gamma(m/2)}\, \ln r \right) = \delta(\ux) \quad , \quad m \ {\rm even}
$$
More generally it holds for all $k \in \mN$ and still for $m$ even, that, (see \cite{dancing}), 
\begin{align*}
E_{m+2k-1} &=   \left(p_{2k-1}\, \ln r + q_{2k-1}   \right)\, \frac{\pi^{\frac{m}{2}+k}}{\Gamma(\frac{m}{2}+k)}\, U_{2k-1} \\[1mm]
E_{m+2k} &= \left(p_{2k}\, \ln r + q_{2k}   \right)\, \frac{\pi^{\frac{m}{2}+k}}{\Gamma(\frac{m}{2}+k)}\, T_{2k}
\end{align*}
the constants $(p_{2k-1} ,  q_{2k-1})$ and $(p_{2k} ,  q_{2k})$ satisfying the recurrence relations
\begin{eqnarray*}
\left\{\begin{array}{lll}
p_{2k} & = & \displaystyle{\frac{1}{2k}}\ p_{2k-1}\\[2mm]
q_{2k} & = & \displaystyle{\frac{1}{2k}\left(q_{2k-1}-\frac{1}{2k}\, p_{2k-1}\right)}
\end{array}\right.
\label{rec1}
\end{eqnarray*}
and
\begin{eqnarray*}
\left\{\begin{array}{lll}
p_{2k+1} & = & \displaystyle{-\, \frac{1}{2\pi}}\ p_{2k}\\[2mm]
q_{2k+1} & = & \displaystyle{-\, \frac{1}{2\pi}\left(q_{2k}-\frac{1}{m+2k}\, p_{2k}\right)}
\end{array}\right.
\label{rec2}
\end{eqnarray*}
with initial values
$$
p_0 = -\, \frac{1}{2^{m-1}\, \pi^m} \quad , \quad q_0 = 0
$$
\\

Through the signum--pair of operators $(\Delta , \bZ)$ or, equivalently, by iteration of the action of the operator $\uD$, we obtain the corresponding formulae for the operator $\bZ$ acting between signum--distributions:

\begin{eqnarray}
\label{signumlapsU} {\bf Z}\, ^{s}U_\lambda &=&    \lambda\,(\lambda+m-2)\,  ^{s}U_{\lambda-2}  \quad , \quad \lambda \neq -m+2, -m, -m-2, \ldots  \\
\label{signumlapsT} {\bf Z}\, ^{s}T_\lambda &=&  (\lambda-1)(\lambda+m-1)\, ^{s}T_{\lambda-2} \quad , \quad \lambda \neq -m+1, -m-1, \ldots
\end{eqnarray}
and
\begin{eqnarray}
\label{signumlapsUspecial}
\bZ\, ^{s}U_{-m-2\ell} &=&   (m+2\ell)(2\ell+2)\, ^{s}U_{-m-2\ell-2} + (-1)^\ell\, a_m\, \frac{(m+4\ell+2)(m+2\ell+2)}{C(m,\ell+1)}\, \uom\, \dirac^{2\ell+2}\, \delta(\ux)\\  \\
\label{signumlapsTspecial}
\bZ\, ^{s}T_{-m-2\ell+1} &=& (m+2\ell)(2\ell)\, ^{s}T_{-m-2\ell-1} + (-1)^{\ell+1}\, a_m\, \frac{m+4\ell}{C(m,\ell)}\, \uom\, \dirac^{2\ell+1}\, \delta(\ux)
\end{eqnarray}
and, in particular, for $\ell=0$
\begin{eqnarray}
\label{signumlapsUm}
\bZ\, ^{s}U_{-m} &=& 2 m\, ^{s}U_{-m-2} + a_m\, \frac{m+2}{2m}\, \uom\, \dirac^2\, \delta(\ux)  \\
\label{signumlapsTm1}
\bZ\, ^{s}T_{-m+1} &=& a_m\, \p_r\, \delta(\ux)
\end{eqnarray}
and also
\begin{eqnarray}
\label{signumlapsUm2}
\bZ\, ^{s}U_{-m+2} &=& -\, (m-2)\, a_m\, \uom\, \delta(\ux)
\end{eqnarray}


\newpage
\subsection{The operators $\mE$, $\Gamma$ and $\p_{\uom}$}


\noindent
The operators $\mE$ and $\Gamma$ are cartesian and their actions on the distributions $T_\la$ and $U_\la$ are uniquely determined. 
Combining the actions of the operators $\ux$ and $\dirac$ we find, through formula (\ref{eulergamma}), the following formulae:
\begin{equation}
\label{eulerT}
\boxed{
\mE\, T_{\la} = \la\, T_\la
}
\quad , \quad
\boxed{
\Gamma\, T_{\la} = 0
}
\quad , \quad \la \neq -m, -m-2,\ldots
\end{equation}
and
\begin{equation}
\label{eulerU}
\boxed{
\mE\, U_{\la} = \la\, U_\la
}
\quad , \quad
\boxed{
\Gamma\, U_{\la} = (m-1)\, U_\la
}
\quad , \quad \la \neq -m+1, -m-1, \ldots
\end{equation}
while for $\ell = 0, 1, 2, \ldots$
\begin{equation}
\label{eulerTspecial}
\boxed{
\mE\, T_{-m-2\ell} = -\, (m+2\ell)\, T_{-m-2\ell} + (-1)^\ell\, \frac{m+2\ell}{C(m,\ell)}\, a_m\, \dirac^{2\ell}\, \delta(\ux)
}
\quad , \quad
\boxed{
\Gamma\, T_{-m-2\ell} = 0
}
\end{equation}
and
\begin{equation}
\label{eulerUspecial}
\boxed{
\mE\, U_{-m-2\ell+1} = -\, (m+2\ell-1)\, U_{-m-2\ell+1} + (-1)^{\ell}\, \frac{1}{C(m,\ell-1)}\, a_m\, \dirac^{2\ell-1}\, \delta(\ux)
}
\quad , \quad
\boxed{
\Gamma\, U_{-m-2\ell+1} = (m-1)\, U_{-m-2\ell+1}
}
\end{equation}
and in particular, for $\ell=0$,
\begin{equation}
\label{eulerTm}
\boxed{
\mE\, T_{-m} = -\, m\, T_{-m} + a_m\, \delta(\ux)
}
\quad , \quad 
\boxed{
\Gamma\, T_{-m} = 0
}
\end{equation}
and
\begin{equation}
\label{eulerUm1}
\boxed{
\mE\, U_{-m+1} = (-m+1)\, U_{-m+1}
}
\quad , \quad 
\boxed{
\Gamma\,U_{-m+1} = (m-1)\, U_{-m+1}
}
\end{equation}\\

\noindent
Notice that the result  $\Gamma\, U_\la = (m-1)\, U_\la$ holds for all $\la \in \mC$. This is no surprise since it was proved in \cite{ffb} that for any signum--radial distribution $S$ it holds that $\Gamma\, S = (m-1)\, S$.\\

\noindent
Through the signum--pair of operators $(\mE , \mE)$ the actions of the Euler operator on the associated signumdistributions are readily obtained:
\begin{equation*}
\boxed{
\mE\, ^{s}U_{\la} = \la\, ^{s}U_\la
}
\quad , \quad \la \neq -m, -m-2,\ldots
\end{equation*}
and
\begin{equation*}
\boxed{
\mE\, ^{s}T_{\la} = \la\, ^{s}T_\la
}
\quad , \quad \la \neq -m+1, -m-1, \ldots
\end{equation*}
while for $\ell = 0, 1, 2, \ldots$
\begin{equation*}
\boxed{
\mE\, ^{s}U_{-m-2\ell} = -\, (m+2\ell)\, ^{s}U_{-m-2\ell} + (-1)^\ell\, \frac{m+2\ell}{C(m,\ell)}\, a_m\, \uom\, \dirac^{2\ell}\, \delta(\ux)
}
\end{equation*}
and
\begin{equation*}
\boxed{
\mE\, ^{s}T_{-m-2\ell+1} = -\, (m+2\ell-1)\, ^{s}T_{-m-2\ell+1} + (-1)^{\ell+1}\, \frac{1}{C(m,\ell-1)}\, a_m\, \uom\, \dirac^{2\ell-1}\, \delta(\ux)
}
\end{equation*}
and in particular, for $\ell=0$,
\begin{equation*}
\boxed{
\mE\, ^{s}U_{-m} = -\, m\, ^{s}U_{-m} + a_m\, \uom\, \delta(\ux)
}
\end{equation*}
and
\begin{equation*}
\boxed{
\mE\, ^{s}T_{-m+1} = (-m+1)\, ^{s}T_{-m+1}
}
\end{equation*}

\noindent
Through the signum--pair of operators $(\Gamma , -\, \p_{\uom}\, \uom)$ the above formulae for the action of the operator $\Gamma$ on the distributions $T_\la$ immediately lead to the following results, valid for all $\la \in \mC$:
\begin{itemize}
\item $\p_{\uom}\, T_\la = 0$
\item $\uom\, \p_{\uom}\, \uom\, ^{s}U_\la = 0$
\item $\p_{\uom}\, \uom\, ^{s}U_\la = 0$
\end{itemize}
the last formula leading to
\begin{equation}
\label{gammasU}
\boxed{ \Gamma\, ^{s}U_\la = (m-1)\, ^{s}U_\la} \quad , \quad \la \in \mC
\end{equation}
Similarly, the action of $\Gamma$ on the distributions $U_\la$ entails for all $\la \in \mC$:
\begin{itemize}
\item $\p_{\uom}\, U_\la = -\, (m-1)\, ^{s}T_\la$
\item $\uom\, \p_{\uom}\, \uom\, ^{s}T_\la = -\, (m-1)\, U_\la$
\item $\p_{\uom}\, \uom\, ^{s}T_\la = -\, (m-1)\, ^{s}T_\la$
\end{itemize}
this last formula leading to
\begin{equation}
\label{gammasT}
\boxed{ \Gamma\, ^{s}T_\la = 0} \quad , \quad \la \in \mC
\end{equation}


\newpage
\subsection{The operators $\p_{\uom}\, \uom$, $\uom\, \p_{\uom}\, \uom$, $\Delta^*$ and $\bZ^*$}


As was stated in Proposition \ref{propangop} the angular operators $\p_{\uom}^2$ and $\Delta^*$ are cartesian and their actions on (signum)distributions are uniquely determined.
As
$$
\p_{\uom}\, \uom = -\, (m-1)\, {\bf 1} + \Gamma
$$
it holds for all $\la \in \mC$ that
$$
\p_{\uom}\, \uom\, T_\la = -\, (m-1)\, T_\la
$$
and through the commutative diagram
$$
\begin{array}{ccccccccc}
  &&    & {}_{\p_{\uom}\, \uom}   &                             &  &&&\\
  &&  T_\lambda & \longrightarrow  &  -\,  (m-1)\, T_\lambda & &&& \\[1mm]
  &&\hspace{-8mm}{}_{- \uom} &\hspace{2mm}{} \hspace{13mm}{}_{ \p_{\uom} }&\hspace{6mm}{}_{-\uom}&&&&\\[-2mm]
  && \uparrow & \diagdown \hspace{-1.3mm} \nearrow & \uparrow &&&&\\[-1.1mm]
  && \downarrow & \diagup \hspace{-1.3mm} \searrow & \downarrow &&&&\\[-2mm]
  &&\hspace{-4.5mm}{}^{\uom} & \hspace{2mm}{} \hspace{15mm}{}^{\uom\, \p_{\uom}\, \uom}&\hspace{5mm}{}^{\uom}&&&&\\[1mm]
  &&  ^{s}U_\lambda & \longrightarrow &  -\, (m-1)\, ^{s}U_\lambda  &&&&\\[-1mm]
  &&               & {}_{\uom\, \p_{\uom} }  &           &&&&
 \end{array}
$$
we obtain for all $\la \in \mC$
\begin{itemize}
\item $\uom\, \p_{\uom}\, \uom\, T_\la = -\, (m-1)\, ^{s}U_\la$
\item $\p_{\uom}\, ^{s}U_\la = -\, (m-1)\, T_\la$
\end{itemize}
meanwhile confirming formula (\ref{gammasU}), viz.
$$
\Gamma\, ^{s}U_\la = (m-1)\, ^{s}U_\la
$$
Similarly it holds for all $\la \in \mC$ that
$$
\p_{\uom}\, \uom\,  U_\la = 0
$$
and hence, through the commutative diagram

$$
\begin{array}{ccccccccc}
  &&    & {}_{\p_{\uom}\, \uom}   &                             &  &&&\\
  &&  U_\lambda & \longrightarrow  &  0  & &&& \\[1mm]
  &&\hspace{-8mm}{}_{- \uom} &\hspace{2mm}{} \hspace{13mm}{}_{ \p_{\uom} }&\hspace{6mm}{}_{-\uom}&&&&\\[-2mm]
  && \uparrow & \diagdown \hspace{-1.3mm} \nearrow & \uparrow &&&&\\[-1.1mm]
  && \downarrow & \diagup \hspace{-1.3mm} \searrow & \downarrow &&&&\\[-2mm]
  &&\hspace{-4.5mm}{}^{\uom} & \hspace{2mm}{} \hspace{15mm}{}^{\uom\, \p_{\uom}\, \uom}&\hspace{5mm}{}^{\uom}&&&&\\[1mm]
  && -\, ^{s}T_\lambda & \longrightarrow &  0  &&&&\\[-1mm]
  &&               & {}_{\uom\, \p_{\uom} }  &           &&&&
 \end{array}
$$

\noindent
that, for all $\la \in \mC$
\begin{itemize}
\item $\uom\, \p_{\uom}\, \uom\, U_\la = 0$
\item $\p_{\uom}\, ^{s}T_\la = 0$
\end{itemize}
meanwhile confirming formula (\ref{gammasT}), viz.
$$
\Gamma\, ^{s}T_\la = 0
$$

\noindent
For the actions of the Laplace--Beltrami operator $\Delta^* = \uom\, \p_{\uom} - \p_{\uom}^2$ and its signum--partner $\bZ^*$, we now find, for all $\la \in \mC$,
\begin{eqnarray}
\Delta^*\, T_\la = 0 \quad &{\rm and}& \quad \bZ^*\, ^{s}U_\la = 0\\
\Delta^*\, U_\la =  -\, (m-1)\, U_\la \quad &{\rm and}& \quad \bZ^*\, ^{s}T_\la = -\, (m-1)\, ^{s}T_\la\\
\Delta^*\, ^{s}T_\la = 0 \quad &{\rm and}& \quad \bZ^*\, U_\la = 0\\
\Delta^*\, ^{s}U_\la =  -\, (m-1)\, ^{s}U_\la \quad &{\rm and}& \quad \bZ^*\, T_\la = -\, (m-1)\, T_\la
\end{eqnarray}


\newpage
\subsection{The operators $\invux$, $\invr$ and $\invrsq$}
\label{operatorsinvux}


Division of a distribution by $\ux$ is a cartesian operation which leads, in general, to an equivalence class of distributions. However, in view of the results of Section \ref{uniqueness}, we know that under the action of the operators $\invux$ and $\invr$ only the distribution $U_{-m+1}$ and the signumdistribution $^{s}T_{-m+1}$ will have a non--unique result. The same holds for $U_{-m+1}$, $^{s}T_{-m+1}$, $T_{-m+2}$ and $^{s}U_{-m+2}$ under the action of the operator $\invrsq$.\\

\noindent
We obtain the formulae
$$
\boxed{\invr\, T_\lambda = \, ^{s} T_{\lambda-1}} \quad {\rm and} \quad \boxed{\invr\, ^{s}U_\lambda =  U_{\lambda-1} } \quad , \quad \la \in \mC
$$
and
$$
\boxed{\invux\, T_\lambda = -\, U_{\lambda-1}} \quad {\rm and} \quad \boxed{\invux\, ^{s}U_\lambda = \, ^{s}T_{\lambda-1}} \quad , \quad \la \in \mC
$$
and
$$
\boxed{\invr\, U_\lambda = \, ^{s} U_{\lambda-1}} \quad {\rm and} \quad \boxed{\invr\, ^{s}T_\lambda =  T_{\lambda-1}} \quad , \quad \la \neq -m+1
$$
and also
$$
\boxed{\invux\, U_\lambda = T_{\lambda-1}} \quad {\rm and} \quad \boxed{{\invux\, ^{s}T_\lambda = -\, ^{s}U_{\lambda-1}}} \quad , \quad \la \neq -m+1
$$
For the exceptional case where $\la = -m+1$ we obtain the following  non--unique results:
$$
\boxed{ \invr\, U_{-m+1} =\,  ^{s}U_{-m} + \uom\, \delta(\ux)\, c }
$$
$$
\boxed{ \invr\, ^{s}T_{-m+1} =\,  T_{-m} +  \delta(\ux)\, c }
$$
and
$$
\boxed{ \invux\, U_{-m+1} =\,  T_{-m} + \delta(\ux)\, c }
$$
$$
\boxed{ \invux\, ^{s}T_{-m+1} =\, -\,  ^{s}U_{-m} -\, \uom\, \delta(\ux)\, c }
$$\\

\noindent
By iteration we find formulae for the divion by $r^2$:
$$
\boxed{\invrsq\, T_\lambda =\, \phantom{^{s}}T_{\la-2}} \qquad \boxed{\invrsq\, ^{s}U_\lambda =\, ^{s}U_{\la-2}}
\quad , \quad \la \neq -m+2
$$
and 
$$
\boxed{\invrsq\, ^{s}T_\lambda =\, ^{s}T_{\la-2}} \qquad \boxed{\invrsq\, U_\lambda =\,  \phantom{^{s}} U_{\la-2}}
\quad , \quad \la \neq -m+1
$$
For the exceptional case $\la = -m+2$ we have
$$
\invrsq\, T_{-m+2} = -\, \invux\, \invux\, T_{-m+2} = \invux\, U_{-m+1}
$$
whence
$$
\boxed{ \invrsq\, T_{-m+2} =\,  T_{-m} + \delta(\ux)\, c }  \qquad \boxed{ \invrsq\, ^{s}U_{-m+2} =\,  ^{s}U_{-m} + \uom\, \delta(\ux)\, c }
$$
For the exceptional case $\la=-m+1$ we have
$$
\invrsq\, U_{-m+1} = -\, \invux\, \invux\, U_{-m+1} = -\, \invux\, (T_{-m} + \delta(\ux)\, c)
$$
whence
$$
\boxed{ \invrsq\, ^{s}T_{-m+1} =\,  ^{s}T_{-m-1} - \frac{1}{m}\,  \p_r\, \delta(\ux)\, c } \qquad \boxed{ \invrsq\, U_{-m+1} =  U_{-m-1} - \frac{1}{m}\, \dirac \delta(\ux)\, c } 
$$
\\

The action of the components of the operator $\invux$, viz. $(-\, \frac{x_j}{r^2}), j=1,\ldots,m$, on all $T_\la$ and $U_\la$ is uniquely determined except for $U_{-m+1}$. We have readily
$$
\boxed{(-\, \frac{x_j}{r^2})\, T_\la = -\, x_j\, T_{\la-2}} \qquad  \boxed{(-\, \frac{x_j}{r^2})\, ^{s}U_\la = -\, x_j\, ^{s}U_{\la-2}}  \quad , \quad \la \in \mC
$$
Further for $\la \neq -m+1$ we have
$$
\boxed{(-\, \frac{x_j}{r^2})\, U_{\la} = -\,  x_j\, U_{\la-2} } \qquad \boxed{(-\, \frac{x_j}{r^2})\, ^{s}T_{\la} = -\, x_j\, ^{s}T_{\la-2}} \quad , \quad \la \neq -m+1
$$
while
for the exceptional case  $\la = -m+1$, it holds that
$$
(-\, \frac{x_j}{r^2})\, U_{-m+1} = -\, x_j\, (U_{-m-1} - \frac{1}{m}\, \dirac \delta(\ux)\, c)
$$
whence
$$
\boxed{(-\, \frac{x_j}{r^2})\,  U_{-m+1}  = -\, x_j\, U_{-m-1}  + \frac{1}{m}\, x_j\, \dirac\, \delta(\ux)\, c}  \qquad  \boxed{(-\, \frac{x_j}{r^2})\, ^{s}T_{-m+1} = -\, x_j\, ^{s}T_{-m-1} -  \omega_j\, \delta(\ux)\, c}
$$
or
$$
\boxed{(-\, \frac{x_j}{r^2})\,  U_{-m+1}  = -\, x_j\, U_{-m-1}  - \frac{1}{m}\, \delta(\ux)\, c\, e_j}  \qquad  \boxed{(-\, \frac{x_j}{r^2})\, ^{s}T_{-m+1} = -\, x_j\, ^{s}T_{-m-1} + \frac{1}{m}\, \uom\, \delta(\ux)\, c\, e_j}
$$


\newpage
\subsection{The operators $\omega_j, j=1,\ldots,m$}
\label{omegaj}


The multiplication operator $\omega_j, j=1,\ldots,m$ is a spherical operator whose actions are uniquely determined, see Definition \ref{multiplicomegaj} in Section \ref{spheroperators}.  It holds for all $\la \in \mC$ that
$$
\omega_j\, T_\la = \{ \uom\, T_\la  \} _j = \{ ^{s}U_\la  \} _j = x_j\, ^{s}T_{\la-1}
$$
and, through signum--association,
$$
\omega_j\, ^{s}U_\la = x_j\, U_{\la-1}
$$
It is verified at once that
$$
\uom\, T_\la = \Sigma_{j=1}^m\, e_j\, \omega_j\, T_\la =  \Sigma_{j=1}^m\, e_j\, x_j\, ^{s}T_{\la-1} = \ux\, ^{s}T_{\la-1} =\,  ^{s}U_\la
$$
and, 
$$
\uom\, ^{s}U_\la = \Sigma_{j=1}^m\, e_j\, \omega_j\, ^{s}U_\la =  \Sigma_{j=1}^m\, e_j\, x_j\, U_{\la-1} = \ux\, U_{\la-1} = -\, T_\la
$$

\noindent
In a similar way we find that
$$
\omega_j\, U_\la = \sum_k\, \omega_j\, \{  U_\la  \}_k\, e_k =   \sum_k\, \omega_j\, x_k\, T_{\la-1}\, e_k =   \sum_k\,  x_k\, x_j\, ^{s}T_{\la-2}\, e_k =
                   x_j\, ^{s}U_{\la-1}
$$
and, through signum--association,
$$
\omega_j\, ^{s}T_\la = x_j\, T_{\la-1}
$$

\noindent
By iteration it follows, for all $\la \in \mC$,  that
\begin{align*}
\omega_k\, \omega_j\, T_\la &= x_j\, x_k\, T_{\la-2}\\
\omega_k\, \omega_j\, U_\la &= x_j\, x_k\, U_{\la-2}\\
\omega_k\, \omega_j\, ^{s}T_\la &= x_j\, x_k\, ^{s}T_{\la-2}\\
\omega_k\, \omega_j\, ^{s}U_\la &= x_j\, x_k\, ^{s}U_{\la-2} 
\end{align*}


\newpage
\subsection{The operators  $\uom\, \p_r$, $\invr\, \p_{\uom}$ and $\p_r$}
\label{operatordiracparts}


As we know how to divide by $\ux$ the distributions $T_\la$ and $U_\la$, we are ready now to compute the actions of the radial and angular parts of the Dirac operator through the procedures (\ref{equirad}) and (\ref{equiang}) respectively. In general the results of those actions are equivalence classes of distributions. However we know from Section \ref{uniqueness} that  we may expect uniquely determined results, except for the distribution $U_{-m+1}$ and the signumdistribution $^{s}T_{-m+1}$.\\

\noindent
First we tackle the distributions $T_\la$.
For $\la \neq -m, -m-2, \ldots$ we find that $(\uom\, \p_r)\, T_\la$ is uniquely determined by
\begin{equation*}
(\uom\, \p_r)\, T_\la = -\, \left[ \invux\, \mE\, T_\la \right] = -\, \left[ \invux\, \la\, T_\la \right] = = -\, \invux\, \la\, T_\la = \la\, U_{\la-1}
\end{equation*}
In view of  (\ref{diracTlambda}) it follows that
$$
(\invr\, \p_{\uom})\, T_\la = 0 \quad , \quad \la \neq -m, -m-2, \ldots
$$
For the exceptional cases we find that the distributions $(\uom\, \p_r)\, T_{-m-2\ell}, \ell=0,1,2,\ldots$, are uniquely determined by
\begin{eqnarray*}
(\uom\, \p_r)\, T_{-m-2\ell} &=& -\, \left[ \invux\, \mE\, T_{-m-2\ell}   \right]\\
&=& -\, \left[ \invux\, \left( -\, (m+2\ell)\, T_{-m-2\ell} + (-1)^{\ell}\, \frac{m+2\ell}{C(m,\ell)}\, a_m\, \dirac^{2\ell}\, \delta(\ux)  \right)  \right]\\
&=& -\, (m+2\ell)\, U_{-m-2\ell-1} + (-1)^{\ell+1}\, \frac{1}{C(m,\ell)}\, a_m\, \dirac^{2\ell+1}\, \delta(\ux)
\end{eqnarray*}
In view of  (\ref{diracTspecial}) it follows that
$$
(\invr\, \p_{\uom})\, T_{-m-2\ell} = 0, \quad , \quad \ell = 0, 1, 2, \ldots
$$
In particular,
for $\la = -m$, we have
\begin{eqnarray*}
(\uom\, \p_r)\, T_{-m} &=& -\, \left[ \invux\, \mE\, T_{-m} \right] = -\, \left[ \invux\, (-m\, T_{-m} + a_m\, \delta(\ux)) \right]\\
&=& (-m)\, U_{-m-1} - \frac{1}{m}\, a_m\, \dirac\delta(\ux)
\end{eqnarray*}
and
$$
(\invr\, \p_{\uom})\, T_{-m} = 0
$$
which are in accordance with (\ref{diracTm}). \\

\noindent
Notice that the formula
$$
\boxed{(\invr\, \p_{\uom})\, T_{\la} = 0}
$$
thus holds for all $\la \in \mC$.\\

\noindent
Through the signum--pair of operators $[\uom\, \p_r , \uom\, \p_r ]$ and the corresponding commutative diagrams we then find the following formulae.\\

\noindent
For $\la \neq -m, -m-2, \ldots$ it holds that
\begin{itemize}
\item $\p_r\, T_\la = \la\, ^{s}T_{\la-1}$
\item $ \p_r\, ^{s}U_\la = \la\, U_{\la-1}$
\item $ (\uom\, \p_r)\, ^{s}U_\la = -\, \la\, ^{s}T_{\la-1}$
\end{itemize}

\noindent
For the exceptional values $\la = -\, m - 2\ell, \ell = 0, 1, 2, \ldots$ it holds that
\begin{itemize}
\item $\p_r\, T_{-m-2\ell} = -\, (m+2\ell)\, ^{s}T_{-m-2\ell-1} + (-1)^\ell\, \frac{1}{C(m,\ell)}\, a_m\, \uom\, \dirac^{2\ell+1}\, \delta(\ux)$
\item $ \p_r\, ^{s}U_{-m-2\ell} = -\, (m+2\ell)\, U_{-m-2\ell-1} + (-1)^{\ell+1}\, \frac{1}{C(m,\ell)}\, a_m\,  \dirac^{2\ell+1}\, \delta(\ux)$
\item $ (\uom\, \p_r)\, ^{s}U_{-m-2\ell} =  (m+2\ell)\, ^{s}T_{-m-2\ell-1} + (-1)^{\ell+1}\, \frac{1}{C(m,\ell)}\, a_m\, \uom\, \dirac^{2\ell+1}\, \delta(\ux)$
\end{itemize}
and in particular, for $\ell=0$,
\begin{itemize}
\item $\p_r\, T_{-m} = (-m)\, ^{s}T_{-m-1} - \frac{1}{m}\, a_m\, \p_r\, \delta(\ux)$
\item $ \p_r\, ^{s}U_{-m} = (-m)\, U_{-m-1} - \frac{1}{m}\, a_m\, \dirac\, \delta(\ux) $
\item $ (\uom\, \p_r)\, ^{s}U_{-m} = m\, ^{s}T_{-m-1} + \frac{1}{m}\, a_m\, \p_r\, \delta(\ux)$
\end{itemize}

\noindent
Through the signum--pair of operators $[\invr\, \p_{\uom} , -\, \invr\, \uom\, \p_{\uom}\, \uom]$ and the corresponding commutative diagrams we obtain that for all $\la \in \mC$ 
\begin{itemize}
\item $ \invr\, \Gamma\, T_\la = 0$
\item $ \invr\, \Gamma\, ^{s}U_{\la} = (m-1)\, U_{\la-1}$
\item$\boxed{ (\invr\, \p_{\uom})\, ^{s}U_\la = -\, (m-1)\, ^{s}T_{\la-1}}$ 
\end{itemize}
\vspace*{5mm}

\noindent
Now let us turn our attention to the distributions $U_\la$.\\

\noindent
For $\la \neq -m+1, -m-1, \ldots$ we find the unique expressions
\begin{equation*}
(\uom\, \p_r)\, U_\la = -\, \left[ \invux\, \mE\, U_\la \right] = -\, \left[ \invux\, \la\, U_\la \right] = -\, \la\, T_{\la-1}
\end{equation*}
and
\begin{equation*}
(\invr\, \p_{\uom})\, U_\la = -\, \left[ \invux\, \Gamma\, U_\la \right] = -\, \left[ \invux\, (m-1)\, U_\la \right] = -\, (m-1)\, T_{\la-1}
\end{equation*}
which are in accordance with (\ref{diracUlambda}).\\
In the exceptional cases we find, using (\ref{eulerUspecial}), for $\ell = 1, 2, \ldots$, the unique expressions
\begin{eqnarray*}
(\uom\, \p_r)\, U_{-m-2\ell+1} &=& -\, \left[  \invux\, \mE\, U_{-m-2\ell+1} \right]\\
&=& -\, \left[ \invux\, \left(  -\, (m+2\ell-1)\, U_{-m-2\ell+1} + (-1)^{\ell}\, \frac{1}{C(m,\ell-1)}\, a_m\, \dirac^{2\ell-1}\, \delta(\ux) \right)  \right]\\
&=& (m+2\ell-1)\, T_{-m-2\ell} + (-1)^{\ell+1}\, \frac{m+2\ell}{C(m,\ell)}\, a_m\, \dirac^{2\ell}\, \delta(\ux)
\end{eqnarray*}
and
\begin{eqnarray*}
(\invr\, \p_{\uom})\, U_{-m-2\ell+1} &=& -\, \left[  \invux\, \Gamma\, U_{-m-2\ell+1} \right]\\
&=& -\, \left[ \invux\, (m-1)\, U_{-m-2\ell+1} \right] \\
&=& -\, (m-1)\, T_{-m-2\ell}
\end{eqnarray*}
which are in accordance with (\ref{diracUspecial}).\\
In the particular case where $\la=-m+1$ we find, using (\ref{eulerUm1}), the entangled expressions
\begin{eqnarray}
\label{raddiracUm1}
(\uom\, \p_r)\, U_{-m+1} &=& -\, \left[ \invux\, \mE\, U_{-m+1} \right] = -\, \left[ \invux\, (-m+1)\, U_{-m+1} \right] \nonumber\\[2mm]
&=& (m-1)\, T_{-m}  + \delta(\ux)\, c_1
\end{eqnarray}
\begin{eqnarray}
\label{angdiracUm1}
(\invr\, \p_{\uom})\, U_{-m+1} &=& -\, \left[ \invux\, \Gamma\, U_{-m+1} \right] = -\, \left[ \invux\, (m-1)\, U_{-m+1} \right]\nonumber\\[2mm]
&=& -\,  (m-1)\, T_{-m}  + \delta(\ux)\, c_2
\end{eqnarray}
where, seen (\ref{diracUm11}), the arbitrary constants $c_1$ and $c_2$ must satisfy the entanglement condition
\begin{equation}
\label{entanglecond}
c_1 + c_2 = -\, a_m
\end{equation}

\begin{remark}
{\rm
Just keep in mind that in subsection 8.8 we will fix a particular choice for the constants $c_1$ and $c_2$ satisfying (\ref{entanglecond}) in this way {\em defining } both $(\uom\, \p_r)\, U_{-m+1}  $  and $ (\invr\, \p_{\uom})\, U_{-m+1} $.
}
\end{remark}\hfill\\

Through the actions of the operators $\uom\, \p_r$ and $\invr\, \p_{\uom}$ it is possible to recover the formuale for the action of the Laplace operator.

\begin{itemize}
\item For $\la \neq -m+2, -m, -m-2, \ldots$ one has
\begin{align*}
(\uom\, \p_r)^2\, T_\la &= (\uom\, \p_r)\, (\la\, U_{\la-1}) = -\, \la\, (\la-1)\, T_{\la-2}\\
(\invr\, \p_{\uom})\, (\uom\, \p_r)\,  T_\la &= (\invr\, \p_{\uom})\, (\la\, U_{\la-1}) = -\, \la\, (m-1)\, T_{\la-2}
\end{align*}
Keeping in mind that $(\invr\, \p_{\uom})\, T_\la = 0$, it is confirmed that
$$
\Delta\, T_\la = \la\, (\la+m-2)\, T_{\la-2}
$$
\item For $\la = -m+2$ one has
\begin{align*}
(\uom\, \p_r)^2\, T_{-m+2} &= (\uom\, \p_r)\, ((-m+2)\, U_{-m+1}) = -\, (m-2)\,( (m-1)\, T_{-m} - a_m\, \delta(\ux) )\\
(\invr\, \p_{\uom})\,  (\uom\, \p_r)\, T_{-m+2} &= (\invr\, \p_{\uom})\, ((-m+2)\, U_{-m+1}) = (m-2) (m-1)\, T_{-m}
\end{align*}
confirming that
$$
\Delta\, T_{-m+2} = -\,  (m-2)\, a_m\, \delta(\ux)
$$
\item For $\la = -m-2\ell$ one has
\begin{align*}
(\uom\, \p_r)^2\, T_{-m-2\ell} &= (\uom\, \p_r)\, \left(-\, (m+2\ell)\, U_{-m-2\ell-1} + (-1)^{\ell+1}\,  \frac{(2\ell+2)(m+2\ell+2)}{C(m, \ell+1)}\, a_m\, \dirac^{2\ell+1} \delta(\ux) \right)\\[2mm]
&= -\, (m+2\ell)(m+2\ell+1)\, T_{-m-2\ell-2} + (-1)^{\ell+1}\, \frac{m+2\ell+2}{C(m,\ell+1)}\, (2m+4\ell+1)\, a_m\, \dirac^{2\ell+2}\, \delta(\ux)\\[3mm]
(\invr\, \p_{\uom})\, (\uom\, \p_r)\,  T_{-m-2\ell} &= (\invr\, \p_{\uom})\, \left(-\, (m+2\ell)\, U_{-m-2\ell-1} + (-1)^{\ell+1}\,  \frac{(2\ell+2)(m+2\ell+2)}{ C(m, \ell+1)}\, a_m\, \dirac^{2\ell+1} \delta(\ux) \right)\\[2mm]
&= -\, (m+2\ell)(1-m)\, T_{-m-2\ell-2} + (-1)^{\ell+1}\, \frac{m+2\ell+2}{C(m, \ell+1)}\, (1-m)\, a_m\, \dirac^{2\ell+2}\, \delta(\ux)
\end{align*}
confirming that
$$
\Delta\, T_{-m-2\ell} =   (m+2\ell)(2\ell+2)\, T_{-m-2\ell-2} + (-1)^\ell\,  \frac{(m+4\ell+2)(m+2\ell+2)}{C(m,\ell+1)}\, a_m\, \dirac^{2\ell+2}\, \delta(\ux)\\
$$
\item For $\la \neq -m+1, -m-1, -m-3, \ldots$ one has
\begin{align*}
(\uom\, \p_r)^2\, U_\la &= (\uom\, \p_r)\, (-\, \la\, T_{\la-1}) = -\, \la\, (\la-1)\, U_{\la-2}\\
(\invr\, \p_{\uom})\, (\uom\, \p_r)\, U_\la &=    (\invr\, \p_{\uom})\,   (-\, \la\, T_{\la-1}) = 0\\
(\uom\, \p_r)\, (\invr\, \p_{\uom})\, U_\la &=(\uom\, \p_r)\, ( -\, (m-1)\, T_{\la-1}) = -\, (m-1)\, (\la-1)\, U_{\la-2}\\
(\invr\, \p_{\uom})\, (\invr\, \p_{\uom})\, U_\la &= (\invr\, \p_{\uom})\, (-\, (m-1)\, T_{\la-1})  = 0
\end{align*}
confirming that
$$
\Delta\, U_\la = (\la-1)(\la+m-1)\, U_{\la-2}
$$
\item For $\la = -m+1$ one has
\begin{align*}
(\uom\, \p_r)^2\, U_{-m+1} &= (\uom\, \p_r)\, ((m-1)\, T_{-m} - a_m\, \delta(\ux)) = -\, (m-1)m\, U_{-m-1} - \frac{2m-1}{m}\, a_m\, \dirac\, \delta(\ux)\\
(\invr\, \p_{\uom})\, (\uom\, \p_r)\, U_{-m+1} &=    (\invr\, \p_{\uom})\,   ((m-1)\, T_{-m} - a_m\, \delta(\ux)) = 0\\
(\uom\, \p_r)\, (\invr\, \p_{\uom})\, U_{-m+1} &=(\uom\, \p_r)\, (-\, (m-1)\, T_{-m}) = (m-1)m\, U_{-m-1} + \frac{m-1}{m}\, a_m\, \dirac\, \delta(\ux)\\
(\invr\, \p_{\uom})\, (\invr\, \p_{\uom})\, U_{-m+1} &= (\invr\, \p_{\uom})\, (-\, (m-1)\, T_{-m})  = 0
\end{align*}
confirming that
$$
\Delta\, U_{-m+1} = a_m\, \dirac\, \delta(\ux)
$$
For $\la = -m-2\ell+1$ one has 
\begin{align*}
(\uom\, \p_r)^2\, U_{-m+1} &=   -\, (m+2\ell-1)(m+2\ell)\, U_{-m-2\ell-1} + (-1)^{\ell+1}\, \frac{2m+4\ell-1}{C(m,\ell)}\, a_m\, \dirac^{2\ell+1}  \\
(\invr\, \p_{\uom})\, (\uom\, \p_r)\, U_{-m+1} &= 0\\
(\uom\, \p_r)\, (\invr\, \p_{\uom})\, U_{-m+1} &=  (m+2\ell)(m-1)\, U_{-m-2\ell-1} - (-1)^{\ell+1}\, \frac{m-1}{C(m,\ell)}\, a_m\, \dirac^{2\ell+1}\, \delta(\ux)\\
(\invr\, \p_{\uom})\, (\invr\, \p_{\uom})\, U_{-m+1} &= 0
\end{align*}
confirming that
$$
\Delta\, U_{-m-2\ell+1} = (m+2\ell)(2\ell)\, U_{-m-2\ell-1} + (-1)^\ell\, a_m\, \frac{m+4\ell}{C(m,\ell)}\, \dirac^{2\ell+1}\, \delta(\ux)
$$
\end{itemize}\hfill\\

Through the signum--pair of operators $[\uom\, \p_r , \uom\, \p_r]$ and the corresponding commutative diagrams, we find the following formulae.\\

\noindent
For $\la \neq -m+1, -m-1, \ldots$ it holds that
\begin{itemize}
\item $\p_r\, U_\la = \la\, ^{s}U_{\la-1}$
\item $\p_r\, ^{s}T_\la = \la\, T_{\la-1}$
\item $(\uom\, \p_r)\, ^{s}T_\la = \la\, ^{s}U_{\la-1}$
\end{itemize}

\noindent
For the exceptional values $\la = -m-2\ell+1, \ell=1,2,\ldots$ it holds that
\begin{itemize}
\item $\p_r\, U_{-m-2\ell+1} = -\, (m+2\ell-1)\, ^{s}U_{-m-2\ell} + (-1)^\ell\, \frac{m+2\ell}{C(m,\ell)}\, a_m\, \uom\, \dirac^{2\ell}\, \delta(\ux)$
\item $\p_r\, ^{s}T_{-m-2\ell+1} = -\,  (m+2\ell-1)\, T_{-m-2\ell} + (-1)^\ell\, \frac{m+2\ell}{C(m,\ell)}\, a_m\,  \dirac^{2\ell}\, \delta(\ux)$
\item $(\uom\, \p_r)\, ^{s}T_{-m-2\ell+1} = -\, (m+2\ell-1)\, ^{s}U_{-m-2\ell} + (-1)^\ell\, \frac{m+2\ell}{C(m,\ell)}\, a_m\, \uom\, \dirac^{2\ell}\, \delta(\ux)$
\end{itemize}
while for $\ell=0$ it holds that
\begin{itemize}
\item $\p_r\, U_{-m+1} = -\, (m-1)\, ^{s}U_{-m} - \uom\, \delta(\ux)\, c_1$
\item $\p_r\, ^{s}T_{-m+1} = -\, (m-1)\, T_{-m} - \delta(\ux)\, c_1$
\item $(\uom\, \p_r)\, ^{s}T_{-m+1} = -\, (m-1)\, ^{s}U_{-m} - \uom\, \delta(\ux)\, c_1$
\end{itemize}
\vspace*{5mm}

\noindent
Through the signum--pair of operators $[\invr\, \p_{\uom} , -\, \invr\, \uom\, \p_{\uom}\, \uom]$ and the corresponding commutative diagrams, we find the following formulae.\\

\noindent
For $\la \neq -m+1, -m-1, \ldots$ it holds that
\begin{itemize}
\item $(-\, \invr\, \Gamma + (m-1)\, \invr)\, ^{s}T_\la = (m-1)\, T_{\la-1}$
\item $\invr\, \Gamma\, U_\la = (m-1)\, ^{s}U_{\la-1}$
\item $(-\, \invr\, \p_{\uom} + (m-1)\, \invr\, \uom)\, ^{s}T_\la = (m-1)\, ^{s}U_{\la-1}$
\end{itemize}
or equivalently
\begin{itemize}
\item $\invr\, ^{s}T_\la = T_{\la-1}$
\item $\invr\, \Gamma\, U_\la = (m-1)\, ^{s}U_{\la-1}$
\item $\invr\, \p_{\uom}\, ^{s}T_\la = 0$
\end{itemize}

\noindent
For the exceptional values $\la = -m-2\ell+1, \ell=1, 2, \ldots$ it holds that
\begin{itemize}
\item $(-\, \invr\, \Gamma + (m-1)\, \invr)\, ^{s}T_{-m-2\ell+1} = (m-1)\, T_{m-2\ell}$
\item $\invr\, \Gamma\, U_{-m-2\ell+1} = (m-1)\, ^{s}U_{m-2\ell}$
\item $(-\, \invr\, \p_{\uom} + (m-1)\, \invr\, \uom)\, ^{s}T_{-m-2\ell+1}= (m-1)\, ^{s}U_{m-2\ell}$
\end{itemize}
or equivalently
\begin{itemize}
\item $\invr\, ^{s}T_{-m-2\ell+1} = T_{m-2\ell}$
\item $\invr\, \Gamma\, U_{-m-2\ell+1}= (m-1)\, ^{s}U_{m-2\ell}$
\item $\invr\, \p_{\uom}\, ^{s}T_{-m-2\ell+1} = 0$
\end{itemize}

For the exceptional value $\la = -m+1$ it holds that
\begin{itemize}
\item $(-\, \invr\, \Gamma + (m-1)\, \invr)\, ^{s}T_{-m+1} = (m-1)\, T_{-m} - \delta(\ux)\, c_2$
\item $\invr\, \Gamma\, U_{-m+1} = (m-1)\, ^{s}U_{-m} + \uom\, \delta(\ux)\, c_2$
\item $(-\, \invr\, \p_{\uom} + (m-1)\, \invr\, \uom)\, ^{s}T_{-m+1} = (m-1)\, ^{s}U_{-m} - \uom\, \delta(\ux)\, c_2$
\end{itemize}
or, equivalently,
\begin{itemize}
\item $ (m-1)\, \invr\, ^{s}T_{-m+1} = (m-1)\, T_{-m} - \delta(\ux)\, c_2$
\item $\invr\, \Gamma\, U_{-m+1} = (m-1)\, ^{s}U_{-m} + \uom\, \delta(\ux)\, c_2$
\item $ (\invr\, \p_{\uom})\, ^{s}T_{-m+1} = 0$
\end{itemize}

\noindent
Notice that the formula
$$
\boxed{ (\invr\, \p_{\uom})\, ^{s}T_{\la} = 0}
$$
is valid for all $\la \in \mC$.\\

As expected $U_{-m+1}$ is the only distribution in the whole of the two families of distributions $T_\la$ and $U_\la$, $\la \in \mC$, for which the actions of the operators $\dirac_{rad}$ and $\dirac_{ang}$ are {\em not} uniquely determined, and a similar observation holds for the signumdistribution $^{s}T_{-m+1}$. Nothing prevents us from defining these actions by choosing ourselves specific values for the arbitrary constants involved.  However,  in order to avoid possible inconsistencies, we will, for the moment, not do so and keep the constants $c_1$ and $c_2$ arbitrary but constraint by the entanglement condition (\ref{entanglecond}), viz.: $c_1 + c_2 = -\, a_m$.
Notice that this is made possible by the fact that $\invux\, U_{-m+1}$ was expressed as an equivalence class of distributions, viz.
$$
\invux\, U_{-m+1} = T_{-m} + \delta(\ux)\, c
$$
which is really necessary since fixing here and now a particular value for this arbitrary constant $c$ would lead to the impossible result $\dirac\, U_{-m+1} = 0$.


\newpage
\subsection{The operators $\uD$ and $\bZ$}


Combining the results obtained in the preceding Subsection \ref{operatordiracparts} for the actions of the radial and angular parts of the Dirac operator, we obtain the action of the operator $\uD$ on the distributions $T_\la$ and $U_\la$ and of the Dirac operator $\dirac$ on the signumdistributions $^{s}T_\la$ and $^{s}U_\la$. These results should match each other through the signum--pair of operators $[\uD , \dirac]$. We expect $U_{-m+1}$ to be the only distribution for which the $\uD$--action is not uniquely determined and the same for the signumdistribution $^{s}T_{-m+1}$ with respect to the $\dirac$--action. Follows an overview of the results, recalling, for completeness' sake, some formulae already obtained in the preceding subsections.\\

First we consider the case of the distributions $T_\la$ and the associated signumdistributions $^{s}U_\la$.\\

\noindent
For $\boxed{\lambda \neq -m, -m-2, \ldots}$ it holds that
\begin{eqnarray}
(\uom\, \p_r)\, T_\la &=& \la\, U_{\la-1}\\
(\invr\, \p_{\uom})\, T_\la &=& 0\\
\dirac\, T_\la &=& \la\, U_{\la-1}\\
\label{uDT} \uD\, T_\lambda &=&  (\lambda+m-1)\, U_{\lambda-1} 
\end{eqnarray}
and
\begin{eqnarray}
(\uom\, \p_r)\, ^{s}U_\la &=& -\, \la\, ^{s}T_{\la-1}\\
(\invr\, \p_{\uom})\, ^{s}U_\la &=& -\, (m-1)\, ^{s}T_{\la-1}\\
\label{uDsU} \uD\, ^{s}U_\la &=& -\, \la\, ^{s}T_{\la-1}\\
\label{diracsU} \dirac\, ^{s}U_\lambda &=& -\,  (\lambda+m-1)\, ^{s}T_{\lambda-1} 
\end{eqnarray}

\noindent
For $\boxed{\lambda = -m}$ it holds that
\begin{eqnarray}
(\uom\, \p_r)\, T_{-m} &=& -m\, U_{-m-1} - \frac{1}{m}\, a_m\, \dirac\, \delta(\ux)\\
(\invr\, \p_{\uom})\, T_{-m} &=& 0\\
\dirac\, T_{-m} &=& -m\, U_{-m-1} - \frac{1}{m}\, a_m\, \dirac\, \delta(\ux)\\
\label{uDT} \uD\, T_{-m} &=&  -\, U_{-m-1} - \frac{1}{m}\, a_m\, \dirac\, \delta(\ux)
\end{eqnarray}
and
\begin{eqnarray}
(\uom\, \p_r)\, ^{s}U_{-m} &=&  m\, ^{s}T_{-m-1} + \frac{1}{m}\, a_m\, \p_r\, \delta(\ux) \\
(\invr\, \p_{\uom})\, ^{s}U_{-m} &=& -\, (m-1)\, ^{s}T_{-m-1}\\
\label{uDsU} \uD\, ^{s}U_{-m} &=& m\, ^{s}T_{-m-1} + \frac{1}{m}\, a_m\, \p_r\, \delta(\ux)\\
\dirac\,  ^{s}U_{-m} &=&  ^{s}T_{-m-1} + \frac{1}{m}\, a_m\, \p_r\, \delta(\ux)
\end{eqnarray}

\noindent
More generally for $\boxed{\la = -m-2\ell}\, , \ell=0,1,2,\ldots$, it holds that
\begin{eqnarray}
(\uom\, \p_r)\, T_{-m-2\ell} &=& -\, (m+2\ell)\, U_{-m-2\ell-1} + (-1)^{\ell+1}\, \frac{1}{C(m,\ell)}\, a_m\, \dirac^{2\ell+1}\, \delta(\ux)\\
(\invr\, \p_{\uom})\, T_{-m-2\ell} &=& 0\\
\dirac\, T_{-m-2\ell} &=& -\, (m+2\ell)\, U_{-m-2\ell-1} + (-1)^{\ell+1}\, \frac{1}{C(m,\ell)}\, a_m\, \dirac^{2\ell+1}\, \delta(\ux)\\
\uD\, T_{-m-2\ell} &=&  -\, (2\ell+1)\, U_{-m-2\ell-1} + (-1)^{\ell+1}\, \frac{1}{C(m,\ell)}\, a_m\, \dirac^{2\ell+1}\, \delta(\ux)
\end{eqnarray}
and
\begin{eqnarray}
(\uom\, \p_r)\, ^{s}U_{-m-2\ell} &=&  (m+2\ell)\, ^{s}T_{-m-2\ell-1} + (-1)^{\ell+1}\, \frac{1}{C(m,\ell)}\, a_m\, \uom\, \dirac^{2\ell+1}\, \delta(\ux)\\
(\invr\, \p_{\uom})\, ^{s}U_{-m-2\ell} &=& -\, (m-1)\, ^{s}T_{-m-2\ell-1}\\
\uD\, ^{s}U_{-m-2\ell} &=&  (m+2\ell)\, ^{s}T_{-m-2\ell-1} + (-1)^{\ell+1}\, \frac{1}{C(m,\ell)}\, a_m\, \uom\, \dirac^{2\ell+1}\, \delta(\ux)\\
\label{diracsUspecial} \dirac\, ^{s}U_{-m-2\ell} &=&  (2\ell+1)\, ^{s}T_{-m-2\ell-1} + (-1)^{\ell+1}\, \frac{1}{C(m,\ell)}\, a_m\, \uom\, \dirac^{2\ell+1}\, \delta(\ux)
\end{eqnarray}\\

Now we treat the case of the distributions $U_\la$ and the associated signumdistributions $^{s}T_\la$.\\

\noindent
For $\boxed{\lambda \neq -m+1, -m-1, \ldots}$ it holds that
\begin{eqnarray}
(\uom\, \p_r)\, U_\la &=& -\,  \la\, T_{\la-1}\\
(\invr\, \p_{\uom})\, U_\la &=& -\, (m-1)\, T_{\la-1}\\
\dirac\, U_\la &=& -\,  (\lambda+m-1)\, T_{\lambda-1} \\
\label{uDU} \uD\, U_\lambda &=& -\,  \la\, T_{\la-1}
\end{eqnarray}
and
\begin{eqnarray}
 (\uom\, \p_r)\, ^{s}T_\la &=&  \la\, ^{s}U_{\la-1}\\
 (\invr\, \p_{\uom})\, ^{s}T_\la &=& 0\\
\label{uDsT} \uD\, ^{s}T_\la &=&   (\lambda+m-1)\, ^{s}U_{\lambda-1}\\
\label{diracsT}  \dirac\, ^{s}T_\lambda &=& \la\, ^{s}U_{\la-1}
\end{eqnarray}

\noindent
For $\boxed{\la = -m-2\ell+1}\, , \ell=1,2,\ldots$, it holds that
\begin{eqnarray}
(\uom\, \p_r)\, U_{-m-2\ell+1} &=& (m+2\ell-1)\, T_{-m-2\ell} + (-1)^{\ell+1}\, \frac{m+2\ell}{C(m,\ell)}\, a_m\, \dirac^{2\ell}\, \delta(\ux)\\
(\invr\, \p_{\uom})\, U_{-m-2\ell+1} &=& -\, (m-1)\, T_{-m-2\ell}\\
\dirac\, U_{-m-2\ell+1} &=&  (2\ell)\, T_{-m-2\ell} + (-1)^{\ell+1}\, \frac{m+2\ell}{C(m,\ell)}\, a_m\, \dirac^{2\ell}\, \delta(\ux)\\
 \uD\, U_{-m-2\ell+1} &=&   (m+2\ell-1)\, T_{-m-2\ell} + (-1)^{\ell+1}\, \frac{m+2\ell}{C(m,\ell)}\, a_m\, \dirac^{2\ell}\, \delta(\ux)
\end{eqnarray}
and
\begin{eqnarray}
(\uom\, \p_r)\, ^{s}T_{-m-2\ell+1} &=& -\,  (m+2\ell-1)\, ^{s}U_{-m-2\ell} + (-1)^{\ell}\, \frac{m+2\ell}{C(m,\ell)}\, a_m\, \uom\, \dirac^{2\ell}\, \delta(\ux)\\
(\invr\, \p_{\uom})\, ^{s}T_{-m-2\ell+1} &=& 0\\
\label{uDsTsing} \uD\, ^{s}T_{-m-2\ell+1} &=& -\, (2\ell)\, ^{s}U_{-m-2\ell} + (-1)^{\ell}\, \frac{m+2\ell}{C(m,\ell)}\, a_m\, \uom\, \dirac^{2\ell}\, \delta(\ux) \\
\label{diracsTsing}  \dirac\, ^{s}T_{-m-2\ell+1} &=&  -\,  (m+2\ell-1)\, ^{s}U_{-m-2\ell} + (-1)^{\ell}\, \frac{m+2\ell}{C(m,\ell)}\, a_m\, \uom\, \dirac^{2\ell}\, \delta(\ux)
\end{eqnarray}

\noindent
For the exceptional case $\boxed{\la = -m+1}$ we find that
\begin{eqnarray*}
 (\uom\, \p_r)\, U_{-m+1} &=& (m-1)\, T_{-m} + \delta(\ux)\, c_1\\
 (\invr\, \p_{\uom})\, U_{-m+1} &=& -\,  (m-1)\, T_{-m} + \delta(\ux)\, c_2\\
 \dirac\, U_{-m+1} &=& -\, a_m\, \delta(\ux) \\
 \uD\, U_{-m+1} &=&  (m-1)\, T_{-m} + \delta(\ux)\, c^*
\end{eqnarray*}
and
\begin{eqnarray*}
(\uom\, \p_r)\, ^{s}T_{-m+1} &=&  -\, (m-1)\, ^{s}U_{-m} - \uom\, \delta(\ux)\, c_1\\
(\invr\, \p_{\uom})\, ^{s}T_{-m+1} &=& 0\\
 \uD\, ^{s}T_{-m+1} &=&   \uom\, \delta(\ux)\, a_m\\
 \dirac\, ^{s}T_{-m+1} &=& -\, (m-1)\, ^{s}U_{-m} - \uom\, \delta(\ux)\, c^*
\end{eqnarray*}
Apparently the constant $c^*$ may be chosen freely, whereas the constants $c_1$ and $c_2$ must satisfy the entanglement condition (\ref{entanglecond}), viz.
$$
c_1 + c_2 = -\, a_m
$$
We will now fix these arbitrary constants ourselves, which means that we are going to {\em define} $(\uom\, \p_r)\, U_{-m+1} $, $(\invr\, \p_{\uom})\, U_{-m+1}$ and $\uD\, U_{-m+1}$. Quite naturally we have to be very cautious and constantly watch over the concistency. We make the following choice:
$$
c^* = c_1 = -\, a_m \quad , \quad c_2 = 0
$$
Are there arguments for this particular choice? Concerning the choice for $c^*$, one can see that for general $\la$ the expressions for $(\uom\, \p_r)\, U_{\la} $ and $\uD\, U_{\la}$ coincide, inspiring the choice $c^* = c_1$. Next, putting $\ell = 0$ in the expression for $(\uom\, \p_r)\, U_{-m-2\ell+1} $ results into $$(m-1)\, T_{-m} - a_m\, \delta(\ux)$$ which then can serve as the definition for $(\uom\, \p_r)\, U_{-m+1}$ making $c_1 = -\, a_m$ and $c_2 = 0$. Hence the following definition.

\begin{definition}
One defines
\begin{eqnarray}
\label{omegadrUm1} (\uom\, \p_r)\, U_{-m+1} &=& (m-1)\, T_{-m} - a_m\, \delta(\ux)\\
\label{invrdomegaUm1} (\invr\, \p_{\uom})\, U_{-m+1} &=& -\,  (m-1)\, T_{-m}\\
\label{diracUm1} \dirac\, U_{-m+1} &=& -\, a_m\, \delta(\ux) \\
\label{uDUm1} \uD\, U_{-m+1} &=&  (m-1)\, T_{-m} - a_m\,  \delta(\ux)
\end{eqnarray}
and
\begin{eqnarray}
(\uom\, \p_r)\, ^{s}T_{-m+1} &=&  -\, (m-1)\, ^{s}U_{-m} + a_m\,  \uom\, \delta(\ux)\\
(\invr\, \p_{\uom})\, ^{s}T_{-m+1} &=& 0\\
\label{uDsTm1} \uD\, ^{s}T_{-m+1} &=&  a_m\,  \uom\, \delta(\ux)\\
\label{diracsTm1} \dirac\, ^{s}T_{-m+1} &=& -\, (m-1)\, ^{s}U_{-m} + a_m\, \uom\, \delta(\ux)
\end{eqnarray}
\end{definition}

\noindent
Notice by the way that  formula (\ref{invrdomegaUm1})  matches the corresponding formula for all other values of the parameter $\la$, viz.
$$
(\invr\, \p_{\uom})\, U_{\la} = -\, (m-1)\, T_{\la-1}
$$

Iteration of the operators $\dirac$ and $\uD$ yield the following formulae involving the Laplace  and $\bZ$ operators. These results should match each other pairwise through the signum--pairs of operators $(\Delta , \bZ)$ and $[\bZ , \Delta]$. Based on the results of Section \ref{uniqueness} we expect the distributions $T_{-m+2}$ and $U_{-m+1}$ to have non--unique results under the action of the operator $\bZ$; the same holds for the signumdistributions $^{s}U_{-m+2}$ and $^{s}T_{-m+1}$ under the action of the Laplace operator $\Delta$.\\

First we consider the case of the distributions $T_\la$ and the associated signumdistributions $^{s}U_\la$.\\

\noindent
For $\boxed{\lambda \neq -m+2, -m, -m-2, \ldots}$ we already found formulae (\ref{lapT}) and (\ref{signumlapsU}), viz.
\begin{eqnarray*}
\Delta\, T_\la &=&  \la\, (\la+m-2)\, T_{\la-2}   \\
\bZ\, ^{s}U_\la &=& \la\, (\la+m-2)\, ^{s}U_{\la-2} 
\end{eqnarray*}
and now it also holds that
\begin{eqnarray}
\label{signumlapT} \bZ\, T_\la &=& (\lambda-1)(\lambda+m-1)\, T_{\lambda-2} \\
\label{lapsU} \Delta\, ^{s}U_\lambda &=&   (\lambda-1)(\lambda+m-1)\, ^{s}U_{\lambda-2} 
\end{eqnarray}

\noindent
For $\boxed{\la = -m+2}$ we already found formulae (\ref{lapTm2}) and (\ref{signumlapsUm2}), viz.
\begin{eqnarray*}
\Delta\, T_{-m+2} &=&   -(m-2)\, a_m\, \delta(\ux)  \\
\bZ\, ^{s}U_{-m+2} &=& -(m-2)\, a_m\, \uom\, \delta(\ux)
\end{eqnarray*}
and now it also holds that
\begin{eqnarray}
\label{signumlapTm2} \bZ\, T_{-m+2} &=&  -(m-1)\, T_{-m}\, + a_m\,   \delta(\ux)\\
\label{lapsUm2} \Delta\, ^{s}U_{-m+2} &=&   -(m-1)\, ^{s}U_{-m} + a_m\,  \uom\, \delta(\ux)
\end{eqnarray}

\noindent
For $\boxed{\la = -m}$ we already found formulae (\ref{lapTm}) and (\ref{signumlapsUm}), viz.
\begin{eqnarray*}
\Delta\, T_{-m} &=&   2m\, T_{-m-2} + (\onehalf + \frac{1}{m})\, a_m\, \dirac^2\, \delta(\ux)   \\
\bZ\, ^{s}U_{-m} &=& 2m\, ^{s}U_{-m-2} + (\onehalf + \frac{1}{m})\, a_m\, \uom\, \dirac^2\, \delta(\ux) 
\end{eqnarray*}
and now it also holds that
\begin{eqnarray}
\label{signumlapTm} \bZ\, T_{-m} &=&  (m+1)\, T_{-m-2} + \frac{m+2}{2m}\, a_m\, \dirac^2\, \delta(\ux)\\
\label{lapsUm} \Delta\, ^{s}U_{-m} &=&    (m+1)\, ^{s}U_{-m-2} + \frac{m+2}{2m}\, a_m\, \uom\, \dirac^2\, \delta(\ux)
\end{eqnarray}

\noindent
More generally for $\boxed{\la = -m-2\ell}\, , \ell= 0, 1, 2, \ldots$ we already found formulae (\ref{lapTspecial}) and (\ref{signumlapsUspecial}), viz.
\begin{eqnarray*}
\Delta\, T_{-m-2\ell} &=&   (m+2\ell)(2\ell+2)\, T_{-m-2\ell-2} + (-1)^{\ell}\, \frac{(m+4\ell+2)(m+2\ell+2)}{C(m,\ell+1)}\, a_m\, \dirac^{2\ell+2}\, \delta(\ux)  \\
\bZ\, ^{s}U_{-m-2\ell} &=&  (m+2\ell)(2\ell+2)\, ^{s}U_{-m-2\ell-2} + (-1)^{\ell}\, \frac{(m+4\ell+2)(m+2\ell+2)}{C(m,\ell+1)}\, a_m\, \uom\, \dirac^{2\ell+2}\, \delta(\ux) 
\end{eqnarray*}
and now it also holds that
\begin{eqnarray}
\label{signumlapTspecial} \bZ\, T_{-m-2\ell} &=&  (m+2\ell+1)(2\ell+1)\, T_{-m-2\ell-2} + (-1)^{\ell}\, \frac{(m+4\ell+2)(m+2\ell+2)}{C(m,\ell+1)}\, a_m\, \dirac^{2\ell+2}\, \delta(\ux)\\
\label{lapsUspecial} \Delta\, ^{s}U_{-m-2\ell} &=&   (m+2\ell+1)(2\ell+1)\, ^{s}U_{-m-2\ell-2} + (-1)^{\ell}\, \frac{(m+4\ell+2)(m+2\ell+2)}{C(m,\ell+1)}\, a_m\, \uom\, \dirac^{2\ell+2}\, \delta(\ux)
\end{eqnarray}

Now we treat the case of the distributions $U_\la$ and the associated signumdistributions $^{s}T_\la$.\\

\noindent
For $\boxed{\lambda \neq -m+1, -m-1, \ldots}$ we already found formulae (\ref{lapU}) and (\ref{signumlapsT}), viz.
\begin{eqnarray*}
\Delta\, U_\la &=&   (\la-1)(\la + m-1)\, U_{\la-2}  \\
\bZ\, ^{s}T_\la &=&  (\la-1)(\la + m-1)\, ^{s}T_{\la-2} 
\end{eqnarray*}
and now it also holds that
\begin{eqnarray}
\label{signumlapU} \bZ\, U_\la &=&   \lambda\,(\lambda+m-2)\, U_{\lambda-2} \\
\label{lapsT} \Delta\, ^{s}T_\lambda &=&  \lambda\,(\lambda+m-2)\, ^{s}T_{\lambda-2}
\end{eqnarray}

\noindent
For $\boxed{\la = -m+1}$ we already found formulae (\ref{lapUm1}) and (\ref{signumlapsTm1}), viz.
\begin{eqnarray*}
\Delta\, U_{-m+1} &=&   a_m\, \dirac\, \delta(\ux)  \\
\bZ\, ^{s}T_{-m+1} &=& a_m\, \p_r\, \delta(\ux)
\end{eqnarray*}
and now it also holds that
\begin{eqnarray}
\label{signumlapUm1} \bZ\, U_{-m+1} &=&   (m-1)\, U_{-m-1} + a_m\, \dirac\, \delta(\ux) \\
\label{lapsTm1} \Delta\, ^{s}T_{-m+1} &=&  (m-1)\, ^{s}T_{-m-1}  - a_m\, \uom\, \dirac\, \delta(\ux)
\end{eqnarray}

\noindent
More generally for $\boxed{\la = -m-2\ell+1}\, , \ell=1,2,\ldots$ we already found formulae (\ref{lapUspecial}) and (\ref{signumlapsTspecial}), viz.
\begin{eqnarray*}
\Delta\, U_{-m-2\ell+1} &=&   (m+2\ell)(2\ell)\, U_{-m-2\ell-1}  + (-1)^{\ell}\, \frac{m+4\ell}{C(m,\ell)}\, a_m\, \dirac^{2\ell+1}\, \delta(\ux)\\
\bZ\, ^{s}T_{-m-2\ell+1} &=&  (m+2\ell)(2\ell)\, ^{s}T_{-m-2\ell-1}  + (-1)^{\ell+1}\, \frac{m+4\ell}{C(m,\ell)}\, a_m\, \uom\,  \dirac^{2\ell+1}\, \delta(\ux)
\end{eqnarray*}
and now it also holds that
\begin{eqnarray}
\label{signumlapUspecial} \bZ\, U_{-m-2\ell+1} &=& (m+2\ell-1)(2\ell+1)\, U_{-m-2\ell-1} + (-1)^{\ell}\, \frac{m+4\ell}{C(m,\ell)}\, a_m\, \dirac^{2\ell+1}\, \delta(\ux)  \\
\label{lapsTspecial} \Delta\, ^{s}T_{-m-2\ell+1} &=&  (m+2\ell-1)(2\ell+1)\, ^{s}T_{-m-2\ell-1} + (-1)^{\ell+1}\, \frac{m+4\ell}{C(m,\ell)}\, a_m\, \uom\, \dirac^{2\ell+1}\, \delta(\ux)
\end{eqnarray}

\noindent
Note that when putting $\ell=0$ in (\ref{signumlapUspecial}) then (\ref{signumlapUm1}) is obtained.


\newpage
\subsection{The operators $\p_r^2$, $\invr\, \p_r$, $\invrsq\, \Delta^*$ and $\invrsq\, \bZ^*$ }


In general the separate actions of the three components of the Laplace operator expressed in spherical co-ordinates
$$
\Delta = \p_r^2 + (m-1)\, \invr\, \p_r + \invrsq\, \Delta^*
$$
and of the associated signum--operator
$$
\bZ = \p_r^2 + (m-1)\, \invr\, \p_r + \invrsq\, \bZ^*
$$
are not uniquely determined, but are entangled instead. However for the specific (signum)distributions under consideration we expect all actions to be uniquely determined except for the distributions $T_{-m+2}$ and $U_{-m+1}$, and the associated signumdistributions $^{s}U_{-m+2}$ and $^{s}T_{-m+1}$.\\

\noindent
For general $\la$ we find:
\begin{eqnarray*}
\p_r^2\, T_\lambda &=& \lambda(\lambda-1)\, T_{\lambda-2} \quad , \quad \lambda \neq -m+2, -m, -m-2,\ldots\\
\invr\, \p_r\, T_\lambda &=& \lambda\, T_{\lambda-2} \quad , \quad \lambda \neq  -m+2, -m, -m-2,\ldots\\
\invrsq\, \Delta^*\, T_\lambda &=& 0\\
\invrsq\, {\bf {\bf Z}}^*\, T_\lambda &=& -\, (m-1)\, T_{\lambda-2} \quad , \quad \la \neq -m+2
\end{eqnarray*}
leading to expressions for $\Delta\, T_\la$ and $\bZ\, T_\la$ confirming formulae (\ref{lapT}) and (\ref{signumlapT}) respectively, 
and
\begin{eqnarray*}
\p_r^2\, U_\lambda &=& \lambda(\lambda-1)\, U_{\lambda-2} \quad , \quad \lambda \neq -m+1, -m-1, \ldots\\
\invr\, \p_r\, U_\lambda &=& \lambda\, U_{\lambda-2} \quad , \quad \lambda \neq  -m+1, -m-1, \ldots\\
\invrsq\, \Delta^*\, U_\lambda &=& -\, (m-1)\, U_{\lambda-2} \\
\invrsq\, {\bf {\bf Z}}^*\, U_\lambda &=& 0
\end{eqnarray*}
leading to expressions for $\Delta\, U_\la$ and $\bZ\, U_\la$ confirming formulae (\ref{lapU}) and (\ref{signumlapU}) respectively. 
Through signum--pairing we obtain
\begin{eqnarray*}
\p_r^2\, ^{s}T_\lambda &=& \lambda(\lambda-1)\, ^{s}T_{\lambda-2} \quad , \quad \lambda \neq -m+1, -m-1, \ldots\\
\invr\, \p_r\, ^{s}T_\lambda &=& \lambda\, ^{s}T_{\lambda-2} \quad , \quad \lambda \neq  -m+1, -m-1, \ldots\\
\invrsq\, \Delta^*\, ^{s}T_\lambda &=& 0\\
\invrsq\, {\bf {\bf Z}}^*\, ^{s}T_\lambda &=& -\, (m-1)\, ^{s}T_{\lambda-2}
\end{eqnarray*}
leading to expressions for $\Delta\, ^{s}T_\la$ and $\bZ\, ^{s}T_\la$ confirming formulae (\ref{lapsT}) and (\ref{signumlapsT}) respectively, 
and
\begin{eqnarray*}
\p_r^2\, ^{s}U_\lambda &=& \lambda(\lambda-1)\, ^{s}U_{\lambda-2} \quad , \quad \lambda \neq -m+2, -m, \ldots\\
\invr\, \p_r\, ^{s}U_\lambda &=& \lambda\, ^{s}U_{\lambda-2} \quad , \quad \lambda \neq  -m+2, -m, \ldots\\
\invrsq\, \Delta^*\, ^{s}U_\lambda &=& -\, (m-1)\, ^{s}U_{\lambda-2}  \quad , \quad \la \neq -m+2\\
\invrsq\, {\bf {\bf Z}}^*\, ^{s}U_\lambda &=& 0
\end{eqnarray*}
leading to expressions for $\Delta\, ^{s}U_\la$ and $\bZ\, ^{s}U_\la$ confirming formulae (\ref{lapsU}) and (\ref{signumlapsU}) respectively.
\\

\noindent
For $\boxed{\la=-m+2}$ we obtain
\begin{eqnarray*}
\p_r^2\, T_{-m+2} &=& (m-1)(m-2)\, T_{-m} - (m-2)\, a_m\, \delta(\ux)   \\
\invr\, \p_r\, T_{-m+2} &=& \invr\, \left( -\, (m-2)\, ^{s}T_{-m+1} \right) = -\, (m-2)\, (T_{-m} + \delta(\ux)\, c_3) \\
\invrsq\, \Delta^*\, T_{-m+2} &=& 0 \\
\invrsq\, \bZ^*\, T_{-m+2} &=& -\, (m-1)\, T_{-m} - (m-1)\, c_4\, \delta(\ux)
\end{eqnarray*}
these expressions being entangled by (\ref{lapTm2}) and (\ref{signumlapTm2}), viz.
\begin{eqnarray*}
\Delta\, T_{-m+2} &=& -\, (m-2)\, a_m\, \delta(\ux)  \\
\bZ\, T_{-m+2} &=& -\, (m-1)\, T_{-m} + a_m\, \delta(\ux)
\end{eqnarray*}
These entanglement conditions fix the values of the arbitrary constanst $c_3$ and $c_4$ to be $\boxed{c_3 = 0}$ and $\boxed{c_4 = -\, a_m}$, eventually leading to
\begin{eqnarray*}
\p_r^2\, T_{-m+2} &=& (m-1)(m-2)\, T_{-m} - (m-2)\, a_m\, \delta(\ux)   \\
\invr\, \p_r\, T_{-m+2} &=& \invr\, \left( -\, (m-2)\, ^{s}T_{-m+1} \right) = -\, (m-2)\, T_{-m} \\
\invrsq\, \Delta^*\, T_{-m+2} &=& 0 \\
\invrsq\, \bZ^*\, T_{-m+2} &=& -\, (m-1)\, T_{-m} +  (m-1)\, a_m\, \delta(\ux)
\end{eqnarray*}

\noindent
For $\boxed{\la = -m+1}$ we obtain
\begin{eqnarray*}
\p_r^2\, U_{-m+1} &=& m(m-1)\, U_{-m-1} + \frac{2m-1}{m}\, a_m\, \dirac \delta(\ux)   \\
\invr\, \p_r\, U_{-m+1}  &=& -\, (m-1)\, U_{-m-1} -  \frac{1}{m}\, a_m\, \dirac  \delta(\ux)  \\
\invrsq\, \Delta^*\, U_{-m+1} &=& -\, (m-1)\, U_{-m-1}  \\
\invrsq\, \bZ^*\, U_{-m+1} &=& 0
\end{eqnarray*}
confirming the formulae (\ref{lapUm1}) and (\ref{signumlapUm1}), viz.
\begin{eqnarray*}
\Delta\, U_{-m+1} &=&  a_m\, \dirac \delta(\ux) \\
\bZ\, U_{-m+1}&=& (m-1)\, U_{-m-1} + a_m\, \dirac \delta(\ux)
\end{eqnarray*}\\

\noindent
Through the signum--pairs of operators $[\p_r^2, \p_r^2]$, $[\invr\, \p_r , \invr\, \p_r]$, $[\invrsq , \invrsq]$ and $(\Delta^* , \bZ^* )$ we then find

\begin{eqnarray*}
\p_r^2\, ^{s}U_{-m+2} &=& (m-1)(m-2)\, ^{s}U_{-m} - (m-2)\, a_m\, \uom\, \delta(\ux)\\
\invr\, \p_r\, ^{s}U_{-m+2} &=&  -\, (m-2)\, ^{s}U_{-m} \\
\invrsq\, \bZ^*\, U_{-m+2} &=& 0  \\
\invrsq\, \Delta^*\, ^{s}U_{-m+2} &=&  -\, (m-1)\, ^{s}U_{-m} + (m-1)\, a_m\, \uom\, \delta(\ux)
\end{eqnarray*}
confirming (\ref{signumlapsUm2}) and (\ref{lapsUm2}), viz.
\begin{eqnarray*}
\bZ\,  ^{s}U_{-m+2} &=& -\, (m-2)\, a_m\, \uom\, \delta(\ux)  \\
\Delta\, ^{s}U_{-m+2} &=& -\, (m-1)\, ^{s}U_{-m} + a_m\, \uom\, \delta(\ux)  
\end{eqnarray*}
and
\begin{eqnarray*}
\p_r^2\, ^{s}T_{-m+1} &=& m(m-1)\, ^{s}T_{-m-1} - \frac{2m-1}{m}\, a_m\, \uom\, \dirac \delta(\ux)   \\
\invr\, \p_r\, ^{s}T_{-m+1} &=&  -\, (m-1)\, ^{s}T_{-m-1} + \frac{1}{m}\, a_m\, \uom\, \dirac \delta(\ux) \\
\invrsq\, \bZ^*\, ^{s}T_{-m+1} &=&   -\, (m-1)\, ^{s}T_{-m-1}\\
\invrsq\, \Delta^*\, ^{s}T_{-m+1} &=&  0
\end{eqnarray*}
confirming (\ref{signumlapsTm1}) and (\ref{lapsTm1}), viz.
\begin{eqnarray*}
\bZ\, ^{s}T_{-m+1} &=& -\, a_m\, \uom\, \dirac \delta(\ux)  \\
\Delta\, ^{s}T_{-m+1} &=& (m-1)\, ^{s}T_{-m-1} - a_m\, \uom\, \dirac \delta(\ux)
\end{eqnarray*}
\\

\noindent
Additionally we also find
\begin{eqnarray*}
\p_r^2\, T_{-m} &=&   m(m+1)\, T_{-m-2} +  \frac{2m+1}{2m}\, a_m\, \dirac^2\, \delta(\ux)  \\
\invr\, \p_r\, T_{-m} &=& -\, m\, T_{-m-2} - \frac{1}{2m}\, a_m\, \dirac^2\, \delta(\ux) \\
\invrsq\, \Delta^*\, T_{-m} &=& 0 \\
\invrsq\, \bZ^*\, T_{-m} &=& -\, (m-1)\, T_{-m-2}
\end{eqnarray*}
confirming (\ref{lapTm}) and (\ref{signumlapTm}), viz.
\begin{eqnarray*}
\Delta\, T_{-m} &=&  2m\, T_{-m-2} + \frac{m+2}{2m}\, a_m\, \dirac^2\, \delta(\ux) \\
\bZ\, T_{-m} &=& (m+1)\, T_{-m-2} + \frac{m+2}{2m}\, a_m\, \dirac^2\, \delta(\ux)
\end{eqnarray*}
and, via signum--pairing,
\begin{eqnarray*}
\p_r^2\, ^{s}U_{-m} &=&  m(m+1)\, ^{s}U_{-m-2} +  \frac{2m+1}{2m}\, a_m\, \uom\, \dirac^2\, \delta(\ux)  \\
\invr\, \p_r\, ^{s}U_{-m} &=&  -\, m\, ^{s}U_{-m-2} - \frac{1}{2m}\, a_m\, \uom\, \dirac^2\, \delta(\ux) \\
\invrsq\, \bZ^*\, U_{-m} &=& 0 \\
\invrsq\, \Delta^*\, ^{s}U_{-m} &=& -\, (m-1)\, ^{s}U_{-m-2}
\end{eqnarray*}
confirming (\ref{signumlapsUm}) and (\ref{lapsUm}), viz.
\begin{eqnarray*}
\bZ\,  ^{s}U_{-m} &=&   2m\, ^{s}U_{-m-2} + \frac{m+2}{2m}\, a_m\, \uom\, \dirac^2\, \delta(\ux)  \\
\Delta\, ^{s}U_{-m} &=& (m+1)\, ^{s}U_{-m-2} + \frac{m+2}{2m}\, a_m\, \uom\, \dirac^2\, \delta(\ux) 
\end{eqnarray*}\\

More generally it holds that
\begin{eqnarray*}
\p_r^2\, T_{-m-2\ell} &=&  (m+2\ell)(m+2\ell+1)\, T_{-m-2\ell-2}  + (-1)^\ell\, \frac{(m+2\ell+2)(2m+4\ell+1)}{C(m,\ell+1)}\, a_m\, \dirac^{2\ell+2}\, \delta(\ux)\\
\invr\, \p_r\, T_{-m-2\ell} &=&  -\, (m+2\ell)\, T_{-m-2\ell-2} + (-1)^{\ell+1}\, \frac{m+2\ell+2}{C(m,\ell+1)}\, a_m\, \dirac^{2\ell+2}\, \delta(\ux)\\
\invrsq\, \Delta^*\, T_{-m-2\ell} &=& 0 \\
\invrsq\, \bZ^*\, T_{-m-2\ell} &=& -\, (m-1)\, T_{-m-2\ell-2}
\end{eqnarray*}
confirming (\ref{lapTspecial}) and (\ref{signumlapTspecial}), and
\begin{eqnarray*}
\p_r^2\, U_{-m-2\ell+1} &=&  (m+2\ell)(m+2\ell-1)\, U_{-m-2\ell-1} + (-1)^\ell\, \frac{2m+4\ell-1}{C(m,\ell)}\, a_m\, \dirac^{2\ell+1}\delta(\ux)\\
\invr\, \p_r\, U_{-m-2\ell+1} &=& -\, (m+2\ell-1)\, U_{-m-2\ell-1} + (-1)^{\ell+1}\, \frac{1}{C(m,\ell)}\, a_m\, \dirac^{2\ell+1}\delta(\ux)\\
\invrsq\, \Delta^*\, U_{-m-2\ell+1} &=& -\, (m-1)\, U_{-m+2\ell-1} \\
\invrsq\, \bZ^*\, U_{-m-2\ell+1} &=& 0
\end{eqnarray*}
confirming (\ref{lapUspecial}) and (\ref{signumlapUspecial}).\\

\noindent
Via signum--pairing it also holds that
\begin{eqnarray*}
\p_r^2\, ^{s}U_{-m-2\ell} &=& (m+2\ell)(m+2\ell+1)\, ^{s}U_{-m-2\ell-2}  + (-1)^\ell\, \frac{(m+2\ell+2)(2m+4\ell+1)}{C(m,\ell+1)}\, a_m\, \uom\, \dirac^{2\ell+2}\, \delta(\ux)\\
\invr\, \p_r\, ^{s}U_{-m-2\ell} &=&  -\, (m+2\ell)\, ^{s}U_{-m-2\ell-2} + (-1)^{\ell+1}\, \frac{m+2\ell+2}{C(m,\ell+1)}\, a_m\, \uom\, \dirac^{2\ell+2}\, \delta(\ux)\\
\invrsq\, \bZ^*\, ^{s}U_{-m-2\ell} &=& 0 \\
\invrsq\, \Delta^*\, ^{s}U_{-m-2\ell} &=& -\, (m-1)\, ^{s}U_{-m-2\ell-2}
\end{eqnarray*}
confirming (\ref{signumlapsUspecial}) and (\ref{lapsUspecial}), and
\begin{eqnarray*}
\p_r^2\, ^{s}T_{-m-2\ell+1} &=&  (m+2\ell)(m+2\ell-1)\, ^{s}T_{-m-2\ell-1} + (-1)^{\ell+1}\, \frac{2m+4\ell-1}{C(m,\ell)}\, a_m\, \uom\, \dirac^{2\ell+1}\delta(\ux)\\
\invr\, \p_r\, ^{s}T_{-m-2\ell+1} &=& -\, (m+2\ell-1)\, ^{s}T_{-m-2\ell-1} + (-1)^{\ell}\, \frac{1}{C(m,\ell)}\, a_m\, \uom\, \dirac^{2\ell+1}\delta(\ux)\\
\invrsq\, \Delta^*\, ^{s}T{-m-2\ell+1} &=&  0\\
\invrsq\, \bZ^*\, ^{s}T{-m-2\ell+1} &=& -\, (m-1)\, ^{s}T_{-m+2\ell-1}
\end{eqnarray*}
confirming (\ref{lapsTspecial}) and (\ref{signumlapsTspecial}).

\begin{remark}
{\rm
To give an idea how compuatations with signumdistributions are performed in practice, we will directly prove formula (\ref{diracsTm1}), viz.
$$
\dirac\, ^{s}T_{-m+1} = -\, (m-1)\, ^{s}U_{-m} + a_m\, \uom\, \delta(\ux)
$$
Let $\uom\, \varphi \in \Omega(\mR^m; \mR^m)$ be a signumdistributional test function, then
\begin{eqnarray*}
\langle \, \uom\, \varphi \, , \, \dirac\, ^{s}T_{-m+1} \, \rangle &=& -\, \langle \, (\uom\, \varphi) \, \dirac \, , \, ^{s}T_{-m+1} \, \rangle \\
&=& -\, \langle \, (1-m)\, \invr\, \varphi + \uom\, \dirac(\varphi) \, , \, ^{s}T_{-m+1} \, \rangle \\
&=& (m-1)\, \langle \,  \invr\, \varphi \, , \, ^{s}T_{-m+1} \, \rangle -\, \langle \,  \uom\, \dirac(\varphi) \, , \, ^{s}T_{-m+1} \, \rangle \\
\end{eqnarray*}
The first term at the right--hand side equals
$$
(m-1)\, \langle\, \invux\, \varphi \, , \, U_{-m+1} \, \rangle = (m-1)\,  \langle\,  \varphi \, , \, \invux\, U_{-m+1} \, \rangle =
(m-1)\,  \langle\,  \varphi \, , \, T_{-m} \, \rangle = -(m-1)\,  \, \langle \, \uom\, \varphi \, , \, ^{s}U_{-m} \, \rangle
$$
while the second term at the right--hand side equals
\begin{eqnarray*}
-\, a_m\, \langle \,  \Sigma^{(0)}\left[ \uom\, \dirac(\varphi)\right]  \, , \, {\rm Fp}\, r_+^0 \, \rangle_r &=& 
-\, a_m\, \langle \,  \Sigma^{(1)}\left[ \dirac(\varphi)\right]  \, , \, {\rm Fp}\, r_+^0 \, \rangle_r\\
&=& -\,  \langle \,   \dirac(\varphi) \, , \, U_{-m+1} \, \rangle\\
&=&  \langle \,   \varphi \, , \, \dirac\, U_{-m+1} \, \rangle\\
&=&  -\, a_m\, \langle \,   \varphi \, , \, \delta(\ux) \, \rangle\\
&=&   a_m\, \langle \,  \uom\,  \varphi \, , \, \uom\, \delta(\ux) \, \rangle
\end{eqnarray*}
and the formula follows readily.
}
\end{remark}


\newpage
\section{Cartesian derivatives of $T_\la$ and $U_\la$}
\label{cartderiv}


In this section we compute the first and second order partial derivatives of the distributions $T_\la$ and $U_\la$ and of their associated signumdistributions $^{s}T_\la$ and $^{s}U_\la$. Recall that the distributions $T_\la$ are scalar valued, while the distributions $U_\la$ are (Clifford) vector valued with scalar components denoted by $U_\la^j, j=1,\ldots,m$. As  it holds,  for all $\la \in \mC$, that $\ux\, T_\la = U_{\la+1}$ , it follows that these components are given by
$$
U_\la^j = x_j\, T_{\la-1} \quad , \quad j=1,\ldots,m
$$

Quite naturally the action of the cartesian derivatives $\p_{x_j}$  on all of the distributions $T_\la$ and $U_\la$ will be uniquely determined. It follows that the action of the operators $d_j$ on the associated signumdistributions $^{s}U_\la$ and $^{s}T_\la$ will be uniquely determined too.\\
Seen the definition of the operators $d_j, j=1,\ldots,m$, viz.
$$
d_j = \invux\, e_j + (-\, \frac{x_j}{r^2}) + \p_{x_j} = -\, e_j \invux - (-\, \frac{x_j}{r^2}) + \p_{x_j}
$$
their action on the distributions $T_\la$ and $U_\la$ will be uniquely determined except for the distribution $U_{-m+1}$, and, by signum--pairing, the cartesian derivatives of the associated signumdistributions $^{s}U_\la$ and $^{s}T_\la$ will be uniquely determined too, except for $^{s}T_{-m+1}$. 

The actions of the second order cartesian derivatives $\p_{x_k} \p_{x_j}$ and $\p^2_{x_j}$ on the  distributions $T_\la$ and $U_\la$ and of the operators $d_k d_j$ and $d_j^2$ on the associated signumdistributions $^{s}U_\la$ and $^{s}T_\la$, will be uniquely determined. Seen the expressions
$$
d_k d_j = \invrsq\, (x_k + e_k\, \ux)\, \p_{x_j} - \invrsq\, \ux\, e_j\, \p_{x_k} + \frac{1}{r^4}\, x_k x_j  +  \frac{1}{r^4}\, (\ux\, x_k\, e_j - e_k\, \ux\, x_k) + \p_{x_k}\p_{x_j}
$$
and
$$
d_j^2 = -\, \invrsq - 4\, e_j\, x_j\, \frac{\ux}{r^4} - 3\, \frac{x_j^2}{r^4} + \p_{x_j}^2
$$
the actions of the operators $d_k d_j$ and $d_j^2$ on the  distributions $T_\la$ and $U_\la$ will be uniquely determined except for the distributions $T_{-m+2}$ and $U_{-m+1}$. Through signum--pairing the actions of the operators $\p_{x_k} \p_{x_j}$ and $\p^2_{x_j}$ on the associated signumdistributions $^{s}U_\la$ and $^{s}T_\la$ will be uniquely determined too except for $^{s}U_{-m+2}$ and $^{s}T_{-m+1}$.\\
Notice that
$$
\Sigma_{j<k}\, e_j e_k\, (d_jd_k -d_kd_j) = -\, \invrsq\, \Gamma
$$
and
$$
\Sigma_{j=1}^m\, d_j^2 = \Delta - (m-1)\, \invrsq
$$


\newpage
\subsection{First order partial derivatives of $T_\la$}


In general, i.e. for $\la \neq -m, -m-2,\ldots$, it holds that
$$
\boxed{\p_{x_j}\, T_\la = \la\, x_j\, T_{\la-2}}
$$
verifying (\ref{diracTlambda}) by
$$
\dirac\, T_\la = \la\, \ux\, T_{\la-2} = \la\, U_{\la-1}
$$
and verifying (\ref{eulerT}) by
$$
\mE\, T_\la = \Sigma_{j=1}^m\, x_j\, \la\, x_j\, T_{\la-2} = -\, \ux^2\, T_{\la-2} = \la\, T_\la
$$
Through the signumpair of operators $(\p_{x_j} , d_j)$ it follows that
$$
d_{j}\, ^{s}U_\la = \la\, x_j\, ^{s}U_{\la-2}
$$
Moreover a straightforward computation shows that
$$
\boxed{d_{j}\, T_\la
= -\, U_{\la-1}\, e_j + (\la-1)\, x_j\, T_{\la-2} = e_j\, U_{\la-1} + (\la+1)\, x_j\, T_{\la-2}}
$$
whence,   through the signumpair of operators $[d_j , \p_{x_j} ]$, 
$$
\p_{x_j}\, ^{s}U_\la =\,   ^{s}T_{\la-1}\, e_j + (\la-1)\, x_j\, ^{s}U_{\la-2}
$$
As a verification we may now compute $\dirac\, ^{s}U_\la$ and we obtain
\begin{eqnarray*}
\dirac\, ^{s}U_\la &=& (-m)\, ^{s}T_{\la-1} + (\la-1)\, \ux\, ^{s}U_{\la-2} = -\, (\la+m-1)\, ^{s}T_{\la-1}
\end{eqnarray*}
confirming (\ref{diracsU}).\\

\noindent
For the singular values $\la = -m, -m-2, \ldots$ it holds that
$$
\boxed{\p_{x_j}\, T_{-m-2\ell} = -\, (m+2\ell)\, x_j\, T_{-m-2\ell-2} - \frac{a_m}{C(m,\ell)}\, \p_{x_j}\,  \Delta^\ell\,  \delta(\ux)}
$$
with, recall,
$$
C(m,\ell) = 2^{2\ell+1}\, \ell!\, \frac{\Gamma(\frac{m}{2}+\ell+1)}{\Gamma(\frac{m}{2})} =  2^{\ell}\, \ell!\, m(m+2)(m+4)\cdots(m+2\ell)
$$
verifying (\ref{diracTspecial}) by
\begin{eqnarray*}
\dirac\, T_{-m-2\ell} &=&  -\, (m+2\ell)\, \ux\, T_{-m-2\ell-2} - \frac{a_m}{C(m,\ell)}\, \dirac\,  \Delta^\ell\,  \delta(\ux)\\
&=& -\, (m+2\ell)\, U_{-m-2\ell-1} + (-1)^{\ell+1}\, \frac{a_m}{C(m,\ell)}\, \dirac^{2\ell+1}\, \delta(\ux)
\end{eqnarray*}
and verifying (\ref{eulerTspecial}) by
\begin{align*}
\mE\, T_{-m-2\ell+1} &=   -\, (m+2\ell)\, r^2\, T_{-m-2\ell-2} - \frac{a_m}{C(m,\ell)}\, \mE\,  \Delta^\ell\,  \delta(\ux) \\
&=  -\, (m+2\ell)\,  T_{-m-2\ell} + (-1)^\ell \frac{a_m}{C(m,\ell)}\, (m+2\ell)\,  \dirac^{2\ell}\,  \delta(\ux) 
\end{align*}

\noindent
In particular, for $\ell=0$, we have
$$
\p_{x_j}\, T_{-m} = -\, m\, x_j\, T_{-m-2} - \frac{1}{m}\, a_m\, \p_{x_j}\delta(\ux)
$$
verifying (\ref{diracTm}) by
\begin{eqnarray*}
\dirac\, T_{-m} &=&  -\, m\, \ux\, T_{-m-2} - \frac{1}{m}\, a_m\,  \dirac\,  \delta(\ux) = -\, m\, U_{-m-1} - \frac{1}{m}\, a_m\,  \dirac\,  \delta(\ux)
\end{eqnarray*}
and verifying (\ref{eulerTm}) by
\begin{align*}
\mE\, T_{-m} &= -\, m\,  T_{-m} +  a_m\, \delta(\ux)
\end{align*}

\noindent
Hence
\begin{eqnarray*}
d_j\, ^{s}U_{-m-2\ell} &=& -\, (m+2\ell)\, x_j\, ^{s}U_{-m-2\ell-2} - \frac{a_m}{C(m,\ell)}\, d_j\, \bZ^\ell\,  (\uom\, \delta(\ux))\\
&=& -\, (m+2\ell)\, x_j\, ^{s}U_{-m-2\ell-2} - \frac{a_m}{C(m,\ell)}\, \uom\, \Delta^{\ell}\, \p_{x_j}\, \delta(\ux)
\end {eqnarray*}
and, in particular, for $\ell=0$,
\begin{eqnarray*}
d_j\, ^{s}U_{-m} &=& -\, m\, x_j\, ^{s}U_{-m-2} - \frac{1}{m}\, a_m\,  d_j\, (\uom \delta(\ux))\\
&=& -\, m\, x_j\, ^{s}U_{-m-2} - \frac{1}{m}\, a_m\,  \uom\, \p_{x_j}\, \delta(\ux)
\end{eqnarray*}
Moreover a direct computation shows that
$$
\boxed{d_{j}\, T_{-m-2\ell} 
= -\, U_{-m-2\ell-1}\, e_j -\, (m+2\ell+1)\, x_j\, T_{-m-2\ell-2} - \frac{a_m}{C(m,\ell)}\, \p_{x_j}\, \Delta^\ell\,  \delta(\ux)}
$$
whence
\begin{eqnarray*}
\p_{x_j}\, ^{s}U_{-m-2\ell} &=&\, ^{s}T_{-m-2\ell-1}\, e_j -\, (m+2\ell+1)\, x_j\, ^{s}U_{-m-2\ell-2} - \frac{a_m}{C(m,\ell)}\, d_j\, \bZ^\ell\,  \uom\,  \delta(\ux)\\
&=&\, ^{s}T_{-m-2\ell-1}\, e_j -\, (m+2\ell+1)\, x_j\, ^{s}U_{-m-2\ell-2} - \frac{a_m}{C(m,\ell)}\, \uom\, \p_{x_j}\, \Delta^{\ell}\,  \delta(\ux)
\end{eqnarray*}
and, in particular for $\ell=0$,
$$
d_j\, T_{-m} =  -\, U_{-m-1}\, e_j -\, (m+1)\, x_j\, T_{-m-2} - \frac{1}{m}\, a_m\, \p_{x_j}\, \delta(\ux)
$$
and so
\begin{eqnarray*}
\p_{x_j}\, ^{s}U_{-m} &=&\, ^{s}T_{-m-1}\, e_j -\, (m+1)\, x_j\, ^{s}U_{-m-2} - \frac{1}{m}\, a_m\, d_j\, \uom\, \delta(\ux)\\
&=&\, ^{s}T_{-m-1}\, e_j -\, (m+1)\, x_j\, ^{s}U_{-m-2} - \frac{1}{m}\, a_m\,  \uom\, \p_{x_j} \delta(\ux)
\end{eqnarray*}

\noindent
As a verification of formula  (\ref{diracsUspecial}) we now may compute $\dirac\, ^{s}U_{-m-2\ell}$ and we obtain
\begin{eqnarray*}
\dirac\, ^{s}U_{-m-2\ell} &=& (-m)\, ^{s}T_{-m-2\ell-1} - (m+2\ell+1)\, \ux\, ^{s}U_{-m-2\ell-2} - \frac{a_m}{C(m,\ell)}\, \Sigma_{\ell=1}^m\, e_j\, \uom\, \p_{x_j}\Delta^\ell\, \delta(\ux)\\
&=& (2\ell+1)\, ^{s}T_{-m-2\ell-1} + (-1)^{\ell+1}\, \frac{a_m}{C(m,\ell)}\, \uom\, \dirac^{2\ell+1}\delta(\ux)
\end{eqnarray*}
since, due to the radial character of $\Delta^\ell \delta(\ux)$, it holds that
\begin{eqnarray*}
\Sigma_{\ell=1}^m\, e_j\, \uom\, \p_{x_j}\Delta^\ell\, \delta(\ux) &=& \Sigma_{\ell=1}^m\, e_j\, \uom\, \omega_j\, \p_{r}\, (-1)^{\ell}\, \dirac^{2\ell}\, \delta(\ux)\\
&=& (-1)^{\ell}\, \Sigma_{\ell=1}^m\, e_j\,  \omega_j\, (\uom\, \p_{r})\,  \dirac^{2\ell}\, \delta(\ux)\\
&=& (-1)^{\ell}\, \uom\, \dirac^{2\ell+1}\delta(\ux)
\end{eqnarray*}


\newpage
\subsection{First order partial derivatives of $U_\la$}


In general, i.e. for $\la \neq -m+1, -m-1, -m-3, \ldots$, it holds that
$$
\boxed{\p_{x_j}\, U_\la = T_{\la-1}\, e_j + (\la-1)\, x_j\, U_{\la-2}}
$$
verifying (\ref{diracUlambda}) by
$$
\dirac\, U_\la  = -\, m\, T_{\la-1} + (\la-1)\, \ux\, U_{\la-2} = -\, (\la+m-1)\, T_{\la-1}
$$
verifying and (\ref{eulerU}) by
$$
\mE\, U_\la = \ux\, T_{\la-1} + (\la-1)\, U_\la = \la\, U_\la
$$

\noindent
Through the signumpair of operators $(\p_{x_j} , d_j)$ it follows that
$$
d_{j}\, ^{s}T_\la = -\, ^{s}U_{\la-1}\, e_j+ (\la-1)\, x_j\ ^{s}T_{\la-2}
$$
Moreover a straightforward computation shows  that
$$
\boxed{d_{j}\, U_\la 
=  \la\, x_j\, U_{\la-2}}
$$
whence,   through the signumpair of operators $[d_j , \p_{x_j} ]$, 
$$
\p_{x_j}\, ^{s}T_\la = \la\, x_j\, ^{s}T_{\la-2}
$$
We may now compute $\dirac\, ^{s}T_\la$ and we obtain
$$
\dirac\, ^{s}T_\la = \la\, \ux\, ^{s}T_{\la-2} = \la\, ^{s}U_{\la-1}
$$
confirming (\ref{diracsT}).\\

\noindent
For the exceptional value $\la=-m+1$ it holds that
$$
\boxed{\p_{x_j}\, U_{-m+1} 
= T_{-m}\, e_j -\, m\, x_j\, U_{-m-1} + \frac{1}{m}\, a_m\, \delta(\ux)\, e_j}
$$
verifying (\ref{diracUm11}) by
$$
\dirac\, U_{-m+1} = -\, m\, T_{-m} - m\, \ux\, U_{-m-1} - a_m\, \delta(\ux) = - a_m\, \delta(\ux) 
$$
and verifying (\ref{eulerUm1}) by
$$
\mE\, U_{-m+1} = \ux\, T_{-m} -\, m\, U_{-m+1} + \frac{1}{m}\, a_m\, \ux\, \delta(\ux) = (-m+1)\, U_{-m+1}
$$

\noindent
Hence
\begin{eqnarray*}
d_j\, ^{s}T_{-m+1} &=&  -\,  ^{s}U_{-m}\, e_j   - m\, x_j\, ^{s}T_{-m-1} +  \frac{1}{m}\, a_m\,     x_j\, \uD\, \uom\, \delta(\ux) \\
&=&  -\,  ^{s}U_{-m}\, e_j   - m\, x_j\, ^{s}T_{-m-1} -  \frac{1}{m}\, a_m\, \uom\, \delta(\ux)\, e_j
\end{eqnarray*}

\noindent
For the action of the operator $d_j$ we expect a non--uniquely determined result, and, indeed, a direct computation shows that
$$
d_j\, U_{-m+1} = -\, (m-1)\, x_j\, U_{-m-1} + c^*\, x_j\, \dirac \delta(\ux)
$$
where we expect the constant $c^*$, to be determined by the formulae (\ref{uDUm1}) and (\ref{diracsTm1}).\\
Indeed we have, through signum--pairing, that
$$
\p_{x_j}\, ^{s}T_{-m+1} = -\, (m-1)\, x_j\, ^{s}T_{-m-1} - c^*\, x_j\, \uom\, \dirac \delta(\ux) = -\, (m-1)\, x_j\, ^{s}T_{-m-1} + c^*\, x_j\, \p_r\,  \delta(\ux)
$$
leading to
\begin{eqnarray*}
\dirac\, ^{s}T_{-m+1} &=& -\, (m-1)\, \ux\, ^{s}T_{-m-1} +  c^*\,  \ux\, \p_r\, \delta(\ux)\\
&=& -\, (m-1)\, ^{s}U_{-m} - m\, c^*\, \uom\, \delta(\ux)
\end{eqnarray*}

\noindent
However we already found, see (\ref{diracsTm1}), that
$$
\dirac\, ^{s}T_{-m+1} = -\, (m-1)\, ^{s}U_{-m} + a_m\,  \uom\, \delta(\ux)
$$
from which it follows that
$
c^* =  -\, \frac{1}{m}\,  a_m
$,
whence
$$
\boxed{d_j\, U_{-m+1} = -\, (m-1)\, x_j\, U_{-m-1} - \frac{1}{m}\, a_m\, x_j\, \dirac \delta(\ux)}
$$
and
$$
\boxed{\p_{x_j}\, ^{s}T_{-m+1} = -\, (m-1)\, x_j\, ^{s}T_{-m-1} + \frac{1}{m}\, a_m\, x_j\, \uom\, \dirac \delta(\ux) } 
$$
\\

\noindent
For the singular values  $\la = -m-1, -m-3, \ldots$ it holds that
$$
\boxed{\p_{x_j}\, U_{-m-2\ell+1} = T_{-m-2\ell}\, e_j   -\, (m+2\ell)\, x_j\, U_{-m-2\ell-1} + (-1)^{\ell-1}\, \frac{a_m}{C(m,\ell)}\, x_j\, \dirac^{2\ell+1}\, \delta(\ux)}
$$
Note that by putting $\ell=0$ in the above expression for $\p_{x_j}\, U_{-m-2\ell+1}$, the formula for $\p_{x_j}\, U_{-m+1}$ is recovered, and that formula (\ref{diracUspecial}) is verified through
\begin{eqnarray*}
\dirac\, U_{-m-2\ell+1} &=& -\, m\, T_{-m-2\ell}   + (m+2\ell)\, T_{-m-2\ell} + (-1)^{\ell-1}\, \frac{a_m}{C(m,\ell)}\, \ux\, \dirac^{2\ell+1}\, \delta(\ux)\\
&=& 2\, \ell\, T_{-m-2\ell} + (-1)^{\ell-1}\, (m+2\ell)\, \frac{a_m}{C(m,\ell)}\, \dirac^{2\ell}\, \delta(\ux)
\end{eqnarray*}
Also (\ref{eulerUspecial}) is verified by
\begin{align*}
\mE\, U_{-m+2\ell+1} &= \ux\, T_{-m-2\ell}   -\, (m+2\ell)\, r^2\, U_{-m-2\ell-1} + (-1)^{\ell}\, \frac{a_m}{C(m,\ell)}\, \ux^2\, \dirac^{2\ell+1}\, \delta(\ux)   \\
&= -\, (m+2\ell-1)\, U_{-m-2\ell+1} + (-1)^{\ell}\, \frac{a_m}{C(m,\ell)}\, (m+2\ell)(2\ell)\, \dirac^{2\ell-1}\, \delta(\ux) \\
&= -\, (m+2\ell-1)\, U_{-m-2\ell+1} + (-1)^{\ell}\, \frac{a_m}{C(m,\ell-1)}\,  \dirac^{2\ell-1}\, \delta(\ux)
\end{align*}

\noindent
Through signum--pairing we obtain
\begin{eqnarray*}
d_j\, ^{s}T_{-m-2\ell+1} &=& \, -\,  ^{s}U_{-m-2\ell}\, e_j   - (m+2\ell)\, x_j\, ^{s}T_{-m-2\ell-1} + (-1)^{\ell}\, \frac{a_m}{C(m,\ell)}\, x_j\, \uD^{2\ell+1}\, \uom\, \delta(\ux)\\
 &=& \, -\,  ^{s}U_{-m-2\ell}\, e_j   - (m+2\ell)\, x_j\, ^{s}T_{-m-2\ell-1} + (-1)^{\ell}\, \frac{a_m}{C(m,\ell)}\, \uom\, x_j\, \dirac^{2\ell+1}\,  \delta(\ux)
\end{eqnarray*}

\noindent
Moreover a direct computation leads to
\begin{eqnarray*}
\boxed{d_{j}\, U_{-m-2\ell+1} =  -\, (m+2\ell-1)\, x_j\, U_{-m-2\ell-1} + (-1)^{l-1}\, \frac{a_m}{C(m,\ell)}\, x_j\, \dirac^{2\ell+1}\,  \delta(\ux)}
\end{eqnarray*}
Note that by putting $\ell=0$ in the above expression for $d_{j}\, U_{-m-2\ell+1} $, the formula for $d_{j}\, U_{-m+1} $ is recovered. Through signum--pairing we  obtain
\begin{eqnarray*}
\p_{x_j}\, ^{s}T_{-m-2\ell+1} &=& \, -\, (m+2\ell-1)\, x_j\, ^{s}T_{-m-2\ell-1} + (-1)^{l}\, \frac{a_m}{C(m,\ell)}\, x_j\, \uD^{2\ell+1}\,  \uom\, \delta(\ux)\\
&=& \, -\, (m+2\ell-1)\, x_j\, ^{s}T_{-m-2\ell-1} + (-1)^{l}\, \frac{a_m}{C(m,\ell)}\, \uom\, x_j\, \dirac^{2\ell+1}\,  \delta(\ux)
\end{eqnarray*}
We may now verify formula (\ref{diracsTsing}) by
\begin{eqnarray*}
\dirac\, ^{s}T_{-m-2\ell+1}  &=& -\, (m+2\ell-1)\, \ux\, ^{s}T_{-m-2\ell-1} + (-1)^\ell\, \frac{a_m}{C(m,\ell)}\, \Sigma_{j=1}^m\, e_j\, \uom\, x_j\, \dirac^{2\ell+1}\delta(\ux)
\end{eqnarray*}
which indeed reduces to 
\begin{eqnarray*}
\dirac\, ^{s}T_{-m-2\ell+1}  &=& -\, (m+2\ell-1)\,  ^{s}U_{-m-2\ell} + (-1)^\ell\, \frac{a_m}{C(m,\ell)}\, (m+2\ell)\, \uom\, \dirac^{2\ell}\delta(\ux)
\end{eqnarray*}
since
\begin{eqnarray*}
\Sigma_{j=1}^\ell\, e_j\, \uom\, x_j\, \dirac^{2\ell+1}\delta(\ux) &=& \Sigma_{j=1}^m\, e_j\, \uom\, r\, \omega_j\, \dirac^{2\ell+1}\delta(\ux)\\
&=& -\, r\, \dirac^{2\ell+1}\delta(\ux) = \uom\, (r \p_{r})\, \dirac^{2\ell}\delta(\ux)\\
&=& (m+2\ell)\, \uom\, \dirac^{2\ell}\delta(\ux)
\end{eqnarray*}


\newpage
\subsection{Second order partial derivatives of $T_\la$}


First we assume that $\la \neq -m+2, -m, -m-2, \ldots$, and we find
\begin{eqnarray*}
\p_{x_k}\p_{x_j}\, T_\la &=& \la\, (\la-2)\, x_k\, x_j\, T_{\la-4} \quad , \quad j \neq k\\[2mm]
\p_{x_j}^2\, T_\la &=& \la\, (\la-2)\, x_j^2\, T_{\la-4} + \la\, T_{\la-2}
\end{eqnarray*}
confirming (\ref{lapT}) through
$$
\Delta\, T_\la = \la\, (\la-2)\, r^2\, T_{\la-4} + m\, \la\, T_{\la-2} = \la\, (\la + m -2  )\, T_{\la-2}
$$

\noindent
Through the signumpartners of $\p_{x_k}$, $\p_{x_j}$ and $T_\la$, it follows that
\begin{eqnarray*}
d_k\, d_j\, ^{s}U_\la &=&  \la\, (\la-2)\, x_k\, x_j\, ^{s}U_{\la-4} \quad , \quad j \neq k\\[2mm]
d_j^2\, ^{s}U_\la &=&  \la\, (\la-2)\,  x_j^2\, ^{s}U_{\la-4} + \la ^{s}U_{\la-2}
\end{eqnarray*}
A direct computation shows that 
\begin{eqnarray*}
d_k d_j\, T_\la &=&  (\la-1)(\la-3)\, x_k\, x_j\, T_{\la-4} - (\la-1)\, U_{\la-3}\, (x_k\, e_j + x_j\, e_k) \quad ,  \quad j \neq k
\end{eqnarray*}
and
\begin{eqnarray*}
 d_j^2\, T_\la &=&  (\la-1)(\la-3)\, x_j^2\, T_{\la-4} - 2(\la-1)\, U_{\la-3}\, x_j\, e_j + (\la-1)\, T_{\la-2}
\end{eqnarray*}
from which it follows at once that
\begin{eqnarray*}
\p_{x_k}\p_{x_j}\, ^{s}U_\la &=&  (\la-1)(\la-3)\, x_k\, x_j\, ^{s}U_{\la-4} + (\la-1)\, ^{s}T_{\la-3}\, (x_k\, e_j + x_j\, e_k) \quad , \quad j \neq k\\[2mm]
\p_{x_j}^2\, ^{s}U_\la &=& (\la-1)(\la-3)\, x_j^2\, ^{s}U_{\la-4} + 2(\la-1)\, ^{s}T_{\la-3}\, x_j\, e_j + (\la-1)\, ^{s}U_{\la-2}
\end{eqnarray*}
confirming (\ref{lapsU}) through
$$
\Delta\, ^{s}U_\la = (\la-1)(\la-3)\, r^2\, ^{s}U_{\la-4} + 2(\la-1)\, ^{s}T_{\la-3}\, \ux + m\, (\la-1)\, ^{s}U_{\la-2} = (\la-1)\, (\la+m-1)\,  ^{s}U_{\la-2} 
$$
\\

\noindent
For the singular values of the parameter $\la$, viz. $\la =  -m+2, -m, -m-2, \ldots$, specific computations are necessary.\\

\noindent
For $\la = -m+2$ we obtain
\begin{eqnarray*}
\p_{x_k}\p_{x_j}\, T_{-m+2} &=& m\, (m-2)\, x_k\, x_j\, T_{-m-2} \quad , \quad j \neq k\\
\p_{x_j}^2\, T_{-m+2} &=& m\, (m-2)\, x_j^2\, T_{-m-2} - (m-2)\, T_{-m} - \frac{m-2}{m}\, a_m\, \delta(\ux)
\end{eqnarray*}
confirming (\ref{lapTm2}) through
$$
\Delta\, T_{-m+2} = m\, (m-2)\, r^2\, T_{-m-2} - m (m-2)\, T_{-m} - (m-2)\, a_m\, \delta(\ux)
= - (m-2)\, a_m\, \delta(\ux)
$$
It follows that
\begin{eqnarray*}
d_k d_j\, ^{s}U_{-m+2} &=& m\, (m-2)\, x_k\, x_j\, ^{s}U_{-m-2} \quad , \quad j \neq k\\
d_j^2\, ^{s}U_{-m+2} &=& m\, (m-2)\, x_j^2\, ^{s}U_{-m-2} - (m-2)\, ^{s}U_{-m} - \frac{m-2}{m}\, a_m\, \uom\, \delta(\ux)
\end{eqnarray*}
A direct computation shows that
\begin{eqnarray*}
d_k d_j\, T_{-m+2} &=& (m-1)(m+1)\, x_j\, x_k\, T_{-m-2} + (m-1)\, U_{-m-1}\, (x_k\, e_j + x_j\, e_k) - \alpha_k\, \delta(\ux)\, e_k\, e_j \quad , \quad j \neq k\\
d_j^2\, T_{-m+2} &=& (m-1)(m+1)\, x_j^2\, T_{-m-2} + 2\,(m-1)\, U_{-m-1}\, x_j\, e_j  + (-m+1)\, T_{-m}  + (\alpha_j -\frac{m-1}{m}\, a_m)\, \delta(\ux)
\end{eqnarray*}
whence also
\begin{eqnarray*}
\p_{x_k}\, \p_{x_j}\, ^{s}U_{-m+2} &=& (m-1)(m+1)\, x_j\, x_k\, ^{s}U_{-m-2} + (-m+1)\, ^{s}T_{-m-1}\, (x_k\, e_j + x_j\, e_k) - \alpha_k\, \uom\, \delta(\ux)\, e_k\, e_j  \quad , \quad j \neq k\\
\p_{x_j}^2\, ^{s}U_{-m+2} &=&  (m-1)(m+1)\, x_j^2\, ^{s}U_{-m-2} + 2\,(-m+1)\, ^{s}T_{-m-1}\, x_j\, e_j  + (-m+1)\, ^{s}U_{-m}  + (\alpha_j -\frac{m-1}{m}\, a_m)\, \uom\, \delta(\ux)
\end{eqnarray*}
leading to
$$
\Delta\, ^{s}U_{-m+2} = -\, (m-1)\, ^{s}U_{-m} + ((\Sigma\, \alpha_j) - (m-1)\, a_m)\uom\, \delta(\ux)
$$
But we already found, see (\ref{lapsUm2}), that
$$
\Delta\, ^{s}U_{-m+2} = -\, (m-1)\, ^{s}U_{-m}  + a_m\, \uom\, \delta(\ux)
$$
whence all $\alpha_j, j=1,\ldots,m$ must be equal to $a_m$, which eventually leads to\\[-2mm]
\begin{eqnarray*}
d_k d_j\, T_{-m+2} &=& (m-1)(m+1)\, x_j\, x_k\, T_{-m-2} + (m-1)\, U_{-m-1}\, (x_k\, e_j + x_j\, e_k) - a_m\, \delta(\ux)\, e_k\, e_j \quad , \quad j \neq k\\
d_j^2\, T_{-m+2} &=& (m-1)(m+1)\, x_j^2\, T_{-m-2} + 2\,(m-1)\, U_{-m-1}\, x_j\, e_j  + (-m+1)\, T_{-m}  + \frac{1}{m}\, a_m\, \delta(\ux)
\end{eqnarray*}
and
\begin{eqnarray*}
\p_{x_k}\, \p_{x_j}\, ^{s}U_{-m+2} &=& (m-1)(m+1)\, x_j\, x_k\, ^{s}U_{-m-2} + (-m+1)\, ^{s}T_{-m-1}\, (x_k\, e_j + x_j\, e_k) - a_m\, \uom\, \delta(\ux)\, e_k\, e_j  \quad , \quad j \neq k\\
\p_{x_j}^2\, ^{s}U_{-m+2} &=&  (m-1)(m+1)\, x_j^2\, ^{s}U_{-m-2} + 2\,(-m+1)\, ^{s}T_{-m-1}\, x_j\, e_j  + (-m+1)\, ^{s}U_{-m}  + \frac{1}{m}\, a_m\, \uom\, \delta(\ux)
\end{eqnarray*}
\\

\noindent
For $\la = -m-2\ell, \ell = 0,1,2,\ldots$ we obtain
\begin{eqnarray*}
\p_{x_k}\p_{x_j}\, T_{-m-2\ell} &=& (m+2\ell)(m+2\ell+2)\, x_k\, x_j\, T_{-m-2\ell-4}  -\,\frac{2m+4\ell+2}{ (m+2\ell+2)C(m,\ell)}\,  a_m\, \p_{x_k}\p_{x_j}\, \Delta^\ell\, \delta(\ux)\quad , \quad j \neq k\\[2mm]
\p_{x_j}^2\, T_{-m-2\ell} &=& (m+2\ell)(m+2\ell+2)\,  x_j^2\, T_{-m-2\ell-4}  -\, \frac{2m+4\ell+2}{ (m+2\ell+2)C(m,\ell)}\, a_m\, \p_{x_j}^2\, \Delta^\ell\, \delta(\ux)  
 -\, (m+2\ell)\, T_{-m-2\ell-2} -\, \frac{m+2\ell}{C(m,\ell+1)}\, a_m\, \Delta^{\ell+1}\, \delta(\ux)
\end{eqnarray*}
confirming (\ref{lapTspecial}) through
\begin{align*}
\Delta\, T_{-m-2\ell} &= (m+2\ell)(m+2\ell+2)\, T_{-m-2\ell-2}  -\, \frac{2m+4\ell+2}{ (m+2\ell+2)C(m,\ell)}\, a_m\,  \Delta^{\ell+1}\, \delta(\ux)  
 -\, m\, (m+2\ell)\, T_{-m-2\ell-2} -\, \frac{m(m+2\ell)}{C(m,\ell+1)}\, a_m\, \Delta^{\ell+1}\, \delta(\ux)\\
&= (m+2\ell)(2\ell+2)\, T_{-m-2\ell-2} + (-1)^\ell\, \frac{(m+4\ell+2)(m+2\ell+2)}{C(m,\ell+1)}\, a_m\, \dirac^{2\ell+2}\, \delta(\ux)
\end{align*}
It follows that
\begin{eqnarray*}
d_k\, d_j\, ^{s}U_{-m-2\ell} &=& (m+2\ell)(m+2\ell+2)\, x_k\, x_j\, ^{s}U_{-m-2\ell-4}  -\,\frac{2m+4\ell+2}{ (m+2\ell+2)C(m,\ell)}\,  a_m\, \uom\, \p_{x_k}\p_{x_j}\, \Delta^\ell\, \delta(\ux)\quad , \quad j \neq k\\[2mm]
d_j^2\, ^{s}U_{-m-2\ell} &=& (m+2\ell)(m+2\ell+2)\,  x_j^2\, ^{s}U_{-m-2\ell-4}  -\, \frac{2m+4\ell+2}{ (m+2\ell+2)C(m,\ell)}\, a_m\, \uom\, \p_{x_j}^2\, \Delta^\ell\, \delta(\ux)  
 -\, (m+2\ell)\, ^{s}U_{-m-2\ell-2} -\, \frac{m+2\ell}{C(m,\ell+1)}\, a_m\, \uom\, \Delta^{\ell+1}\, \delta(\ux)
\end{eqnarray*}
In particular, for $\la=-m$, we have
\begin{eqnarray*}
\p_{x_k}\p_{x_j}\, T_{-m} &=& m(m+2)\, x_k\, x_j\, T_{-m-4}  -\, \frac{2m+2}{m(m+2)}\, a_m\, \p_{x_k}\p_{x_j}\,  \delta(\ux)\quad , \quad j \neq k\\[2mm]
\p_{x_j}^2\, T_{-m} &=& m(m+2)\,  x_j^2\, T_{-m-4}  -\, \frac{2m+2}{m(m+2)}\, a_m\, \p_{x_j}^2\,  \delta(\ux) - 
m\, T_{-m-2} -\, \frac{1}{2(m+2)}\, a_m\, \Delta\, \delta(\ux)
\end{eqnarray*}
confirming (\ref{lapTm}) through
\begin{align*}
\Delta\, T_{-m} &= m(m+2)\,  T_{-m-2}  -\, \frac{2m+2}{m(m+2)}\, a_m\, \Delta\,  \delta(\ux) - 
m^2\, T_{-m-2} -\, \frac{m}{2(m+2)}\, a_m\, \Delta\, \delta(\ux)\\
&= 2\, m\, T_{-m} - \frac{m+2}{m}\, a_m\, \Delta\, \delta(\ux)
\end{align*}
It follows that
\begin{eqnarray*}
d_k\, d_j\, ^{s}U_{-m} &=& m(m+2)\, x_k\, x_j\, ^{s}U_{-m-4}  -\, \frac{2m+2}{m(m+2)}\, a_m\, \uom\, \p_{x_k}\p_{x_j}\,  \delta(\ux)\quad , \quad j \neq k\\[2mm]
d_j^2\, ^{s}U_{-m} &=& m(m+2)\,  x_j^2\, ^{s}U_{-m-4}  -\, \frac{2m+2}{m(m+2)}\, a_m\, \uom\, \p_{x_j}^2\,  \delta(\ux) - 
m\, ^{s}U_{-m-2} -\, \frac{1}{2(m+2)}\, a_m\, \uom\, \Delta\, \delta(\ux)
\end{eqnarray*}
A direct computation shows that
\begin{eqnarray*}
d_k\, d_j\, T_{-m-2\ell} &=& (m+2\ell+1)\, U_{-m-2\ell-3}\, (x_k\, e_j + x_j\, e_k) + (m+2\ell+1)(m+2\ell+3)\, x_j\, x_k\, T_{-m-2\ell-4}\\
&+& (-1)^{\ell+1}\, \frac{1}{(m+2\ell+2)C(m,\ell)}\, a_m\, (\p_{x_k} \dirac^{2\ell+1} \delta(\ux)\, e_j + \p_{x_j} \dirac^{2\ell+1}\delta(\ux)\, e_k) + 2(-1)^{\ell+1}\, \frac{a_m}{C(m,\ell)}\, \p_{x_j}\p_{x_k}\dirac^{2\ell} \delta(\ux)
\quad , \quad j \neq k\\[2mm]
d_j^2\, T_{-m-2\ell} &=& 2(m+2\ell+1)\, U_{-m-2\ell-3}\, x_j\, e_j - (m+2\ell+1)\, T_{-m-2\ell-2} + (m+2\ell+1)(m+2\ell+3)\, x_j^2\, T_{-m-2\ell-4}\\
&+& (-1)^\ell\, \frac{m+2\ell}{C(m,\ell+1)}\, a_m\, \dirac^{2\ell+2}\delta(\ux) + 2(-1)^{\ell+1}\, \frac{1}{(m+2\ell+2)C(m,\ell)}\, a_m\, \p_{x_j}\, \dirac^{2\ell+1}\delta(\ux)\, e_j + 2(-1)^{\ell+1}\, \frac{1}{C(m,\ell)}\, a_m\, \p_{x_j}^2\, \dirac^{2\ell}\delta(\ux)
\end{eqnarray*}
whence also
\begin{eqnarray*}
\p_{x_k}\p_{x_j}\, ^{s}U_{-m-2\ell} &=&  -(m+2\ell+1)\, ^{s}T_{-m-2\ell-3}\, (x_k\, e_j + x_j\, e_k) + (m+2\ell+1)(m+2\ell+3)\, x_j\, x_k\, ^{s}U_{-m-2\ell-4}\\
&+& (-1)^{\ell+1}\, \frac{1}{(m+2\ell+2)C(m,\ell)}\, a_m\, \uom\, (\p_{x_k} \dirac^{2\ell+1} \delta(\ux)\, e_j + \p_{x_j} \dirac^{2\ell+1}\delta(\ux)\, e_k) + 2(-1)^{\ell+1}\, \frac{a_m}{C(m,\ell)}\, \uom\, \p_{x_j}\p_{x_k}\dirac^{2\ell} \delta(\ux)
\quad , \quad j \neq k\\[2mm]
\p_{x_j}^2\, ^{s}U_{-m-2\ell} &=& -2(m+2\ell+1)\, ^{s}T_{-m-2\ell-3}\, x_j\, e_j - (m+2\ell+1)\, ^{s}U_{-m-2\ell-2} + (m+2\ell+1)(m+2\ell+3)\, x_j^2\, ^{s}U_{-m-2\ell-4}\\
&+& (-1)^\ell\, \frac{m+2\ell}{C(m,\ell+1)}\, a_m\, \uom\, \dirac^{2\ell+2}\delta(\ux) + 2(-1)^{\ell+1}\, \frac{1}{(m+2\ell+2)C(m,\ell)}\, a_m\, \uom\, \p_{x_j}\, \dirac^{2\ell+1}\delta(\ux)\, e_j + 2(-1)^{\ell+1}\, \frac{1}{C(m,\ell)}\, a_m\, \uom\, \p_{x_j}^2\, \dirac^{2\ell}\delta(\ux)
\end{eqnarray*}
confirming (\ref{lapsUspecial}) through
$$
\Delta\, ^{s}U_{-m-2\ell} = (m+2\ell+1)(2\ell+1)\, ^{s}U_{-m-2\ell-2} + (-1)^\ell\, (m+2\ell+2)(m+4\ell+2)\, \frac{a_m}{C(m,\ell+1)}\, \uom\, \dirac^{2\ell+2}\delta(\ux)
$$
In particular, for $\la=-m$, we have
\begin{eqnarray*}
d_k\, d_j\,\, T_{-m} &=&  (m+1)\, U_{-m-3}\, (x_k\, e_j + x_j\, e_k) + (m+1)(m+3)\, x_j\, x_k\, T_{-m-4}\\
&-& \frac{1}{m(m+2)}\, a_m\, (\p_{x_k} \dirac \delta(\ux)\, e_j + \p_{x_j} \dirac \delta(\ux)\, e_k) -  \frac{2}{m}\, a_m\, \p_{x_j}\p_{x_k} \delta(\ux)
\quad , \quad j \neq k\\[2mm]
d_j^2\, T_{-m} &=& 2(m+1)\, U_{-m-3}\, x_j\, e_j - (m+1)\, T_{-m-2} + (m+1)(m+3)\, x_j^2\, T_{-m-4}\\
&+&  \frac{1}{2(m+2)}\, a_m\, \dirac^{2}\delta(\ux) - 2\, \frac{1}{m(m+2)}\, a_m\, \p_{x_j}\, \dirac \delta(\ux)\, e_j - 2\, \frac{1}{m}\, a_m\, \p_{x_j}^2\, \delta(\ux)
\end{eqnarray*}
whence also
\begin{eqnarray*}
\p_{x_k}\p_{x_j}\, ^{s}U_{-m} &=&  -(m+1)\, ^{s}T_{-m-3}\, (x_k\, e_j + x_j\, e_k) + (m+1)(m+3)\, x_j\, x_k\, ^{s}U_{-m-4}\\
&-& \frac{1}{m(m+2)}\, a_m\, \uom\, (\p_{x_k} \dirac \delta(\ux)\, e_j + \p_{x_j} \dirac \delta(\ux)\, e_k) -  \frac{2}{m}\, a_m\, \uom\, \p_{x_j}\p_{x_k} \delta(\ux)
\quad , \quad j \neq k\\[2mm]
\p_{x_j}^2\, ^{s}U_{-m} &=& -2(m+1)\, ^{s}T_{-m-3}\, x_j\, e_j - (m+1)\, ^{s}U_{-m-2} + (m+1)(m+3)\, x_j^2\, ^ {s}U_{-m-4}\\
&+&  \frac{1}{2(m+2)}\, a_m\, \uom\, \dirac^{2}\delta(\ux) - 2\, \frac{1}{m(m+2)}\, a_m\, \uom\, \p_{x_j}\, \dirac \delta(\ux)\, e_j - 2\, \frac{1}{m}\, a_m\, \uom\, \p_{x_j}^2\, \delta(\ux)
\end{eqnarray*}
which confirms (\ref{lapsUm}). 


\newpage
\subsection{Second order partial derivatives of $U_\la$}


First we compute $\p_{x_k}\p_{x_j}\, U_\la$ assuming that $\la \neq -m+3, -m+1, -m-1, \ldots$, and we find
\begin{eqnarray*}
\p_{x_k}\p_{x_j}\, U_\la &=&  (\la-1)(\la-3)\, x_k\, x_j\, U_{\la-4} + (\la-1)\, T_{\la-3}\, (x_k\, e_j + x_j\, e_k) \quad , \quad j \neq k\\[2mm]
\p_{x_j}^2\, U_\la &=&  (\la-1)(\la-3)\,  x_j^2\, U_{\la-4} + 2(\la-1)\, T_{\la-3}\, x_j\, e_j + (\la-1)\, U_{\la-2}
\end{eqnarray*}
confirming (\ref{lapU}) through
\begin{align*}
\Delta\, U_\la &=  (\la-1)(\la-3)\,   U_{\la-2} + 2(\la-1)\, T_{\la-3}\, \ux + m(\la-1)\, U_{\la-2}\\
&= (\la-1)(\la+m-1)\, U_{\la-2}
\end{align*}
It follows at once that
\begin{eqnarray*}
d_k\, d_j\, ^{s}T_\la &=& (\la-1)(\la-3)\, x_k\, x_j\, ^{s}T_{\la-4} - (\la-1)\, ^{s}U_{\la-3}\, (x_k\, e_j + x_j\, e_k) \quad , \quad j \neq k\\[2mm]
d_j^2\, ^{s}T_\la &=& (\la-1)(\la-3)\,  x_j^2\, ^{s}T_{\la-4} - 2(\la-1)\, ^{s}U_{\la-3}\, x_j\, e_j + (\la-1)\, ^{s}T_{\la-2}
\end{eqnarray*}
A direct computation shows that
\begin{eqnarray*}
d_k\, d_j\, U_\la &=& \la\, (\la-2)\, x_k\, x_j\, U_{\la-4} \quad , \quad j \neq k\\[2mm]
 d_j^2\, U_\la &=& \la\, (\la-2)\,  x_j^2\, U_{\la-4} + \la\, U_{\la-2}
\end{eqnarray*}
from which it follows that
\begin{eqnarray*}
\p_{x_k}\, \p_{x_j}\, ^{s}T_\la &=& \la\, (\la-2)\, x_k\, x_j\, ^{s}T_{\la-4} \quad , \quad j \neq k\\[2mm]
 \p_{x_j}^2\, ^{s}T_\la &=& \la\, (\la-2)\, x_j^2\, ^{s}T_{\la-4} + \la\, ^{s}T_{\la-2}
\end{eqnarray*}
confirming (\ref{lapsT}) through
\begin{align*}
\Delta\, ^{s}T_\la &= \la\, (\la-2)\,  ^{s}T_{\la-2} + m\, \la\, ^{s}T_{\la-2}\\
&= \la(\la+m-2)\, ^{s}T_{\la-2}
\end{align*}

\noindent
For $\la = -m-2\ell+1, \ell=1,2,\ldots$ we obtain
\begin{eqnarray*}
\p_{x_k}\, \p_{x_j}\, U_{-m-2\ell+1} &=& -(m+2\ell)\, T_{-m-2\ell-2}\, (x_k\, e_j + x_j\, e_k) + (m+2\ell)(m+2\ell+2)\, x_j\, x_k\, U_{-m-2\ell-3}\\ 
&&   + (-1)^{\ell+1}\,  \frac{m+2\ell}{C(m,\ell+1)}\, a_m\, x_j\, x_k\, \dirac^{2\ell+3}\delta(\ux)  +  (-1)^{\ell+1}\, (2\ell)\, \frac{a_m}{C(m,\ell)}\, \p_{x_k}\, \p_{x_j}\, \dirac^{2\ell-1}\delta(\ux)
 \quad , \quad j \neq k\\[2mm]
\p_{x_j}^2\, U_{-m-2\ell+1} &=& -2(m+2\ell)\, x_j\, e_j\, T_{-m-2\ell-2} -  (m+2\ell)\, U_{-m-2\ell-1} +  (m+2\ell)(m+2\ell+2)\, x_j^2\, U_{-m-2\ell-3} \\
&& + (-1)^{\ell+1}\, \frac{m+2\ell}{C(m,\ell+1)}\, a_m\, x_j^2\, \dirac^{2\ell+3}\delta(\ux) +  (-1)^{\ell+1}\, \frac{2\ell}{C(m,\ell)}\, a_m\, \p_{x_j}^2\, \dirac^{2\ell-1}\delta(\ux)
\end{eqnarray*}
confirming (\ref{lapUspecial}) through
\begin{align*}
\Delta\, U_{-m-2\ell+1} &=  -2(m+2\ell)\, \ux\, T_{-m-2\ell-2} - m(m+2\ell)\, U_{-m-2\ell-1} + (m+2\ell)(m+2\ell+2)\, U_{-m-2\ell-1)} \\
&\phantom{m+2\ell}+ (-1)^{\ell}\, \frac{m+2\ell}{C(m,\ell+1)}\, a_m\, \ux^2\, \dirac^{2\ell+3}\delta(\ux) + (-1)^\ell\, \frac{2\ell}{C(m,\ell)}\, a_m\, \dirac^{2\ell+1}\delta(\ux)\\
&= (m+2\ell)(2\ell)\, U_{-m-2\ell-1} + (-1)^\ell\, \frac{m+4\ell}{C(m,\ell)}\, a_m\, \dirac^{2\ell+1}\delta(\ux)
\end{align*}

\noindent
Hence also
\begin{eqnarray*}
d_k\, d_j\, ^{s}T_{-m-2\ell+1} &=& (m+2\ell)\, ^{s}U_{-m-2\ell-2}\, (x_k\, e_j + x_j\, e_k) + (m+2\ell)(m+2\ell+2)\, x_j\, x_k\, ^{s}T_{-m-2\ell-3}\\ 
&&   + (-1)^{\ell}\,  \frac{m+2\ell}{C(m,\ell+1)}\, a_m\, x_j\, x_k\, \uom\, \dirac^{2\ell+3}\delta(\ux)  +  (-1)^{\ell}\, (2\ell)\, \frac{a_m}{C(m,\ell)}\, \uom\,  \p_{x_k}\, \p_{x_j}\, \dirac^{2\ell-1}\delta(\ux)
 \quad , \quad j \neq k\\[2mm]
d_j^2\, ^{s}T_{-m-2\ell+1} &=& 2(m+2\ell)\, x_j\, e_j\, ^{s}U_{-m-2\ell-2} -  (m+2\ell)\, ^{s}T_{-m-2\ell-1} +  (m+2\ell)(m+2\ell+2)\, x_j^2\, ^{s}T_{-m-2\ell-3} \\
&& + (-1)^{\ell}\, \frac{m+2\ell}{C(m,\ell+1)}\, a_m\, x_j^2\, \uom\, \dirac^{2\ell+3}\delta(\ux) +  (-1)^{\ell}\, \frac{2\ell}{C(m,\ell)}\, a_m\, \uom\, \p_{x_j}^2\, \dirac^{2\ell-1}\delta(\ux)
\end{eqnarray*}

\noindent
A direct computation shows that
\begin{eqnarray*}
d_k d_j\, U_{-m-2\ell+1} &=& (m+2\ell-1)(m+2\ell+1)\, x_j\, x_k\, U_{-m-2\ell-3} + (-1)^{\ell}\, \frac{2(m+2\ell)(2\ell+2)}{C(m,\ell+1)}\, a_m\, (e_k\, \p_{x_j}\, \dirac^{2\ell} \delta(\ux) + e_j\, \p_{x_k}\, \dirac^{2\ell} \delta(\ux))\\
&+& (-1)^{\ell+1}\, \frac{2(m+2\ell)(2\ell+2)(2\ell)}{C(m,\ell+1)}\, a_m\, \p_{x_j}\, \p_{x_k}\, \dirac^{2\ell-1} \delta(\ux)
 \quad , \quad j \neq k\\[2mm]
d_j^2\, U_{-m-2\ell+1} &=&  (m+2\ell-1)(m+2\ell+1)\, x_j^2\, U_{-m-2\ell-3} - (m+2\ell-1)\, U_{-m-2\ell-1} + (-1)^\ell\, \frac{4(m+2\ell)(2\ell+2)}{C(m,\ell+1)}\,  a_m\, e_j\, \p_{x_j}\, \dirac^{2\ell} \delta(\ux)\\
&+&   (-1)^{\ell+1}\, \frac{(2\ell+2)(2\ell)(2m+4\ell)}{C(m,\ell+1)}\, a_m\, \p_{x_j}^2\, \dirac^{2\ell-1} \delta(\ux) + (-1)^\ell\, \frac{(2\ell+2)(m+2\ell-2)}{C(m,\ell+1)}\, a_m\, \dirac^{2\ell+1} \delta(\ux)
\end{eqnarray*}
whence also
\begin{eqnarray*}
\p_{x_k}\, \p_{x_j}\, ^{s}T_{-m-2\ell+1} &=& (m+2\ell-1)(m+2\ell+1)\, x_j\, x_k\, ^{s}T_{-m-2\ell-3} + (-1)^{\ell+1}\, \frac{2(m+2\ell)(2\ell+2)}{C(m,\ell+1)}\, a_m\, \uom\, (e_k\, \p_{x_j}\, \dirac^{2\ell} \delta(\ux) + e_j\, \p_{x_k}\, \dirac^{2\ell} \delta(\ux))\\
&+& (-1)^{\ell}\, \frac{2(m+2\ell)(2\ell+2)(2\ell)}{C(m,\ell+1)}\, a_m\, \uom\, \p_{x_j}\, \p_{x_k}\, \dirac^{2\ell-1} \delta(\ux)
 \quad , \quad j \neq k\\[2mm]
\p_{x_j}^2\, ^{s}T_{-m-2\ell+1} &=&  (m+2\ell-1)(m+2\ell+1)\, x_j^2\, ^{s}T_{-m-2\ell-3} - (m+2\ell-1)\, ^{s}T_{-m-2\ell-1} + (-1)^{\ell+1}\, \frac{4(m+2\ell)(2\ell+2)}{C(m,\ell+1)}\,  a_m\, \uom\, e_j\, \p_{x_j}\, \dirac^{2\ell} \delta(\ux)\\
&+&   (-1)^{\ell}\, \frac{(2\ell+2)(2\ell)(2m+4\ell)}{C(m,\ell+1)}\, a_m\, \p_{x_j}^2\, \uom\, \dirac^{2\ell-1} \delta(\ux) + (-1)^{\ell+1}\, \frac{(2\ell+2)(m+2\ell-2)}{C(m,\ell+1)}\, a_m\, \uom\, \dirac^{2\ell+1} \delta(\ux)
\end{eqnarray*}
which confirms (\ref{lapsTspecial}).\\

\noindent
For $\la=-m+1$, it holds that
\begin{eqnarray*}
\p_{x_k}\, \p_{x_j}\, U_{-m+1} &=& (-m)\, T_{-m-2}\, (x_k\, e_j + x_j\, e_k) + m(m+2)\, x_j\, x_k\, U_{-m-3} - \frac{a_m}{2(m+2)}\, x_j\, x_k\, \dirac^3 \delta(\ux) \quad , \quad j \neq k\\[2mm]
\p_{x_j}^2\, U_{-m+1} &=& 2(-m)\, x_j\, e_j\, T_{-m-2} - m\, U_{-m-1} + m(m+2)\, x_j^2\, U_{-m-3} - \frac{a_m}{2(m+2)}\, x_j^2\, \dirac^3 \delta(\ux)
\end{eqnarray*}
confirming (\ref{lapUm1}) through
\begin{align*}
\Delta\, U_{-m+1} &=   2(-m)\, \ux\, T_{-m-2} - m^2\, U_{-m-1} + m(m+2)\, r^2\, U_{-m-3} - \frac{1}{2(m+2)}\, a_m\, r^2\, \dirac^3\delta(\ux)\\
&= a_m\, \dirac\delta(\ux)
\end{align*}

\noindent
Hence also
\begin{eqnarray*}
d_k\, d_j\, ^{s}T_{-m+1} &=& m\, ^{s}U_{-m-2}\, (x_k\, e_j + x_j\, e_k) + m(m+2)\, x_j\, x_k\, ^{s}T_{-m-3} + \frac{a_m}{2(m+2)}\, x_j\, x_k\, \uom\, \dirac^3 \delta(\ux) \quad , \quad j \neq k\\[2mm]
d_j^2\, ^{s}T_{-m+1} &=& 2\, m\, x_j\, e_j\, ^{s}U_{-m-2} - m\, ^{s}T_{-m-1} + m(m+2)\, x_j^2\, ^{s}T_{-m-3} + \frac{a_m}{2(m+2)}\, x_j^2\, \uom\,  \dirac^3 \delta(\ux)
\end{eqnarray*}

\noindent
We know that a direct computation of $d_k\, d_j\, U_{-m+1}$ and $d_j^2\, U_{-m+1}$ will involve arbitrary constants to be determined through the calculation of $\Delta\, ^{s}T_{-m+1}$. We can circumvent these computations by putting $\ell=0$ in the expressions for 
$d_k\, d_j\, U_{-m-2\ell+1}$ and $d_j^2\, U_{-m-2\ell+1}$ obtaining in this way:
\begin{eqnarray*}
d_k\, d_j\, U_{-m+1} &=&  (m-1)(m+1)\, x_k\, x_j\, U_{-m-3} + \frac{2}{m+2}\, a_m\, (e_k\, \p_{x_j} \delta(\ux) +  e_j\, \p_{x_k} \delta(\ux) )               \quad , \quad j \neq k\\[2mm]
d_j^2\, U_{-m+1} &=& (m-1)(m+1)\, xj^2\, U_{-m-3} - (m-1)\, U_{-m-1} + \frac{4}{m+2}\, a_m\,  e_j\, \p_{x_j} \delta(\ux)  + \frac{m-2}{m(m+2)}\, a_m\, \dirac \delta(\ux)
\end{eqnarray*}
whence also
\begin{eqnarray*}
\p_{x_k}\, \p_{x_j}\, ^{s}T_{-m+1} &=&   (m-1)(m+1)\, x_k\, x_j\, ^{s}T_{-m-3} - \frac{2}{m+2}\, a_m\, \uom\, (e_k\, \p_{x_j} \delta(\ux) +  e_j\, \p_{x_k} \delta(\ux) )        \quad , \quad j \neq k\\[2mm]
\p_{x_j}^2\, ^{s}T_{-m+1} &=& (m-1)(m+1)\, xj^2\, ^{s}T_{-m-3} - (m-1)\, ^{s}T_{-m-1} - \frac{4}{m+2}\, a_m\, \uom\,  e_j\, \p_{x_j} \delta(\ux)  - \frac{m-2}{m(m+2)}\, a_m\, \uom\, \dirac \delta(\ux) 
\end{eqnarray*}
which indeed confirms (\ref{lapsTm1}).



\newpage
\section{Powers of the vector variable $\ux$}
\label{powersux}


In this section we will explore the intimate connection between the $T_\la$ and $U_\la$ distributions on the one hand and powers of the Clifford vector variable $\ux$ on the other. 

\subsection{Integer powers of $\ux$ as regular (signum)distributions}

The integer powers $\ux^k, k \in \mZ$, of the vector variable $\ux$ remain locally integrable in $\mR^m$ as long as $k \geq -m+1$, and so may be interpreted as regular distributions, viz. for all test functions $\varphi(\ux) \in \mcD(\mR^m;\mC)$ one has
$$
\langle \ \ux^k \ ,  \ \varphi(\ux) \ \rangle =  \int_{\mR^m}\, \ux^k\, \varphi(\ux)\, d\ux \quad , \quad k \geq -m+1
$$
In particular for even integer powers: $k = 2\ell \geq -m+1$ we have
\begin{align*}
\langle \ \ux^{2\ell} \ ,  \ \varphi(\ux) \ \rangle \ &= \  \int_{\mR^m}\, \ux^{2\ell}\, \varphi(\ux)\, d\ux  \ = \ (-1)^\ell\, \int_{0}^{+\infty}\, r^{m+2\ell-1}\, dr\, \int_{S^{m-1}}\,  \varphi(r\uom)\, dS_{\uom}  \\[2mm]
&= \  (-1)^\ell\, a_m\, \int_0^{+\infty}\, r^{m+2\ell-1}\, \Sigma^{(0)} [\varphi]\, dr \ =  \ (-1)^\ell\, a_m\, \langle \ r_+^{m+2\ell-1} \ ,  \ \Sigma^{(0)} [\varphi] \ \rangle \\
\end{align*}
implying that
$$
\boxed{\ux^{2\ell} = (-1)^\ell\, T_{2\ell} \quad , \quad 2\ell \geq -m+1}
$$
which clearly are radial regular distributions.\\

\noindent
Similarly, for odd integer powers $k = 2\ell+1 \geq -m+1$ we find
\begin{align*}
\langle \ \ux^{2\ell+1} \ ,  \ \varphi(\ux) \ \rangle \ &=  \ \int_{\mR^m}\, \ux^{2\ell+1}\, \varphi(\ux)\, d\ux \  = \  (-1)^\ell\, \int_{0}^{+\infty}\, r^{m+2\ell}\, dr\, \int_{S^{m-1}}\, \uom\, \varphi(r\uom)\, dS_{\uom}  \\
&= \ (-1)^\ell\, a_m\, \int_0^{+\infty}\, r^{m+2\ell}\, \Sigma^{(1)} [\varphi]\, dr \ = \ (-1)^\ell\, a_m\, \langle \ r_+^{m+2\ell} \ ,  \ \Sigma^{(1)} [\varphi] \ \rangle \\
\end{align*}
implying that
$$
\boxed{\ux^{2\ell+1} = (-1)^\ell\, U_{2\ell+1} \quad , \quad 2\ell \geq -m}
$$
which clearly are signum--radial regular distributions.\\

\noindent
The associated signumdistributions of these integer powers of $\ux$ are:
$$
(\ux^{2\ell})^{\vee} = \uom\, \ux^{2\ell} = (-1)^\ell\, ^{s}U_{2\ell}
$$
and
$$
(\ux^{2\ell+1} )^{\vee} = \uom\, \ux^{2\ell+1}  = (-1)^{\ell+1}\, ^{s}T_{2\ell+1}
$$
but these signumdistributions are by no means powers of the vector variable $\ux$.\\

\noindent
As to derivation we have
\begin{align*}
\dirac\, \ux^{2\ell} &= -\, (2\ell)\, \ux^{2\ell-1} \qquad & \uD\, \ux^{2\ell} &= -\, (m+2\ell-1)\, \ux^{2\ell-1} \quad &, \quad &2\ell-1 \geq -m+1\\[2mm]
\dirac\, \ux^{2\ell+1} &= -\, (m+2\ell)\, \ux^{2\ell} \qquad  & \uD\, \ux^{2\ell+1} &= -\, (2\ell+1)\, \ux^{2\ell}  \quad &, \quad &2\ell \geq -m+1
\end{align*}
and the corresponding formulae for the associated signumdistributions:
\begin{align*}
\uD\, (\uom\, \ux^{2\ell} )&= -\, (2\ell)\, (\uom\, \ux^{2\ell-1} ) \qquad & \dirac\, ( \uom\, \ux^{2\ell})  &= -\, (m+2\ell-1)\, (\uom\, \ux^{2\ell-1} ) \quad &, \quad &2\ell-1 \geq -m+1\\[2mm]
\uD\, (\uom\, \ux^{2\ell+1}) &= -\, (m+2\ell)\, ( \uom\, \ux^{2\ell} ) \qquad  & \dirac\, ( \uom\, \ux^{2\ell+1})  &= -\, (2\ell+1)\, (\uom\, \ux^{2\ell} ) \quad &, \quad &2\ell \geq -m+1
\end{align*}\hfill\\

We know that locally integrable functions may serve as regular signumdistributions as well: for all test functions $\uom\, \varphi(\ux) \in \Omega(\mR^m; \mC)$ one has
$$
\langle \ ^{s}\ux^k \ ,  \ \uom\, \varphi(\ux) \ \rangle =  \int_{\mR^m}\, \ux^k\, \uom\, \varphi(\ux)\, d\ux \quad , \quad k \geq -m+1
$$
In particular for even integer powers: $k = 2\ell \geq -m+1$ we have
\begin{align*}
\langle \ ^{s}\ux^{2\ell} \ ,  \ \uom\, \varphi(\ux) \ \rangle \ &= \  \int_{\mR^m}\, \ux^{2\ell}\, \uom\, \varphi(\ux)\, d\ux  \ = \ (-1)^\ell\, \int_{0}^{+\infty}\, r^{m+2\ell-1}\, dr\, \int_{S^{m-1}}\,  \uom\, \varphi(r\uom)\, dS_{\uom}  \\
&=  \ (-1)^\ell\, a_m\, \int_0^{+\infty}\, r^{m+2\ell-1}\, \Sigma^{(1)} [\varphi]\, dr \ = \  (-1)^\ell\, a_m\, \langle \ r_+^{m+2\ell-1} \ ,  \ \Sigma^{(1)} [\varphi] \ \rangle \\
\end{align*}
implying that
$$
\boxed{^{s}\ux^{2\ell} = (-1)^\ell\, ^{s}T_{2\ell} \quad , \quad 2\ell \geq -m+1}
$$
which clearly are radial regular signumdistributions.\\

\noindent
Similarly, for odd integer powers $k = 2\ell+1 \geq -m+1$ we find
\begin{align*}
\langle \ ^{s}\ux^{2\ell+1} \ ,  \ \uom\, \varphi(\ux) \ \rangle \ &= \ \int_{\mR^m}\, \ux^{2\ell+1}\, \uom\, \varphi(\ux)\, d\ux  \ = \ (-1)^{\ell+1}\, \int_{0}^{+\infty}\, r^{m+2\ell}\, dr\, \int_{S^{m-1}}\,  \varphi(r\uom)\, dS_{\uom}  \\
&= \ (-1)^{\ell+1}\, a_m\, \int_0^{+\infty}\, r^{m+2\ell}\, \Sigma^{(0)} [\varphi]\, dr \ = \  (-1)^{\ell+1}\, a_m\, \langle \ r_+^{m+2\ell} \ ,  \ \Sigma^{(0)} [\varphi] \ \rangle \\
\end{align*}
implying that
$$
\boxed{^{s}\ux^{2\ell+1} = (-1)^\ell\, ^{s}U_{2\ell+1} \quad , \quad 2\ell \geq -m}
$$
which clearly signum--radial regular signumdistributions.\\

\noindent
The associated distributions of these integer power signumdistributions are:
$$
(^{s}\ux^{2\ell})^{\wedge} = -\, \uom\, ^{s} \ux^{2\ell} = (-1)^{\ell+1}\, U_{2\ell}
$$
and
$$
(^{s} \ux^{2\ell+1} )^{\wedge} = -\, \uom\, ^{s}\ux^{2\ell+1}  = (-1)^{\ell}\, T_{2\ell+1}
$$
but these distributions are by no means powers of the vector variable $\ux$.\\

\noindent
As to derivation we have
\begin{align*}
\dirac\, ^{s}\ux^{2\ell} &= -\, (2\ell)\, ^{s}\ux^{2\ell-1} \qquad & \uD\, ^{s}\ux^{2\ell} &= -\, (m+2\ell-1)\, ^{s}\ux^{2\ell-1} \quad &, \quad &2\ell-1 \geq -m+1\\[2mm]
\dirac\, ^{s}\ux^{2\ell+1} &= -\, (m+2\ell)\, ^{s}\ux^{2\ell} \qquad  & \uD\, ^{s}\ux^{2\ell+1} &= -\, (2\ell+1)\, ^{s}\ux^{2\ell}  \quad &, \quad &2\ell \geq -m+1
\end{align*}
and the corresponding formulae for the associated distributions:
\begin{align*}
\uD\, (\uom\, ^{s}\ux^{2\ell} )&= -\, (2\ell)\, (\uom\, ^{s}\ux^{2\ell-1} ) \qquad & \dirac\, ( \uom\, ^{s}\ux^{2\ell})  &= -\, (m+2\ell-1)\, (\uom\, ^{s}\ux^{2\ell-1} ) \quad &, \quad &2\ell-1 \geq -m+1\\[2mm]
\uD\, (\uom\, ^{s}\ux^{2\ell+1}) &= -\, (m+2\ell)\, ( \uom\, ^{s}\ux^{2\ell} ) \qquad  & \dirac\, ( \uom\, ^{s}\ux^{2\ell+1})  &= -\, (2\ell+1)\, (\uom\, ^{s}\ux^{2\ell} ) \quad &, \quad &2\ell \geq -m+1
\end{align*}\hfill\\

\subsection{Integer powers of $\ux$  as finite part (signum)distributions}

For $k < -m+1$ the functions $\ux^k$ are no longer locally integrable in $\mR^m$ and thus no longer  regular (signum)distributions. Whence the need for a definition.

\begin{definition}
One defines the distributions
\begin{align*}
{\rm Fp}\, \ux^{2\ell} &= (-1)^\ell\, T_{2\ell} \quad , \quad 2\ell < -m+1 \\[2mm]
{\rm Fp}\, \ux^{2\ell+1} &= (-1)^\ell\, U_{2\ell+1} \quad , \quad 2\ell < -m
\end{align*}
and the signumdistributions
\begin{align*}
{\rm Fp}\, ^{s}\ux^{2\ell} &= (-1)^\ell\, ^{s}T_{2\ell} \quad , \quad 2\ell < -m+1 \\[2mm]
{\rm Fp}\, ^{s}\ux^{2\ell+1} &= (-1)^\ell\, ^{s}U_{2\ell+1} \quad , \quad 2\ell < -m
\end{align*}
\end{definition}

Note that 
$$
\left( {\rm Fp}\, \ux^{2\ell} \right)^{\vee} =  (-1)^\ell\, ^{s}U_{2\ell}  \qquad \left( {\rm Fp}\, \ux^{2\ell+1} \right)^{\vee} =  (-1)^{\ell+1}\, ^{s}T_{2\ell+1}
$$
and
$$
\left( {\rm Fp}\, ^{s}\ux^{2\ell} \right)^{\wedge} =  (-1)^{\ell+1}\, U_{2\ell}  \qquad \left( {\rm Fp}\, ^{s}\ux^{2\ell+1} \right)^{\wedge} =  (-1)^{\ell}\, T_{2\ell+1}
$$
and observe that none of those four associated (signum)distributions are integer powers of $\ux$.

\begin{remark}
{\rm
We will keep the notations ${\rm Fp}\, \ux^k$ and ${\rm Fp}\, ^{s}\ux^k$ even when $\ux^k$ is locally integrable.
}
\end{remark}

As to division by the vector variable $\ux$ the following formul\ae\  may be easily calculated:
\begin{itemize}
\item $\invux\, {\rm Fp}\, \ux^{2\ell} = {\rm Fp}\, \ux^{2\ell-1}$
\item $\invux\, {\rm Fp}\, \ux^{2\ell+1} = {\rm Fp}\, \ux^{2\ell} \quad , \quad 2\ell \neq -m$
\item $\invux\, {\rm Fp}\, \ux^{-m+1} = [{\rm Fp}\, \ux^{-m}] = {\rm Fp}\, \ux^{-m} + \delta(\ux)\, c \quad , \quad m \ {\rm even}$
\end{itemize}\hfill\\

\subsection{Complex powers of $\ux$ as distributions}

Through the $T_\la$ and $U_\la$ distributions it is possible to define {\em complex} powers of the vector variable $\ux$ as distributions.

\begin{definition}
For $\la \in \mC$ one puts
$$
(i\, \ux)^\la = \onehalf\, (1 + \exp{(i\pi\la)})\, T_\la + \onehalf\, i\,  (1 - \exp{(i\pi\la)})\, U_\la
$$
\end{definition}

\begin{remark}
{\rm
It has to be checked that in the particular case where $\la$ is an integer the distributions introduced above are recovered. We indeed have:
\begin{itemize}
\item for $\la = 2\ell, \ell \in \mZ$ we get
$$
(i\, \ux)^{2\ell} =T_{2\ell } = (-1)^\ell\, {\rm Fp}\, \ux^{2\ell}
$$
\item for $\la = 2\ell+1, \ell \in \mZ$ we get
$$
(i\, \ux)^{2\ell+1} = i\, U_{2\ell+1 } = i\, (-1)^\ell\, {\rm Fp}\, \ux^{2\ell+1}
$$
\end{itemize}
}
\end{remark}

\noindent
Now we investigate the holomorphy of $(i\, \ux)^\la$ in the complex $\la$--plane. Recall that
\begin{itemize}
\item $T_\la$ is a meromorphic function of $\la$ showing simple poles at $\la = -m-2\ell, \ell=0,1,2,\ldots$ 
\item $U_\la$ is a meromorphic function of $\la$ showing simple poles at $\la = -m-2\ell-1, \ell=0,1,2,\ldots$ 
\end{itemize}

\noindent
So it becomes clear that possible singularities of the distributions $(i\, \ux)^\la$ depend on the parity of the dimension $m$.\\[2mm]
If the dimension $m$ is {\em odd}, then
$$
(i\, \ux)^{-m-2\ell} = i\, U_{-m-2\ell}
$$
and
$$
(i\, \ux)^{-m-2\ell-1} = T_{-m-2\ell-1}
$$
implying that the distribution $(i\, \ux)^\la$ then is an {\em entire} function of $\la \in \mC$ showing no singularities at all.\\[1mm]

\noindent
If the dimension $m$ is {\em even}, then
$$
(i\, \ux)^{-m-2\ell} = T_{-m-2\ell}
$$
and
$$
(i\, \ux)^{-m-2\ell-1} = i\, U_{-m-2\ell-1}
$$
implying that the distribution $(i\, \ux)^\la$ then is a meromorphic function of $\la \in \mC$ showing simple poles at  $\la = -m, -m-1, -m-2, \ldots$. But note that at those singular points the distribution  $(i\, \ux)^\la$ still is defined through the monomial pseudofunctions, however {\em not} turning this distribution into an entire function of $\la$.\\

Now we compute the Dirac derivative of the distribution $(i\, \ux)^\la$ . We  obtain in general
\begin{align}
\dirac\, (i\, \ux)^\la &= \onehalf\, (1 - \exp{(i \pi (\la-1))})\, \dirac\, T_\la + \onehalf\, i\, (1 + \exp{(i \pi (\la-1))})\, \dirac\, U_\la{\nonumber} \\[2mm]
\label{diracixla}&= (-i)\, \la\, (i\, \ux)^{\la-1} - \onehalf\, i\, (m-1)(1 - \exp{(i \pi\, \la)})\, T_{\la-1}
\end{align}

\noindent
This formula (\ref{diracixla}) is valid in the whole complex $\la$--plane when the dimension $m$ is odd, and valid in $\mC \backslash \{-m+1, -m, -m-1, \ldots\}$ when $m$ is even. 
For the singular values of the parameter $\la$ we obtain the following results.
\begin{itemize}
\item $\boxed{\la = -m+1}$\\[2mm]
We consider the distribution
$$
(i\, \ux)^{-m+1} = \onehalf\, (1 - (-1)^m)\, T_{-m+1} + \onehalf\, i\, (1 + (-1)^m)\, U_{-m+1}
$$
For odd dimension $m$ we do not expect singularities to appear and indeed we get
$$
(i\, \ux)^{-m+1} = T_{-m+1}
$$
with
$$
\dirac\, (i\, \ux)^{-m+1} =  (-m+1)\, U_{-m}
$$
whence
$$
\dirac\, {\rm Fp}\, \ux^{-m+1} = (m-1)\, {\rm Fp}\, \ux^{-m}
$$
For even dimension $m$ we get
$$
(i\, \ux)^{-m+1} = i\, U_{-m+1}
$$
with
$$
\dirac\, (i\, \ux)^{-m+1} =  i\, (- a_m)\, \delta(\ux)
$$
whence
$$
\dirac\, {\rm Fp}\, \ux^{-m+1} = (-1)^{m/2+1}\, a_m\, \delta(\ux)
$$

\item $\boxed{\la = -m}$\\[2mm]
We consider the distribution
$$
(i\, \ux)^{-m} = \onehalf\, (1 + (-1)^m)\, T_{-m} + \onehalf\, i\, (1 - (-1)^m)\, U_{-m}
$$
For odd dimension $m$ we do not expect singularities to appear and indeed we get
$$
(i\, \ux)^{-m} = i\, U_{-m}
$$
with
$$
\dirac\, (i\, \ux)^{-m} =  i\, T_{-m-1}
$$
whence
$$
\dirac\, {\rm Fp}\, \ux^{-m} =  {\rm Fp}\, \ux^{-m-1}
$$
For even dimension $m$ we get
$$
(i\, \ux)^{-m} = T_{-m}
$$
with
$$
\dirac\, (i\, \ux)^{-m} =  (-m)\, U_{-m-1} - \frac{1}{m}\,  a_m\, \dirac\, \delta(\ux)
$$
whence
$$
\dirac\, {\rm Fp}\, \ux^{-m} = m\, {\rm Fp}\, \ux^{-m-1} + (-1)^{m/2+1}\, \frac{1}{m}\, a_m\, \dirac\, \delta(\ux)
$$

\item $\boxed{\la=-m-2\ell}$\\[2mm]
We consider the distribution
$$
(i\, \ux)^{-m-2\ell} = \onehalf\, (1 + (-1)^m)\, T_{-m-2\ell} + \onehalf\, i\, (1 - (-1)^m)\, U_{-m-2\ell}
$$
For odd dimension $m$ we do not expect singularities to appear and indeed we get
$$
(i\, \ux)^{-m-2\ell} = i\, U_{-m-2\ell}
$$
with
$$
\dirac\, (i\, \ux)^{-m-2\ell} =  i\, (2\ell+1)\, T_{-m-2\ell-1}
$$
whence
$$
\dirac\, {\rm Fp}\, \ux^{-m-2\ell} =  (2\ell+1)\, {\rm Fp}\, \ux^{-m-2\ell-1}
$$
For even dimension $m$ we get
$$
(i\, \ux)^{-m-2\ell} = T_{-m-2\ell}
$$
with
$$
\dirac\, (i\, \ux)^{-m-2\ell} =  (-m-2\ell)\, U_{-m-2\ell-1} + (-1)^{\ell+1}\,  \frac{1}{C(m,\ell)}\,  a_m\, \dirac^{2\ell+1}\, \delta(\ux)
$$
whence
$$
\dirac\, {\rm Fp}\, \ux^{-m-2\ell} = (m+2\ell)\,  {\rm Fp}\, \ux^{-m-2\ell-1} + (-1)^{m/2+1}\, \frac{1}{C(m,\ell)}\, a_m\, \dirac^{2\ell+1}\, \delta(\ux)
$$

\item $\boxed{\la=-m-2\ell-1}$\\[2mm]
We consider the distribution
$$
(i\, \ux)^{-m-2\ell-1} = \onehalf\, (1 + (-1)^{m+1})\, T_{-m-2\ell-1} + \onehalf\, i\, (1 - (-1)^{m+1})\, U_{-m-2\ell-1}
$$
For odd dimension $m$ we do not expect singularities to appear and indeed we get
$$
(i\, \ux)^{-m-2\ell-1} = T_{-m-2\ell-1}
$$
with
$$
\dirac\, (i\, \ux)^{-m-2\ell-1} =  -\,  (m+2\ell+1)\, U_{-m-2\ell-2}
$$
whence
$$
\dirac\, {\rm Fp}\, \ux^{-m-2\ell-1} =  (m+2\ell+1)\, {\rm Fp}\, \ux^{-m-2\ell-2}
$$
For even dimension $m$ we get
$$
(i\, \ux)^{-m-2\ell-1} = i\, U_{-m-2\ell-1}
$$
with
$$
\dirac\, (i\, \ux)^{-m-2\ell-1} =  i\, (2\ell+2)\, T_{-m-2\ell-2} + i\,  (-1)^{\ell}\,  \frac{m+2\ell+2}{C(m,\ell+1)}\,  a_m\, \dirac^{2\ell+2}\, \delta(\ux)
$$
whence
$$
\dirac\, {\rm Fp}\, \ux^{-m-2\ell-1} = (2\ell+2)\,  {\rm Fp}\, \ux^{-m-2\ell-2} + (-1)^{m/2+1}\, \frac{m+2\ell+2}{C(m,\ell+1)}\, a_m\, \dirac^{2\ell+2}\, \delta(\ux)
$$
\end{itemize}\hfill\\

We may also compute the action of the operator $\uD$ on the complex powers of the vector variable $\ux$. Recall that the operator $\uD$ is the cartesian signum--partner operator to the Dirac operator $\dirac$. We  obtain in general
\begin{align}
\uD\, (i\, \ux)^\la &= \onehalf\, (1 + \exp{(i \pi \la)})\, \uD\, T_\la + \onehalf\, i\, (1 - \exp{(i \pi \la)})\, \uD\, U_\la{\nonumber} \\[2mm]
\label{uDixla}&= \onehalf\, (1 + \exp{(i \pi \la)})\, (\la+m-1)\, U_{\la-1} - \onehalf\, i\, (1 - \exp{(i \pi \la)})\, \la\, T_{\la-1}
\end{align}

\noindent
This formula (\ref{uDixla}) is valid in the whole complex $\la$--plane when the dimension $m$ is odd, and valid in $\mC \backslash \{-m+1, -m, -m-1, \ldots\}$ when $m$ is even. For the singular values of the parameter $\la$ we obtain the following expressions.

\begin{itemize}
\item \boxed{\la = -m+1}\\[2mm]
For odd dimension $m$ we do not expect singularities to appear and indeed we find
\begin{align*}
\uD\, (i \ux)^{-m+1} = \uD\, T_{-m+1} = 0
\end{align*}
whence
$$
\uD\, {\rm Fp}\, \ux^{-m+1} = 0
$$
For even dimension $m$ we find
$$
\uD\, (i \ux)^{-m+1} = i\, \uD\, U_{-m+1} = i\, (m-1)\, T_{-m} - i\, a_m\, \delta(\ux)
$$
whence
$$
\uD\, {\rm Fp}\, \ux^{-m+1} = (m-1)\, {\rm Fp}\, \ux^{-m} + (-1)^{-m/2+1}\, a_m\, \delta(\ux)
$$

\item \boxed{\la = -m}\\[2mm]
For odd dimension $m$ we do not expect singularities to appear and indeed we find
$$
\uD\, (i \ux)^{-m} = i\, \uD\, U_{-m} = i\, m\, T_{-m-1}
$$
whence
$$
\uD\, {\rm Fp}\, \ux^{-m} = m\, {\rm Fp}\, \ux^{-m-1}
$$
For even dimension $m$ we find
$$
\uD\, (i \ux)^{-m} = \uD\, T_{-m} = -\, U_{-m-1} - \frac{1}{m}\, a_m\, \dirac \delta(\ux)
$$
whence
$$
\uD\, {\rm Fp}\, \ux^{-m} = {\rm Fp}\, \ux^{-m-1} + (-1)^{-m/2+1}\, \frac{1}{m}\,  a_m\, \dirac \delta(\ux)
$$

\item \boxed{\la = -m-2\ell}\\[2mm]
For odd dimension $m$ we do not expect singularities to appear and indeed we find
$$
\uD\, (i \ux)^{-m-2\ell} = i\, \uD\, U_{-m-2\ell} = i\, (m+2\ell)\, T_{-m-2\ell-1}
$$
whence
$$
\uD\, {\rm Fp}\, \ux^{-m-2\ell} = (m+2\ell)\, {\rm Fp}\, \ux^{-m-2\ell-1}
$$
For even dimension $m$ we find
$$
\uD\, (i \ux)^{-m-2\ell} = \uD\, T_{-m-2\ell} = -\, (2\ell+1)\, U_{-m-2\ell-1} + \frac{(-1)^{\ell+1}}{C(m,\ell)}\, a_m\, \dirac^{2\ell+1}\delta(\ux)
$$
whence
$$
\uD\, {\rm Fp}\, \ux^{-m-2\ell} = (2\ell+1)\, {\rm Fp}\, \ux^{-m-2\ell-1} + (-1)^{-m/2 + 1}\, \frac{1}{C(m,\ell)}\, a_m\, \dirac^{2\ell+1}\delta(\ux)
$$

\item \boxed{\la = -m-2\ell-1}\\[2mm]
For odd dimension $m$ we do not expect singularities to appear and indeed we find
$$
\uD\, (i \ux)^{-m-2\ell-1} = \uD\, T_{-m-2\ell-1} = -\, (2\ell+2)\, U_{-m-2\ell-2}
$$
whence
$$
\uD\, {\rm Fp}\, \ux^{-m-2\ell-1} = (2\ell+2)\, {\rm Fp}\, \ux^{-m-2\ell-2}
$$
For even dimension $m$ we find
$$
\uD\, (i \ux)^{-m-2\ell-1} = i\, \uD\, U_{-m-2\ell-1} = i\, (m+2\ell+1)\, T_{-m-2\ell-2} + i\, (-1)^{\ell}\, \frac{m+2\ell+2}{C(m,\ell+1)}\, a_m\, \dirac^{2\ell+2}\delta(\ux)
$$
whence
$$
\uD\, {\rm Fp}\, \ux^{-m-2\ell-1} = (m+2\ell+1)\, {\rm Fp}\, \ux^{-m-2\ell-2} + (-1)^{-m/2-1}\, \frac{m+2\ell+2}{C(m,\ell+1)}\, a_m\, \dirac^{2\ell+2}\delta(\ux)
$$

\end{itemize}\hfill\\

\subsection{Complex powers of $\ux$ as signumdistributions}

Similarly as in the preceding subsection we can define signumdistributions involving complex powers of $\ux$.

\begin{definition}
For $\la \in \mC$ one puts
$$
^{s}(i\, \ux)^\la = \onehalf\, (1 + \exp{(i\pi\la)})\, ^{s}T_\la + \onehalf\, i\,  (1 - \exp{(i\pi\la)})\, ^{s}U_\la
$$
\end{definition}

\begin{remark}
{\rm
It has to be checked that for integer values of the complex parameter $\la$ we recover signumdistributions already introduced above. We indeed have
\begin{itemize}
\item for $\la = 2\ell, \ell \in \mZ$
$$
^{s}(i\, \ux)^{2\ell} =\,  ^{s}T_{2\ell } = (-1)^\ell\, {\rm Fp}\, ^{s}\ux^{2\ell}
$$
\item for $\la = 2\ell+1, \ell \in \mZ$
$$
^{s}(i\, \ux)^{2\ell+1} = i\, ^{s}U_{2\ell+1 } = i\, (-1)^\ell\, {\rm Fp}\, ^{s}\ux^{2\ell+1}
$$
\end{itemize}
}
\end{remark}

\noindent
Now we investigate the holomorphy of $^s(i\, \ux)^\la$ in the complex $\la$--plane. Recall that
\begin{itemize}
\item $^sT_\la$ is a meromorphic function of $\la$ showing simple poles at $\la = -m-2\ell-1, \ell=0,1,2,\ldots$ 
\item $^sU_\la$ is a meromorphic function of $\la$ showing simple poles at $\la = -m-2\ell, \ell=0,1,2,\ldots$ 
\end{itemize}

\noindent
So it becomes clear that possible singularities of the signumdistributions $^s(i\, \ux)^\la$ depend on the parity of the dimension $m$.\\[2mm]
If the dimension $m$ is {\em even}, then
$$
^s(i\, \ux)^{-m-2\ell} =\, ^sT_{-m-2\ell}
$$
and
$$
^s(i\, \ux)^{-m-2\ell-1} =\, i\, ^sU_{-m-2\ell-1}
$$
implying that the signumdistribution $^s(i\, \ux)^\la$ then is an {\em entire} function of $\la \in \mC$ showing no singularities at all.\\[1mm]

\noindent
If the dimension $m$ is {\em odd}, then
$$
^s(i\, \ux)^{-m-2\ell} = i\, ^sU_{-m-2\ell}
$$
and
$$
^s(i\, \ux)^{-m-2\ell-1} =\, ^sT_{-m-2\ell-1}
$$
implying that the signumdistribution $^s(i\, \ux)^\la$ then is a meromorphic function of $\la \in \mC$ showing simple poles at  $\la = -m, -m-1, -m-2, \ldots$. But note that at those singular points the signumdistribution  $^s(i\, \ux)^\la$ still is defined through the monomial pseudofunctions, however {\em not} turning this signumdistribution into an entire function of $\la$.\\

Now we compute the Dirac derivative of the signumdistribution $^s(i\, \ux)^\la$. In general we obtain
\begin{equation}\label{diracsixla}
\dirac\, ^s(i \ux)^\la = \onehalf\, (1 + \exp{(i \pi \la)})\, \la\, ^sU_{\la-1} - \onehalf\, i\, (1 - \exp{(i \pi \la)})\, (\la+m-1)\, ^sT_{\la-1}
\end{equation}
This formula (\ref{diracsixla}) is valid in the whole complex $\la$--plane when the dimension $m$ is even, while for odd dimension $m$ it is valid for $\la \neq -m+1, -m, -m-1, -m-2, \ldots$. For those singular values of the parameter $\la$ we obtain the following results.

\begin{itemize}
\item \boxed{\la = -m+1}\\[2mm]
For even dimension $m$ we do not expect singularities to appear and indeed we find
$$
\dirac\, ^s(i \ux)^{-m+1} = i\, \dirac\, ^sU_{-m+1} = 0
$$
whence
$$
\dirac\, {\rm Fp}\, ^s \ux^{-m+1} = 0
$$
For odd dimension $m$ we find
$$
\dirac\, ^s(i \ux)^{-m+1} = \dirac\, ^sT_{-m+1} = -\, (m-1)\, ^sU_{-m} + a_m\, \uom\, \delta(\ux)
$$
whence
$$
\dirac\, {\rm Fp}\, ^s\ux^{-m+1} = (m-1)\, {\rm Fp}\, ^s\ux^{-m} + (-1)^{\frac{-m+1}{2}}\, a_m\, \uom\, \delta(\ux)
$$

\item \boxed{\la = -m}\\[2mm]
For even dimension $m$ we do not expect singularities to appear and indeed we find
$$
\dirac\, ^s(i \ux)^{-m} = \dirac\, ^sT_{-m} = -\, m\, ^sU_{-m-1}
$$
whence
$$
\dirac\, {\rm Fp}\, ^s \ux^{-m} = m\, {\rm Fp}\, ^s\ux^{-m-1}
$$
For odd dimension $m$ we find
$$
\dirac\, ^s(i \ux)^{-m} = i\, \dirac\, ^sU_{-m} = i\, ^sT_{-m-1} - i\, \frac{1}{m}\, a_m\, \uom\, \dirac\delta(\ux)
$$
whence
$$
\dirac\, {\rm Fp}\, ^s\ux^{-m} = {\rm Fp}\, ^s\ux^{-m-1} + (-1)^{\frac{-m+1}{2}}\, \frac{1}{m}\, a_m\, \uom\, \dirac\delta(\ux)
$$

\item \boxed{\la = -m-2\ell}\\[2mm]
For even dimension $m$ we do not expect singularities to appear and indeed we find
$$
\dirac\, ^s(i \ux)^{-m-2\ell} = \dirac\, ^sT_{-m-2\ell} = (- m-2\ell)\, ^sU_{-m-2\ell-1}
$$
whence
$$
\dirac\, {\rm Fp}\, ^s \ux^{-m-2\ell} = (m+2\ell)\, {\rm Fp}\, ^s\ux^{-m-2\ell-1}
$$
For odd dimension $m$ we find
$$
\dirac\, ^s(i \ux)^{-m-2\ell} = i\, \dirac\, ^sU_{-m-2\ell} = i\, \left( (2\ell+1)\, ^sT_{-m-2\ell-1}  + (-1)^{\ell+1}\, \frac{1}{C(m,\ell)}\, a_m\, \uom\, \dirac^{2\ell+1}\delta(\ux)  \right)
$$
whence
$$
\dirac\, {\rm Fp}\, ^s\ux^{-m-2\ell} = (2\ell+1)\, {\rm Fp}\, ^s\ux^{-m-2\ell-1} + (-1)^{\frac{-m+1}{2}}\, \frac{1}{C(m,\ell)}\, a_m\, \uom\, \dirac^{2\ell+1}\delta(\ux)
$$

\item \boxed{\la = -m-2\ell-1}\\[2mm]
For even dimension $m$ we do not expect singularities to appear and indeed we find
$$
\dirac\, ^s(i \ux)^{-m-2\ell-1} = i\, \dirac\, ^sU_{-m-2\ell-1} = i\, (2\ell+2)\, ^sT_{-m-2\ell-2}
$$
whence
$$
\dirac\, {\rm Fp}\, ^s \ux^{-m-2\ell-1} = (2\ell+2)\, {\rm Fp}\, \ux^{-m-2\ell-2}
$$
For odd dimension $m$ we find
$$
\dirac\, ^s(i \ux)^{-m-2\ell-1} = \dirac\, ^sT_{-m2\ell-1} = -\, (m+2\ell+1)\, ^sU_{-m-2\ell-2} + 
(-1)^{\ell+1}\, \frac{m+2\ell+2}{C(m,\ell+1)}\, a_m\, \uom\, \dirac^{2\ell+2}\delta(\ux)
$$
whence
$$
\dirac\, {\rm Fp}\, ^s\ux^{-m-2\ell-1} = (m+2\ell+1)\, {\rm Fp}\, ^s\ux^{-m-2\ell-2} + (-1)^{\frac{-m+1}{2}}\, \frac{m+2\ell+2}{C(m,\ell+1)}\, a_m\, \uom\, \dirac^{2\ell+2}\delta(\ux)
$$

\end{itemize}

Finally we can compute the action of the operator  $\uD$ on the signumdistribution $^s(i\, \ux)^\la$. In general we obtain
\begin{equation}\label{uDsixla}
\uD\, ^s(i \ux)^\la = \onehalf\, (1 + \exp{(i \pi \la)})\, (\la+m-1)\, ^sU_{\la-1} - \onehalf\, i\, (1 - \exp{(i \pi \la)})\, \la\, ^sT_{\la-1}
\end{equation}
This formula (\ref{uDsixla}) is valid in the whole complex $\la$--plane when the dimension $m$ is even, while for odd dimension $m$ it is valid for $\la \neq -m+1, -m, -m-1, -m-2, \ldots$. For those singular values of the parameter $\la$ we obtain the following results.

\begin{itemize}
\item \boxed{\la = -m+1}\\[2mm]
For even dimension $m$ we do not expect singularities to appear and indeed we find
$$
\uD\, ^s(i \ux)^{-m+1} = i\, \uD\, ^sU_{-m+1} = i\, (m-1)\, ^sT_{-m}
$$
whence
$$
\uD\, {\rm Fp}\, ^s \ux^{-m+1} = (m-1)\, {\rm Fp}\, ^s\ux^{-m}
$$
For odd dimension $m$ we find
$$
\uD\, ^s(i \ux)^{-m+1} = \uD\, ^sT_{-m+1} = a_m\, \uom\, \delta(\ux)
$$
whence
$$
\uD\, {\rm Fp}\, ^s\ux^{-m+1} = (-1)^{\frac{-m+1}{2}}\, a_m\, \uom\, \delta(\ux)
$$

\item \boxed{\la = -m}\\[2mm]
For even dimension $m$ we do not expect singularities to appear and indeed we find
$$
\uD\, ^s(i \ux)^{-m} =  \uD\, ^sT_{-m} = -\, ^sU_{-m-1}
$$
whence
$$
\uD\, {\rm Fp}\, ^s \ux^{-m} =  {\rm Fp}\, ^s\ux^{-m-1}
$$
For odd dimension $m$ we find
$$
\uD\, ^s(i \ux)^{-m} = i\, \uD\, ^sU_{-m} = i\, \left(  m\, ^sT_{-m-1} - \frac{1}{m}\,   a_m\, \uom\, \dirac \delta(\ux)  \right)
$$
whence
$$
\uD\, {\rm Fp}\, ^s\ux^{-m} = m\, {\rm Fp}\, ^s\ux^{-m-1} + (-1)^{\frac{-m+1}{2}}\, \frac{1}{m}\,  a_m\, \uom\, \dirac \delta(\ux)
$$

\item \boxed{\la = -m-2\ell}\\[2mm]
For even dimension $m$ we do not expect singularities to appear and indeed we find
$$
\uD\, ^s(i \ux)^{-m-2\ell} =  \uD\, ^sT_{-m-2\ell} = -\, (2\ell+1)\, ^sU_{-m-2\ell-1}
$$
whence
$$
\uD\, {\rm Fp}\, ^s \ux^{-m-2\ell} =  (2\ell+1)\, {\rm Fp}\, ^s\ux^{-m-2\ell-1}
$$
For odd dimension $m$ we find
$$
\uD\, ^s(i \ux)^{-m-2\ell} = i\, \uD\, ^sU_{-m-2\ell} = i\, \left(  (m+2\ell)\, ^sT_{-m-2\ell-1} + (-1)^{\ell+1}\, \frac{1}{C(m,\ell)}\,   a_m\, \uom\, \dirac^{2\ell+1} \delta(\ux)  \right)
$$
whence
$$
\uD\, {\rm Fp}\, ^s\ux^{-m-2\ell} = (m+2\ell)\, {\rm Fp}\, ^s\ux^{-m-2\ell-1} + (-1)^{\frac{-m+1}{2}}\, \frac{1}{C(m,\ell)}\,  a_m\, \uom\, \dirac^{2\ell+1} \delta(\ux)
$$

\item \boxed{\la = -m-2\ell-1}\\[2mm]
For even dimension $m$ we do not expect singularities to appear and indeed we find
$$
\uD\, ^s(i \ux)^{-m-2\ell-1} = i\, \uD\, ^sU_{-m-2\ell-1} = i\, (m+2\ell+1)\, ^sT_{-m-2\ell-2}
$$
whence
$$
\uD\, {\rm Fp}\, ^s \ux^{-m-2\ell-1} = (m+2\ell+1)\, {\rm Fp}\, ^s\ux^{-m-2\ell-2}
$$
For odd dimension $m$ we find
$$
\uD\, ^s(i \ux)^{-m-2\ell-1} = \uD\, ^sT_{-m-2\ell-1} = -\, (2\ell+2)\, ^sU_{-m-2\ell-2} + 
(-1)^{\ell+1}\, \frac{m+2\ell+2}{C(m,\ell+1)}\, a_m\, \uom\, \dirac^{2\ell+2} \delta(\ux)
$$
whence
$$
\uD\, {\rm Fp}\, ^s\ux^{-m-2\ell-1} = (2\ell+2)\, {\rm Fp}\, ^s\ux^{-m-2\ell-2}  +
(-1)^{\frac{-m+1}{2}}\, \frac{m+2\ell+2}{C(m,\ell+1)}\, a_m\, \uom\, \dirac^{2\ell+2} \delta(\ux)
$$

\end{itemize}


\newpage
\section{Normalization of the distributions $T_\la$ and $U_\la$}
\label{regularization}


As we know  the distributions $T_\la$, considered as functions of the complex parameter $\la$, show simple poles at $\la = -m, -m-2, -m-4, \ldots$, while the distributions $U_\la$ show simple poles at $\la = -m-1, -m-3, -m-5, \ldots$. In these singular points the distributions $T_\la$ and $U_\la$ were defined through the monomial pseudofunctions.\\

\noindent
There is however a second option to cope with these singularities: removing the singularities by dividing the distributions $T_\la$ and $U_\la$ by an appropriate Gamma-function. This gives rise to the so--called {\em normalized distributions} $T^{*}_\lambda$ and $U^{*}_\lambda$. Their definition runs as follows.

\noindent
The normalized distributions $T^{*}_\lambda$ are defined by
\begin{eqnarray*}
	\left \{
\begin{array}{ll}
\displaystyle{T_\lambda^* = \pi^{\frac{\lambda+m}{2}} \frac{T_\lambda}{\Gamma \left ( \frac{\lambda+m}{2} \right )}}, & \lambda \ne -m-2\ell\\[5mm]
\displaystyle{T_{-m-2\ell}^* = \frac{\pi^{\frac{m}{2}-\ell}}{2^{2l} \Gamma \left ( \frac{m}{2} + \ell \right )} \dirac^{2\ell} \delta (\ux)}, & \ell \in \mN_0
\end{array}
\right . 
\end{eqnarray*}

\noindent
The normalized distributions $U^{*}_\lambda$ are defined by
\begin{eqnarray*}
\left \{
\begin{array}{ll}
\displaystyle{U_\lambda^* = \pi^{\frac{\lambda+m+1}{2}} \, \frac{U_\lambda}{\Gamma \left ( \frac{\lambda + m + 1}{2} \right )}}, & \lambda \ne -m-2\ell-1\\[5mm]
\displaystyle{U_{-m-2\ell-1}^* = - \frac{\pi^{\frac{m}{2}-\ell}}{2^{2\ell+1} \, \Gamma \left ( \frac{m}{2} + \ell + 1 \right )} \; \dirac^{2\ell+1} \delta(\ux)}, & \ell \in \mN_0
\end{array}
\right . 
\end{eqnarray*}

\noindent
The normalized distributions $T_\lambda^*$ and $U_\lambda^*$ turn out to be holomorphic mappings from $\lambda \in \mC$ to the space $\mathcal{S}'(\mR^m)$ of tempered distributions. They are intertwined by the actions of the multiplication operator $\ux$ and of the Dirac operator, according to the following formul\ae: for all $\lambda \in \mC$ one has
\begin{itemize}
\item[(i)] $\ux \; T_\lambda^* = \frac{\lambda+m}{2\pi} \; U_{\lambda+1}^*$; \quad
$\ux \; U_\lambda^* = U_\lambda^* \; \ux = - T_{\lambda+1}^*$
\item[(ii)] $\dirac \; T_\lambda^* = \lambda \; U_{\lambda-1}^*$; \quad 
$\dirac \; U_\lambda^* = U_\lambda^* \; \dirac = - 2\pi \; T_{\lambda-1}^*$
\item [(iii)] $\Delta T_\lambda^* = 2 \pi \lambda T_{\lambda-2}^*$ ; \quad
$\Delta U_\lambda^* = 2 \pi (\lambda-1) U_{\lambda-2}^*$
\item[(iv)] ${\mathcal F} \left [ T_\lambda^* \right ] = T_{-\lambda-m}^*$ ; \quad
${\mathcal F} \left [ U_\lambda^* \right ] = -i \; U_{-\lambda-m}^*$.
\end{itemize}
where, in property (iv),  the following definition of the Fourier transformation has been adopted:
$$
{\mathcal F}[f(\ux)](\uy) = \int_{\mR^m} f(\ux) \exp \left ( -2 \pi i \, \langle \ux,\uy \rangle \right ) \, d\ux .
$$
For an in-depth study of the normalized distributions $T_\la^*$ and $U_\la^*$ we refer to e.g. \cite{differential}. Additionally we can investigate the behaviour of the $T_\la^*$ and $U_\la^*$ distributions under the action of the operators $\invux$ and $\uD$.

\begin{proposition}
\label{invuxstar}
One has\\[2mm]
(i) $\invux\, T^*_\la = -\, U^*_{\la-1}$,  for all $\la \in \mC$\\[2mm]
(ii) $\invux\, U^*_\la = \frac{2\pi}{\la+m-1}\, T^*_{\la-1}$, for all $\la \in \mC$ except for $\la = -m+1$
\end{proposition}

\pf
(i) For $\la \neq -m-2\ell$ we have
$$
\invux\, T^*_\la = \pi^{\frac{\lambda+m}{2}} \frac{-\, U_{\la-1}}{\Gamma \left ( \frac{\lambda+m}{2} \right )} = -\, U^*_{\la-1}
$$
(i') For $\la = -m-2\ell$ we have
$$
\invux\, T^*_{-m-2\ell} = \frac{\pi^{\frac{m}{2}-\ell}}{2^{2l} \Gamma \left ( \frac{m}{2} + \ell \right )} \invux\, \dirac^{2\ell} \delta (\ux) = \frac{\pi^{\frac{m}{2}-\ell}}{2^{2l} \Gamma \left ( \frac{m}{2} + \ell \right )}\, \frac{1}{m+2\ell}\, \dirac^{2\ell+1} \delta (\ux) = -\, U^*_{-m-2\ell-1}
$$
(ii) For $\la \neq -m-2\ell+1$ we have
$$
\invux\, U^*_\la =  \pi^{\frac{\lambda+m+1}{2}} \, \frac{T_{\la-1}}{\Gamma \left ( \frac{\lambda + m + 1}{2} \right )} = \frac{2\pi}{\la+m-1}\, T^*_{\la-1}
$$
(ii') For $\la = -m-2\ell-1$ we have
$$
\invux\, U_{-m-2\ell-1} = - \frac{\pi^{\frac{m}{2}-\ell}}{2^{2\ell+1} \, \Gamma \left ( \frac{m}{2} + \ell + 1 \right )} \, \invux\,  \dirac^{2\ell+1} \delta(\ux) = - \frac{\pi^{\frac{m}{2}-\ell}}{2^{2\ell+1} \, \Gamma \left ( \frac{m}{2} + \ell + 1 \right )} \, \frac{1}{2\ell+2}\,  \dirac^{2\ell+2} \delta(\ux) = \frac{2\pi}{-2\ell-2}\, T^*_{-m-2\ell-2}
$$
\eop

\begin{remark}
{\rm
Up to a constant the distribution $U^*_{-m+1}$ equals the distribution $U_{-m+1}$, which, recall, is a locally integrable function. We know that $\invux\, U_{-m+1}$ equals an equivalence class of distributions, viz.
$$
\invux\, U_{-m+1} = T_{-m} + \delta(\ux)\, c = T_{-m} + T^*_{-m}\, c'
$$
which clearly, due to the appearance of $T_{-m}$, do not belong tot the $T^*_\la$ family of distributions. This is the reason why we have left out the case $\la=-m+1$ in the above proposition.
}
\end{remark}

\begin{proposition}
\label{uDstar}
One has\\[2mm]
(i) $\uD\, T^*_\la = (\la+m-1)\, U^*_{\la-1}$,  for all $\la \in \mC$\\[2mm]
(ii) $\uD\, U^*_\la = -\, \frac{2 \pi \la}{\la+m-1}\, T^*_{\la-1}$, for all $\la \in \mC$ except for $\la = -m+1$
\end{proposition}

\pf
(i) As $T^*_\la$ is a radial distribution we have for all $\la \in \mC$ that
$$
\uD\, T^*_\la = \left( \uom\, \p_r - (m-1)\, \invux  \right)\, T^*_\la = \dirac\, T^*_\la - (m-1)\, \invux\, T^*_\la =  (\la+m-1)\, U^*_{\la-1}
$$
(ii) For $\la \neq -m-2\ell-1$ we have
$$
\uD\, U^*_{\la} = \left(  \uom\, \p_r - \invr\, \p_{\uom} - (m-1)\, \invux \right)\, U^*_\la =
\frac{\pi^{\frac{\la+m+1}{2}}}{\Gamma{(\frac{\la+m+1}{2})}}\, (-\la+m-1)\, T_{\la-1} - (m-1)\, \frac{2\pi}{\la+m-1}\, T^*_{\la-1} =  -\, \frac{2 \pi \la}{\la+m-1}\, T^*_{\la-1}
$$
(ii') For $\la = -m-2\ell-1$ we have
$$
\uD\, U^*_{-m-2\ell-1} = - \frac{\pi^{\frac{m}{2}-\ell}}{2^{2\ell+1} \, \Gamma \left ( \frac{m}{2} + \ell + 1 \right )} \; \dirac^{2\ell+1}  \left( \frac{m+2\ell+1}{2\ell+2} + \frac{m-1}{2\ell+2} - \frac{m-1}{2\ell+2} \right)\, \dirac^{2\ell+2} \delta(\ux) = -\, 2\pi\, \frac{m+2\ell+1}{2\ell+2}\, T^*_{-m-2\ell-2}
$$
\eop

In a similar way the signumdistributions  $^sT_\la$ and $^sU_\la$ may now be normalized. Recall that $^sT_\la$ shows simple poles at $\la = -m-1, -m-3, \ldots$, while $^sU_\la$ shows simple poles at $\la = -m, -m-2, \ldots$. We define:
\begin{eqnarray*}
\left \{
\begin{array}{ll}
\displaystyle{^sT_\lambda^* = \pi^{\frac{\lambda+m+1}{2}} \, \frac{^sT_\lambda}{\Gamma \left ( \frac{\lambda + m + 1}{2} \right )}}, & \lambda \ne -m-2\ell-1\\[5mm]
\displaystyle{^sT_{-m-2\ell-1}^* =  \frac{\pi^{\frac{m}{2}-\ell}}{2^{2\ell+1} \, \Gamma \left ( \frac{m}{2} + \ell + 1 \right )} \, \uom\, \dirac^{2\ell+1} \delta(\ux)}, & \ell \in \mN_0
\end{array}
\right . 
\end{eqnarray*}
and
\begin{eqnarray*}
	\left \{
\begin{array}{ll}
\displaystyle{^sU_\lambda^* = \pi^{\frac{\lambda+m}{2}} \frac{^sU_\lambda}{\Gamma \left ( \frac{\lambda+m}{2} \right )}}, & \lambda \ne -m-2\ell\\[5mm]
\displaystyle{^sU_{-m-2\ell}^* = \frac{\pi^{\frac{m}{2}-\ell}}{2^{2l} \Gamma \left ( \frac{m}{2} + \ell \right )}\, \uom\,  \dirac^{2\ell} \delta (\ux)}, & \ell \in \mN_0
\end{array}
\right . 
\end{eqnarray*}

\begin{lemma}
For all $\la \in \mC$ one has that
$$
(T_\la^*)^\vee = \uom\, T_\la^* =\,  ^sU_\la^*
$$
and
$$
(U_\la^*)^\vee = \uom\, U_\la^* = -\,  ^sT_\la^*
$$
\end{lemma}

\begin{proposition}
For all $\la \in \mC$ one has that
$$
\ux\, ^sT_\la^* =\,  ^sU^*_{\la+1}
$$
and
$$
\ux\, ^sU_\la^* = -\, \frac{\la+m}{2\pi}\, ^sT^*_{\la+1}
$$
\end{proposition}

\pf
Obviously these result follow by simple transition to the associated signumdistributions of the corresponding formul\ae \ for the distributions $T^*_\la$ and $U^*_\la$, the multiplication operator $\ux$ being a signum--self--adjoint operator.  However let us give a direct proof.\\[2mm]
(i) For $ \lambda \ne -m-2\ell-1$ one has $\ux\, ^sT_\lambda^* = \pi^{\frac{\lambda+m+1}{2}} \, \frac{\ux\, ^sT_\lambda}{\Gamma \left ( \frac{\lambda + m + 1}{2} \right )} = \pi^{\frac{\lambda+m+1}{2}} \, \frac{ ^sU_\la+1}{\Gamma \left ( \frac{\lambda + m + 1}{2} \right )} =\, ^sU^*_{\la+1}$.\\[2mm]
(ii) $\ux\, ^sT_{-m-2\ell-1}^* =  \frac{\pi^{\frac{m}{2}-\ell}}{2^{2\ell+1} \, \Gamma \left ( \frac{m}{2} + \ell + 1 \right )} \, \ux\, \uom\, \dirac^{2\ell+1} \delta(\ux) = \frac{\pi^{\frac{m}{2}-\ell}}{2^{2\ell+1} \, \Gamma \left ( \frac{m}{2} + \ell + 1 \right )} \, \uom\, (m+2\ell)\, \dirac^{2\ell} \delta(\ux) = \frac{\pi^{\frac{m}{2}-\ell}}{2^{2\ell} \, \Gamma \left ( \frac{m}{2} + \ell  \right )} \,  \uom\, \dirac^{2\ell} \delta(\ux) =\,  ^sU^*_{-m-2\ell} $\\\\[2mm]
(iii) For $ \lambda \ne -m-2\ell$ one has $\ux\, ^sU_\lambda^* = \pi^{\frac{\lambda+m}{2}} \frac{\ux\, ^sU_\lambda}{\Gamma \left ( \frac{\lambda+m}{2} \right )} = -\, \pi^{\frac{\lambda+m}{2}} \frac{ ^sT_\lambda+1}{\Gamma \left ( \frac{\lambda+m}{2} \right )} = -\, \frac{\la+m}{2\pi}\, ^sT^*_{\la+1}$.\\\\[2mm]
(iv) $\ux\, ^sU_{-m-2\ell}^* = \frac{\pi^{\frac{m}{2}-\ell}}{2^{2l} \Gamma \left ( \frac{m}{2} + \ell \right )}\, \ux\, \uom\,  \dirac^{2\ell} \delta (\ux) = \frac{\pi^{\frac{m}{2}-\ell}}{2^{2l} \Gamma \left ( \frac{m}{2} + \ell \right )}\, (2\ell)\, \uom\,  \dirac^{2\ell-1} \delta (\ux) = \frac{\ell}{\pi}\, ^sT^*_{-m-2\ell+1}$.
\eop
\\

\begin{proposition}
For all $\la \in \mC$ one has that
$$
\uD\, ^sT_\la^* = 2 \pi\, ^sU_{\la-1}^*
$$
and
$$
\uD\, ^sU_\la^* = -\, \la\, ^sT^*_{\la-1}
$$
\end{proposition}

\pf
These results simply follow from the corresponding formul\ae \ for the distributions $T^*_\la$ and $U^*_\la$, the operator $\uD$ being the signum--adjoint operator of the Dirac operator $\dirac$.
\eop\\

\noindent
In the similar way the results of the propositions \ref{invuxstar} and \ref{uDstar} lead to the following formul\ae.

\begin{proposition}\hfill\\[2mm]
(i)  For $\la \neq -m+1$ one has\\[3mm]
$\invux\, ^sT^*_\la = -\, \frac{2\pi}{\la+m-1}\, ^sU^*_{\la-1}$ \quad and \quad $\dirac\, ^sT^*_\la = \frac{2\pi\la}{\la+m-1}\, ^sU^*_{\la-1} $.\\[3mm]
(ii) For all $\la \in \mC$ one has\\[3mm]
$\invux\, ^sU^*_\la = ^sT^*_{\la-1}$ \quad and \quad $\dirac\, ^sU^*_\la = -\, (\la+m-1)\, ^sT^*_{\la-1}$.
\end{proposition}



\newpage
\section{Physics in threedimensional Euclidean space}
\label{threedim}


Consider standard three dimensional Euclidean space $(m=3)$. We will constantly switch between Clifford algebra and vector algebra notation, identifying the Clifford one--vector $\underline{v}$ with  the algebraic vector $\vec{v}$. Note that for the scalar product of two algebraic vectors it holds that $\vec{v} \circ \vec{w} = -\, \underline{v} \cdot \underline{w}$. In particular we identify the Dirac operator $\dirac$ with the gradient operator $\vec{\nabla}$.   Alternative notations used are  $\uomr$ for $\uom$,  $r^\la$ for $T_{\la}$  and $\uomr\, r^\la$ for $U_\la$. In this section we will discuss some well--known formulae from physics.

\subsection{}

A vector field $\vec{F}$ in $\mR^3$ is called {\em rotation free} (or {\em irrotational}) in an open region $\Omega$ if it is continuously differentiable in $\Omega$ and satisfies
$$
\vec{\nabla} \times \vec{F} = 0
$$
Invoking Stokes's Theorem it holds that the line integral of a rotation free vector field over a closed smooth path $\mcC$ in $\Omega$ vanishes:
$$
\oint_{\mcC}\, \vec{F}\, \circ \vec{dP} = 0
$$
If $\vec{F}$ stands for a force field then the above line integral represents the {\em work} done by this field to move an object around the path $\mcC$. This explains why e.g. a planet's orbiting around the sun in the Newtonian gravitational field, which is proportional to
$$
\vec{F}_{grav} = \frac{\uomr}{r^2} = \frac{\vec{x}}{|\vec{x}|^3}
$$
and is rotation free in $\mR^3 \backslash \{  O\}$, does not cost any energy.\\[1mm]
In Clifford algebra notation, this gravitational field $\vec{F}_{grav}$ takes the form, with dimension $m=3$,
$$
\underline{F}_{grav} = \frac{\uom}{r^2} = U_{-m+1}
$$
The irrotationality of  $\vec{F}_{grav} $ in $\mR^3 \backslash \{  O\}$ corresponds with
$$
\dirac \wedge U_{-m+1} = 0
$$
which follows from
$$
\dirac\, U_{-m+1} = -\, a_m\, \delta(\ux)
$$
the result being a scalar expression. As an aside note that in the complement of the origin it holds that
$$
\dirac \cdot U_{-m+1} =0
$$
which means that the gravitational field $\vec{F}_{grav}$ is also {\em divergence free} in $\mR^3 \backslash \{  O\}$.\\

But, upon inspection of the formulae obtained in Section \ref{actions}, it becomes clear that, for all $\la \in \mC$, $\dirac\, U_\la$ is a scalar expression, in other words:
$$
\dirac \wedge U_\la = 0 \quad , \quad \forall \la \in \mC
$$
which implies that all vector fields of the form
$$
\vec{F} = \frac{\uomr}{r^\la} 
$$
are rotation free in  $\mR^3 \backslash \{  O\}$. This is particularly interesting for e.g.  {\em Modified Newtonian Dynamics}, or MOND for short, which aims, as an alternative to the so--called {\em dark matter} model, at explaining the high velocities of the stars at the outskirts of galaxies. In this MOND model it is proposed to have, from a certain treshold on in the acceleration,  the gravitational field to be proportional to 
$$
\frac{\uomr}{r}
$$
As has been showed here, it is clear  that whatever power of the radial distance is used in a gravitational field, it always will  obey the zero energy orbiting principle.

\subsection{}

In three dimensional space the formulae for the second order cartesian derivatives of the fundamental solution of the Laplace operator $T_{-m+2}$, obtained in subsection 9.3, viz.
$$
\p_{x_k}\p_{x_j}\, T_{-m+2}  = m(m-2)\, x_k\, x_j\, T_{-m-2} \quad , \quad j \neq k
$$
and
$$
\p_{x_j}^2\, T_{-m+2} = m(m-2)\, x_j^2\ T_{-m-2} - (m-2)\, T_{-m} - \frac{m-2}{m}\, a_m\ \delta(\ux)
$$
turn into
\begin{align*}
\p_{x_k}\p_{x_j}\, T_{-1}  &= 3\, x_k\, x_j\, T_{-5} \quad , \quad j \neq k\\
\p_{x_j}^2\, T_{-1} &= 3\, x_j^2\, T_{-5} -  T_{-3} - \frac{4}{3}\,  \pi\ \delta(\ux)
\end{align*}
or, in the new notation,
\begin{align}
\label{frahmjk}
\p_{x_k}\p_{x_j}\, \invr  &= 3\, x_k\, x_j\, \frac{1}{r^5} \quad , \quad j \neq k\\
\label{frahmjj}
\p_{x_j}^2\, \invr &= 3\, x_j^2\, \frac{1}{r^5}  - \frac{4}{3}\,  \pi\, \delta(\ux) -  \frac{1}{r^3}
\end{align}

These expressions for the second order partial derivatives of the Coulomb potential of a unit point charge are of course well--known in physics, see e.g. \cite{frahm, hnidzo, adkins}.
In \cite{franklin} it is argued that, with respect to test functions that are not smooth at the origin, the above formulae should be replaced by
\begin{align}
\label{franklinjk}
\p_{x_k}\p_{x_j}\, \invr  &= 3\, x_k\, x_j\, \frac{1}{r^5} - 4\, \pi\, \frac{x_j\, x_k}{r^2}\, \delta(\ux) \quad , \quad j \neq k\\
\label{franklinjj}
\p_{x_j}^2\, \invr &= 3\, x_j^2\, \frac{1}{r^5}  - 4\,  \pi\,  \frac{x_j^2}{r^2}\, \delta(\ux) -  \frac{1}{r^3}
\end{align}
Mathematically speaking it is not clear what in \cite{franklin} is meant by ''test functions that are not smooth at the origin``. Nevertheless, taking into account that the delta distribution and its derivatives are only defined on differentiable test functions, it is readily observed that 
formulae (\ref{franklinjk}) and (\ref{franklinjj}) simply {\em coincide} with formulae (\ref{frahmjk}) and (\ref{frahmjj}), seen the following results.

\begin{lemma}
For the delta distribution $\delta(\ux)$ in $\mR^m$ one has\\[2mm]
(i) $x_j^2\, \Delta\, \delta(\ux) = 2\, \delta(\ux)$\\[2mm]
(ii) $x_j\, x_k\, \Delta\, \delta(\ux) = 0 \quad , \quad j \neq k$\\[2mm]
(iii) $r^2\, \Delta\, \delta(\ux) = 2\, m\, \delta(\ux)$\\[2mm]
(iv) $\invrsq\, \delta(\ux) = \frac{1}{2m}\, \Delta\, \delta(\ux)$\\[2mm]
(v) $\frac{x_j^2}{r^2}\, \delta(\ux) = \frac{1}{m}\, \delta(\ux)$\\[2mm]
(vi) $\frac{x_j\, x_k}{r^2}\, \delta(\ux) = 0 \quad , \quad j \neq k$
\end{lemma}

\pf
{\em (i)} For any test function $\varphi(\ux)$ and $j=1,\ldots,m$ it holds that
$$
\langle \  x_j^2\, \Delta \delta \ , \ \varphi \  \rangle  = \langle \ \delta \ , \ \Delta\,  (x_j^2\, \varphi) \  \rangle = \langle \ \delta  \ , \ 2\, \varphi + 4\, x_j\, \p_{x_j}\varphi + x_j^2\, \Delta\varphi\  \rangle = 2\, \langle \ \delta \ , \ \varphi \  \rangle
$$
{\em (ii)} One has
$$
\langle \  x_j\, x_k\, \Delta \delta \ , \ \varphi \  \rangle  = \langle \ \delta \ , \ \Delta\,  (x_j\, x_k\, \, \varphi) \  \rangle = \langle \ \delta  \ , 2\, x_k\, \p_{x_j}\varphi + 2\, x_j\,  \p_{x_k}\varphi  + x_j\, x_k\, \, \Delta\varphi\  \rangle = 0
$$
{\em(iii)} It follows from {\em (i)} that
$$
r^2\, \Delta\, \delta =  \Sigma_{j=1}^m\, x_j^2\, \Delta\, \delta = 2\, m\, \delta
$$
{\em (iv)} Taking into account the homogeinity of the delta distribution, it follows from {\em (iii)} that
$$
 \Delta\, \delta = 2\, m\, \invrsq\, \delta
$$
{\em (v)} It follows from {\em (iv)} and {\em (i)} that
$$
\frac{x_j^2}{r^2}\, \delta(\ux) = x_j^2\, \frac{1}{2m}\, \Delta\, \delta(\ux) = \frac{1}{m}\, \delta(\ux)
$$
{\em (vi)} Follows from {\em (iv)} and {\em (ii)}.
\eop\\

\noindent
As a corollary it then follows from {\em (vi)} that the extra term in (\ref{franklinjk}) compared to (\ref{frahmjk}) is indeed zero. Putting $m=3$ in {\em (v)} shows that 
$$
4\,  \pi\,  \frac{x_j^2}{r^2}\, \delta(\ux)  = \frac{4}{3}\, \pi\, \delta(\ux)
$$
thus making (\ref{franklinjj}) coincide with (\ref{frahmjj}).\\

\subsection{}

An interesting formula in \cite{franklin} is 
\begin{equation}
\label{franklin2}
\vec{\nabla}\, \frac{1}{r^3} = 4\, \pi\, \invr\, \uomr\, \delta(\ux) - 3\, \uomr\, \frac{1}{r^4}
\end{equation}
How does it fit into the present theory? Start with formula (\ref{diracTm}), viz.
$$
\dirac\, T_{-m} = -\, m\, U_{-m-1} - \frac{1}{m}\, a_m\, \dirac\, \delta(\ux)
$$
Put $m=3$ and make use of
$$
\frac{1}{m}\, \dirac\, \delta(\ux) = \invux\, \delta(\ux) = - \frac{\uom}{r}\, \delta(\ux)
$$
to obtain
$$
\dirac\, T_{-3} = -\, 3\, U_{-4} + 4\, \pi\, \frac{\uom}{r}\, \delta(\ux)
$$
which in the vector notation reads
$$
\vec{\nabla}\, \frac{1}{r^3} = - 3\, \uomr\, \frac{1}{r^4} +4\, \pi\, \invr\, \uomr\, \delta(\ux) 
$$
which is exactly (\ref{franklin2}).  In \cite{franklin} this result is also written as
$$
\vec{\nabla}\, \frac{1}{r^3} = g(r)\, \uomr
$$
with
$$
g(r) = \uomr \circ \vec{\nabla}\, \frac{1}{r^3} = - 3\, \frac{1}{r^4} + 4\, \pi\, \invr\, \delta(\ux) 
$$
Clearly $g(r)$ is the radial signumdistribution
$$
- 3\, ^{s}T_{-4} + 4\, \pi\, \invr\, \delta(\ux) = -\, \uom\, (\dirac\, T_{-3})
$$

\subsection{}

The electric potential induced by the electric dipole moment $\vec{p}$ is given by 
$$
\phi(\vec{x}) = \vec{p}  \circ \frac{\vec{x}}{r^3}
$$
and the corresponding electric field is 
$$
\vec{E} = -\, \vec{\nabla}\, \phi
$$
In \cite{frahm} it is shown that this electric field may be expressed directly in terms of the electric dipole moment as
$$
\vec{E} = -\, \frac{4\pi}{3}\, \delta(\ux)\, \vec{p} +  3 \left( \vec{p} \circ \frac{\vec{x}}{r^5}\right) \vec{x} - \frac{1}{r^3}\, \vec{p}
$$
It is our aim now to recover this expression through the formulae obtained in the underlying paper. Let $\up$ be the constant Clifford vector corresponding with the vector $\vec{p}$. Then, in general dimension $m$, it holds that
$$
\phi = \vec{p}  \circ \frac{\vec{x}}{r^m} = -\, \up \cdot U_{-m+1} 
$$
From
$$
\p_{x_i}\, (\up\, U_{-m+1}) = \up\, \p_{x_i} U_{-m+1} = \up\, (T_{-m}\, e_i - m\, x_i\, U_{-m-1} + \frac{1}{m}\, a_m\, \delta(\ux)\, e_i)
$$
it follows that 
$$
\p_{x_i}\, \phi = -\,  \up \cdot \p_{x_i}\, U_{-m+1}  = p_i\, T_{-m} + m\, x_i\,  (\up \cdot U_{-m-1})  + \frac{1}{m}\, a_m\, p_i\, \delta(\ux)
$$
whence
$$
\uE = -\, \up\, T_{-m} - m\, (\up \cdot \ux)\, U_{-m-1} - \frac{1}{m}\, a_m\, \delta(\ux)\, \up
$$
since
$$
(\up \cdot U_{-m-1})\, \ux = (\up \cdot \ux)\, U_{-m-1}
$$
For $m=3$ and in vector notation we indeed obtain
$$
\vec{E} = -\, \frac{4\pi}{3}\, \delta(\ux)\, \vec{p} + 3 \left( \vec{p} \circ \frac{\vec{x}}{r^5}\right) \vec{x} - \frac{1}{r^3}\, \vec{p}
$$
In $\mR^m \backslash \{ O\}$ the electric field $\vec{E}$ is divergence and rotation free, since it is the gradient field of a harmonic scalar field $\phi$. So if we compute the divergence and the rotation of $\vec{E}$ as a distribution in $\mR^m$, we expect the results to be distributions supported at the origin.  We can compute both at the same time by acting with the Dirac operator on $\uE$. We obtain
\begin{align*}
\dirac\, \uE &=    -\, \left(-\, m\, U_{-m-1} - \frac{1}{m}\, a_m\, \dirac \delta(\ux) \right)\, \up   + m\, \up\, U_{-m-1} -\, m\, (\up \cdot \ux)\, \left(  2\, T_{-m-2} + \frac{1}{2m}\, a_m\, \dirac^2 \delta(\ux)\right) - \frac{1}{m}\, a_m\, \dirac \delta(\ux)\, \up\\
&= 2\, m\, \up \cdot U_{-m-1} - 2\, m\, (\up \cdot \ux)\, T_{-m-2} -\, \onehalf\, a_m\, (\up \cdot \ux)\, \dirac^2 \delta(\ux)\\
&=   \onehalf\, a_m\, \Sigma\, p_j\, x_j\, \Delta \delta(\ux)\\
&= -\, a_m\, (\up \cdot \dirac) \delta(\ux)
\end{align*}
from which it follows that
$$
\dirac \cdot \uE = -\, a_m\, (\up \cdot \dirac) \delta(\ux)  \quad {\rm and } \quad \dirac \wedge \uE = 0
$$
or in dimension 3 and in vector analysis language:
$$
\vec{\nabla} \circ \vec{E} = -\, 4\pi (\vec{p} \circ \vec{\nabla}) \delta(\ux) \quad {\rm and } \quad \vec{\nabla} \times \vec{E} = 0
$$
showing that the electric field of an electric dipole is rotation free in the whole space even if it is not differentiable at the origin. It's divergence may also be obtained via
$$
\vec{\nabla} \circ \vec{E}  = -\, \Delta\, \phi = \up \cdot \Delta\, U_{-m+1} = a_m\, (\up \cdot \dirac) \delta(\ux) = - a_m\, (\vec{p} \circ \vec{\nabla}) \delta(\ux)
$$
\\

In a similar way if $\vec{M}$ stands for the magnetic dipole moment, the magnetic vector potential is given by
$$
\vec{A} = \vec{M}\, \times\, \frac{\vec{x}}{r^3}
$$
and the corresponding magnetic field is
$$
\vec{B} = \vec{\nabla}\, \times\, \vec{A}
$$
We will now recover the formula, appearing  in \cite{frahm},  expressing the magnetic field in terms of the magnetic dipole moment:
$$
\vec{B} = \frac{8\pi}{3}\, \delta(\ux)\, \vec{M} + 3\, (\vec{M} \circ \vec{x})\, \frac{\vec{x}}{r^5} - \frac{1}{r^3}\, \vec{M}
$$
First we observe that
$$
\vec{B} = \vec{\nabla}\, \times\, \vec{A} = \vec{\nabla}\, \times\, (\vec{M}\, \times\, \frac{\vec{x}}{r^3}) = \left( \vec{\nabla} \circ \frac{\vec{x}}{r^3} \right)\, \vec{M} - \left( \vec{M} \circ \vec{\nabla} \right)\, \frac{\vec{x}}{r^3}
$$
In general dimension $m$ we have
$$
\vec{\nabla} \circ \frac{\vec{x}}{r^m} = -\, \left[ \dirac\, U_{-m+1} \right]_0 = a_m\, \delta(\ux)
$$
and
\begin{align*}
\left( \vec{M} \circ \vec{\nabla} \right)\, \frac{\vec{x}}{r^m} &= \left( \Sigma_{j=1}^m\, M_j\, \p_{x_j} \right)\, U_{-m+1} = \Sigma_{j=1}^m\, M_j\, \left(  T_{-m}\, e_j - m\, x_j\, U_{-m-1} + \frac{1}{m}\, a_m\, \delta(\ux)\, e_j \right)\\
&= T_{-m}\, \vec{M} - m\, \left(\Sigma_{j=1}^m\, M_j\, x_j\right)\, U_{-m-1} + \frac{1}{m}\, a_m\, \delta(\ux)\, \vec{M}
\end{align*}
whence 
$$
\underline{B} = \frac{m-1}{m}\, a_m\, \delta(\ux)\, \underline{M} - T_{-m}\, \underline{M} +m\, (\vec{M} \circ \vec{x})\, U_{-m-1}
$$
and so, putting $m=3$,
\begin{align*}
\vec{B} &= 4\pi\, \delta(\ux)\, \vec{M} - \frac{1}{r^3}\, \vec{M} + 3\, (\vec{M} \circ \vec{x})\, \frac{\vec{x}}{r^5} - \frac{4\pi}{3}\, \delta(\ux)\, \vec{M}\\
&= \frac{8\pi}{3}\, \delta(\ux)\, \vec{M} + 3\, (\vec{M} \circ \vec{x})\, \frac{\vec{x}}{r^5} - \frac{1}{r^3}\, \vec{M}
\end{align*}
This magnetic field $\vec{B}$ satisfies the Maxwell equation
$$
\vec{\nabla} \circ \vec{B} = 0
$$
which would be trivial since $\vec{B}$ is a rotation field, were it not that $\vec{B}$ is not differentiable at the origin. To show that the magnetic field is indeed divergence-free we act with the Dirac operator on $\underline{B}$ and obtain
\begin{align*}
\dirac\, \uB &=  \frac{m-1}{m}\, a_m\, \dirac \delta(\ux) \uM - \left( -\, m\, U_{-m-1} - \frac{1}{m}\, a_m\, \dirac \delta(\ux) \right)\, \uM  + m\, \Sigma\, e_j\, M_j\, U_{-m-1} - m\, (\uM \cdot \ux) \left( 2\, T_{-m-2} + \frac{m+2}{2m(m+2)}\, a_m\, \dirac^2 \delta(\ux) \right) \\ 
&= a_m\, \dirac \delta(\ux)\, \uM + m\, \left(  U_{-m-1}\, \uM + \uM\, U_{-m-1} \right) - 2\, m\, (\uM \cdot \ux)\, T_{-m-2} + a_m\,  \Sigma\, M_j\, \p_{x_j}\, \delta(\ux) \\       
&= a_m\, \dirac \delta(\ux)\,  \uM  - a_m\, \dirac \delta(\ux) \cdot \uM \\
&=  a_m\, \dirac \delta(\ux) \wedge \uM
\end{align*}
from which it follows that
$$
\dirac \cdot \uB = 0 \quad {\rm and} \quad \dirac \wedge \uB = a_m\, \dirac \delta(\ux) \wedge \uM
$$
or, in the vector analysis language of three dimensional space,
$$
\vec{\nabla} \circ \vec{B} = 0 \quad {\rm and} \quad \vec{\nabla} \times \vec{B} = 4\pi\, \vec{\nabla} \delta(\ux) \times \vec{M} = -\, 4\pi\, (\vec{M} \times \vec{\nabla}) \delta(\ux)
$$

\subsection{}

In \cite{adkins}  a number of identities for the delta distribution are given as corollaries to the computation of integrals of spherical Bessel functions. It is shown e.g. that
$$
\p_{x_i}\, \p_{x_j}\, \invrsq - \onethird\, \delta_{ij}\, \Delta\, \invrsq = \frac{8}{r^4}\, \left( \omega_i\, \omega_j\, - \onethird\ \delta_{ij} \right)
$$
or, taking into account that in general $\Delta\, T_{-m+1} = (m-1)\, T_{-m-1}$, and thus in three dimensional space: $\Delta\, T_{-2} = 2\, T_{-4}$,
$$
\p_{x_i}\, \p_{x_j}\, \invrsq = \frac{8}{r^4}\,  \omega_i\, \omega_j\, - 2\, \delta_{ij}\,  \frac{1}{r^4} 
$$
We verify these formulae as follows. If $i \neq j$ then for the distribution $T_{-m+1}$ it holds that
$$
\p_{x_i}\, \p_{x_j}\, T_{-m+1}= (-m+1)(-m-1)\, x_i\, x_j\, T_{-m-3}
$$
which for $m=3$ turns into
$$
\p_{x_i}\, \p_{x_j}\, T_{-2}= 8\, x_i\, x_j\, T_{-6} = 8\, \omega_i\, \omega_j\, T_{-4}
$$
If $i=j$ then we have
$$
\p_{x_i}^2\, T_{-m+1} =  (-m+1)(-m-1)\, x_i^2\, T_{-m-3} + (-m+1)\, T_{-m-1}
$$
which for $m=3$ turns into
$$
\p_{x_i}^2\, T_{-2} =  8\, x_i^2\, T_{-6} - 2\, T_{-4} = 8\, \omega_i^2\, T_{-4} - 2\, T_{-4}
$$

\noindent
It is also shown in \cite{adkins} that
$$
\p_{x_i}\, \p_{x_j}\, \delta(\ux) - \onethird\, \delta_{ij}\, \Delta\, \delta(\ux) = 15\, \invrsq\, \omega_i\, \omega_j\, \delta(\ux) - 5\, \invrsq\, \delta_{ij}\, \delta(\ux)
$$
or, taking into account that in general $\invrsq\, \delta(\ux) = \frac{1}{2m}\, \Delta\, \delta(\ux)$, and thus in three dimensional space: $\invrsq\, \delta(\ux) = \frac{1}{6}\, \Delta\, \delta(\ux)$,
$$
\p_{x_i}\, \p_{x_j}\, \delta(\ux) + \onehalf\, \delta_{ij}\, \Delta\, \delta(\ux) = 15\, \invrsq\, \omega_i\, \omega_j\, \delta(\ux)
$$
We verify these formulae as follows. If $i \neq j$ then
\begin{align*}
\omega_i\, \omega_j\, \invrsq\, \delta(\ux) &= x_i\, x_j\, \frac{1}{r^4}\, \delta(\ux) = x_i\, x_j\, \frac{1}{8m(m+2)}\, \dirac^4\, \delta(\ux)\\
&= \onehalf\, \frac{1}{m(m+2)}\, x_i\, \p_{x_j}\, \dirac^2\, \delta(\ux) = \frac{1}{m(m+2)}\, \p_{x_i}\, \p_{x_j}\, \delta(\ux)
\end{align*}
which for $m=3$ reduces to
$$
\omega_i\, \omega_j\, \invrsq\, \delta(\ux)  = \frac{1}{15}\, \p_{x_i}\, \p_{x_j}\, \delta(\ux)
$$
For $i=j$ we find
\begin{align*}
\omega_i\, \omega_i\, \invrsq\, \delta(\ux) &= x_i\, x_i\, \frac{1}{r^4}\, \delta(\ux) = x_i\, x_i\, \frac{1}{8m(m+2)}\, \dirac^4\, \delta(\ux)\\
&= \onehalf\, \frac{1}{m(m+2)}\,  x_i\,  \p_{x_i}\, \dirac^2\, \delta(\ux)  = \frac{1}{2m(m+2)}\, \left( 2\,  \p_{x_i}\, \p_{x_i}\, \delta(\ux)  + \Delta\, \delta(\ux)\right)\\
&= \frac{1}{2m(m+2)}\, \Delta\, \delta(\ux) + \frac{1}{m(m+2)}\, \p_{x_i}^2\, \delta(\ux)
\end{align*}
which for $m=3$ reduces to
$$
\omega_i^2\,  \invrsq\, \delta(\ux) = \frac{1}{30}\, \Delta\, \delta(\ux) + \frac{1}{15}\, \p_{x_i}^2\, \delta(\ux)
$$

\subsection{}

In \cite{bowen} a formula is proved for 
$$
\p_{x_i}\, \left( \omega_{j_1} \cdots \omega_{j_n}\, \invrsq  \right)
$$
the inverse square field $\invrsq$ being at the heart of the discussion of idealized point sources in three dimensional space. An elegant proof based on more general formulae is given in \cite{estrada}. We show how this formula fits into our theory focussing on two special cases:
\begin{align*}
\p_{x_i}\, \left( \omega_{j}\, \invrsq  \right) &= \delta_{ij}\, \invrcub  - 3\, \omega_i\, \omega_j\, \invrcub + \frac{4\pi}{3}\, \delta_{ij}\, \delta(\ux)\\\
\p_{x_i}\, \left( \omega_{j}\, \omega_{k}\, \invrsq  \right) &= \delta_{ij}\, \omega_k\,  \invrcub + \delta_{ik}\, \omega_j\,  \invrcub - 4\, \omega_i\, \omega_j\, \omega_k\, \invrcub
\end{align*}
or, in the language of the underlying paper:
\begin{align}
\label{bowen1} \p_{x_i}\, \left( \omega_{j}\, ^{s}T_{-2}  \right) &= \delta_{ij}\, T_{-3}  - 3\, \omega_i\, \omega_j\, T_{-3} + \frac{4\pi}{3}\, \delta_{ij}\, \delta(\ux)\\\
\label{bowen2}\p_{x_i}\, \left( \omega_{j}\, \omega_{k}\, T_{-2}  \right) &= \delta_{ij}\, \omega_k\,  ^{s}T_{-3} + \delta_{ik}\, \omega_j\,  ^{s}T_{-3} - 4\, \omega_i\, \omega_j\, \omega_k\, ^{s}T_{-3}
\end{align}
It is important to remark that in (\ref{bowen2}) the expression $\omega_{j}\, \omega_{k}\, T_{-2} $ is a distribution. Indeed, it holds for the distribution $T_{-2}$ that
$$
\omega_k\, T_{-2} = \{\uom\, T_{-2}  \}^{(k)} = \{^{s}U_{-2}  \}^{(k)} = x_k\, ^{s}T_{-3}
$$
whence
$$
\omega_j\, \omega_k\, T_{-2} = \{ \uom\, x_k\, ^{s}T_{-3}\}^{(j)} = \{ x_k\, U_{-3}  \}^{(j)} = x_j\, x_k\, T_{-3}
$$
However, if in (\ref{bowen1}) $T_{-2}$ would be interpreted as a distribution, then the expression $\omega_j\, T_{-2}$ becomes a signumdistribution. In order to obtain the proposed result (\ref{bowen1}), $\omega_j\, T_{-2}$ should be a distribution whence $\invrsq$ has to be interpreted as the signumdistribution $^{s}T_{-2}$. It then holds indeed that
$$
\omega_j\, ^{s}T_{-2} = \{ \uom\, ^{s}T_{-2}   \}^{(j)} = \{  U_{-2} \}^{(j)} = x_j\, T_{-3}
$$
clearly a distribution.\\

\noindent
To prove  (\ref{bowen1}) consider the signumdistribution $^{s}T_{-m+1}$ for which it holds that 
$$
\omega_j\, ^{s}T_{-m+1} = x_j\, T_{-m}
$$
For $i \neq j$ we then have
$$
\p_{x_i} \left( \omega_j\ ^{s}T_{-m+1} \right) = \p_{x_i}\, \left( x_j\ T_{-m} \right) = x_j\, \left( -\, m\, x_i\, T_{-m-2} - \frac{1}{m}\, a_m\, \p_{x_i}\delta(\ux)  \right) = -\, m\, x_i\, x_j\, T_{-m-2}
$$
since $x_j\, \p_{x_i}\delta(\ux) = 0$. It follows that
$$
\p_{x_i} \left(  \omega_j\, ^{s}T_{-m+1} \right)  = -\, m\, \omega_i\, \omega_j\, T_{-m}
$$
which for $m=3$ turns into
$$
\p_{x_i} \left(  \omega_j\, ^{s}T_{-2} \right)  = -\, 3\, \omega_i\, \omega_j\, T_{-3}
$$
For $i = j$ we have
$$
\p_{x_i} \left( \omega_i\ ^{s}T_{-m+1} \right) = \p_{x_i}\, \left( x_i\ T_{-m} \right) = T_{-m} + x_i\, \left( -\, m\, x_i\, T_{-m-2} - \frac{1}{m}\, a_m\, \p_{x_i}\delta(\ux)  \right) = T_{-m} -\, m\, x_i^2\, T_{-m-2} + \frac{1}{m}\, a_m\, \delta(\ux)
$$
since $x_i\, \p_{x_i}\, \delta(\ux) = -\, \delta(\ux)$. Putting $m=3$ we obtain
$$
\p_{x_i} \left( \omega_i\ ^{s}T_{-m+1} \right) = T_{-3} - 3\, \omega_i^2\, T_{-3} + \frac{4\pi}{3}\, \delta(\ux)
$$
Formula (\ref{bowen2}) is proven in a similar way. If $i \neq j$ and $i \neq k$ then
$$
\p_{x_i} \left( \omega_k\, \omega_j\ ^{s}T_{-m+1} \right) = \p_{x_i}\, \left( x_k\, x_j\ T_{-m-1} \right) = x_k\, x_j\, (-m-1)\, x_i\, T_{-m-3} = -\, (m+1)\, x_i\, x_j\, x_k\, T_{-m-3} = -\, (m+1)\, \omega_i\, \omega_j\, \omega_k\, ^{s}T_{-m}
$$
since it easily seen that $\omega_i\, \omega_j\, \omega_k\, ^{s}T_{-m} = x_i\, x_j\, x_k\, T_{-m-3}$.
For $m=3$ we obtain
$$
\p_{x_i} \left( \omega_k\, \omega_j\ ^{s}T_{-2} \right) = -\, 4\, \omega_i\, \omega_j\, \omega_k\, ^{s}T_{-3}
$$
If $i = j$ and $i \neq k$ then
$$
\p_{x_i} \left( \omega_k\, \omega_i\ ^{s}T_{-m+1} \right) = \p_{x_i}\, \left( x_k\, x_i\ T_{-m-1} \right) = x_k\, T_{-m-1} + x_k\, x_i\, (-m-1)\, x_i\, T_{-m-3} = \omega_k\, ^{s}T_{-m} -\, (m+1)\, \omega_i^2\, \omega_k\, ^{s}T_{-m}
$$
which for $m=3$ reduces to
$$
\p_{x_i} \left( \omega_k\, \omega_i\ ^{s}T_{-2} \right) = \omega_k\, ^{s}T_{-3} -\, 4\, \omega_i^2\, \omega_k\, ^{s}T_{-3}
$$
Finaly if $i = j = k$ then
$$
\p_{x_i} \left( \omega_i^2\, T_{-m+1}\right) = \p_{x_i}\left(  x_i^2\, T_{-m-1}\right) = 2\, x_i\, T_{-m-1} + x_i^2\, (-m-1)\, x_i\, T_{-m-3} = 2\, \omega_i\, ^{s}T_{-m} - (m+1)\, \omega_i^3\, T_{-m}
$$
which for $m=3$ turns into
$$
\p_{x_i} \left( \omega_i^2\, T_{-2}\right) = 2\, \omega_i\, ^{s}T_{-3} - 4\, \omega_i^3\, T_{-3}
$$
These considerations show that formulae such as (\ref{bowen1}) and (\ref{bowen2}) should be handled with great care, especially with regard to the distributional or signumdistributional character of negative powers of the radial distance $r$. In this respect it holds more generally that
$$
\omega_{j_1} \cdots \omega_{j_{2n}}\, T_\la = x_{j_1} \cdots x_{j_{2n}}\, T_{\la-2n}
$$
is a distribution, while
$$
\omega_{j_1} \cdots \omega_{j_{2n+1}}\, T_\la = x_{j_1} \cdots x_{j_{2n+1}}\, ^{s}T_{\la-2n-1}
$$
is a signumdistribution.\\


\newpage
\section{Conclusion}


Expressing distributions in Euclidean space in terms of spherical co-ordinates inevitably leads to the introduction of so--called signumdistributions, i.e. bounded linear functionals on a space of indefinitely differentiable test functions showing a non--removable singularity at the origin. Operators acting on distributions appear in two kinds. The cartesian operators can be expressed using cartesian co-ordinates solely. They map distributions to distributions and their action may be either uniquely determined, as is the case for e.g. derivation with respect to the cartesian co-ordinates, multiplication by analytic functions, the Euler operator, the $\Gamma$ operator, the Laplace operator, etc., or may take the form of an equivalence class of distributions, as is the case for e.g. division by analytic functions, the operators $\uom\, \p_r$ and $\invr\, \p_{\uom}$, the operators $\uD$ and $d_j, j=1,\ldots,m$, etc. The so--called spherical operators map distributions to signumdsitributions, and their action may be either uniquely determined, as is the case for the spherical derivatives $\p_r$ and $\p_{\uom}$, or may result into an equivalence class of signumdistributions, as is the case for division by the spherical distance $r$. Due to their homogeneity and their (signum--)radial character, the distributions $T_\la$ and $U_\la$ behave specifically under the action of all those operators, in the sense that the results of these actions all are uniquely determined, except for a small number of cases, viz. $\uom\, \p_r\, U_{-m+1}$, $\invr\, \p_{\uom}\, U_{-m+1}$, $\uD\, U_{-m+1}$, $d_j\, U_{-m+1}, j=1,\ldots,m$, $\invrsq\, T_{-m+2}$ and $\invrsq\, U_{-m+1}$, which all can be traced back to the non--uniquely determined action $\invux\, U_{-m+1}$. Nevertheless we were able to uniquely determine the forementioned actions conserving overall consistency. Crucial to obtaining this result is the possibility to have a free choice of the arbitrary constant in the equivalence class  $\left[\invux\, U_{-m+1}\right]$, the non--uniquely determined character of which should be maintained.


\newpage
\section{Appendix 1: The distribution ${\rm Fp}\, x_+^\mu$ on the real line}
\label{appendixfinitepart}


Let $\mu$ be a complex parameter and $x$ a real variable. We consider the function
$$
x_+^\mu = x^\mu\, Y(x)
$$
where $Y(x)$ is the Heaviside step function.\\

\noindent
If $\Re\,  \mu > -1$ then $x_+^\mu$ is locally integrable and hence a regular distribution given, for all test functions $\phi(x) \in \mcD(\mR)$, by
$$
\langle \ x_+^\mu \ , \ \phi(x) \ \rangle = \int_0^{+ \infty}\, x^\mu\, \phi(x)\, dx \quad , \quad \Re\, \mu > -1
$$
As a function of $\mu$ the distribution $x_+^\mu$ is holomorphic in the halfplane $\Re\, \mu > -1$.\\

\noindent
For $\Re\, \mu \leq -1$ the function $x_+^\mu$ is no longer locally integrable. One introduces the so--calles {\em finite part} distribution ${\rm Fp}\, x_+^\mu$ by putting for $-n-1 < \Re\, \mu < -n$ and all $n \in \mN$
\begin{align}
\langle \ {\rm Fp}\, x_+^\mu \ , \phi(x) \ \rangle &= \int_0^{+ \infty}\, x^\mu\, \left(  \phi(x) - \phi(0) - \frac{\phi'(0)}{1!}\, x  - \cdots - \frac{\phi^{(n-1)}(0)}{(n-1)!}\, x^{n-1}\right)\, dx{\nonumber}\\[3mm]
\label{finitepart1dim} &= \lim_{\epsilon \rightarrow 0+}\, \left(  \int_\epsilon ^{+\infty}\, x^\mu\, \phi(x)\, dx + \phi(0)\, \frac{\epsilon^{\mu+1}}{\mu+1}  + \cdots + \frac{\phi^{(n-1)}(0)}{(n-1)!}\, \frac{\epsilon^{\mu+n}}{\mu+n}\right)
\end{align}

\noindent
Note that for $\mu > -1$ the defining expression for  ${\rm Fp}\, x_+^\mu$   \Eqref{finitepart1dim} reduces to that for $x_+^\mu$, and by analytic continuation  ${\rm Fp}\, x_+^\mu$ becomes holomorphic in $\mC \backslash \{-1, -2, -3, \ldots   \}$. Its singularities $\mu = -1, -2, -3,\ldots$ are simple poles with residue
$$
{\rm Res}_{\mu = -n}\, {\rm Fp}\, x_+^\mu = \frac{(-1)^{n-1}}{(n-1)!}\, \delta^{(n-1)}(x)
$$
It is still possible to define  ${\rm Fp}\, x_+^\mu$ for $\mu = -1, -2, -3, \ldots$ through the so--called {\em monomial pseudofunctions}, but it should be emphasized that by this additional defintion no entire function of $\mu$ is obtained. One puts, for all $n \in \mC$,
$$
\langle \ {\rm Fp}\, x_+^{-n} \ , \phi(x) \ \rangle = \lim_{\epsilon \rightarrow 0+}\, \left(  \int_\epsilon ^{+\infty}\, x^{-n}\, \phi(x)\, dx + \phi(0)\, \frac{\epsilon^{-n+1}}{-n+1}  + \cdots + \frac{\phi^{(n-2)}(0)}{(n-2)!}\, \frac{\epsilon^{-1}}{-1} + \frac{\phi^{(n-1)}(0)}{(n-1)!}\, \ln \epsilon\right)
$$

\begin{lemma}
The distribution ${\rm Fp}\, x_+^\mu$ obeys the multiplication rules
\begin{itemize}
\item $x\, {\rm Fp}\, x_+^\mu = {\rm Fp}\, x_+^{\mu+1} \quad , \quad \mu \neq -1, -2, -3, \dots$
\item $x\, {\rm Fp}\, x_+^{-1} = Y(x)$
\item $x\, {\rm Fp}\, x_+^{-n} =  {\rm Fp}\, x_+^{-n+1} \quad , \quad n = 2, 3, \ldots$
\end{itemize}
\end{lemma}

\begin{lemma}
The distribution ${\rm Fp}\, x_+^\mu$ obeys the differentiation rules
\begin{itemize}
\item $\p_x\, {\rm Fp}\, x_+^\mu = \mu\,  {\rm Fp}\, x_+^{\mu-1} \quad , \quad \mu \neq 0, -1, -2,  \dots$
\item $\p_x\, {\rm Fp}\, x_+^{0} = \delta(x)$
\item $\p_x\, {\rm Fp}\, x_+^{-n} =  (-n)\, {\rm Fp}\, x_+^{-n-1} + (-1)^n\, \frac{1}{n!}\, \delta^{(n)}(x)\quad , \quad n = 1, 2, 3, \ldots$
\end{itemize}
\end{lemma}

\begin{lemma}
If the test function $\phi(x)$ is such that $\phi(0) = \phi'(0) = \ldots = \phi^{(k)}(0) = 0$ then
\begin{itemize}
\item $\langle \  {\rm Fp}\, x_+^\mu \ , \ \frac{1}{x^k}\, \phi(x) \ \rangle = \langle \  {\rm Fp}\, x_+^{\mu-k} \ , \ \phi(x) \ \rangle \quad , \quad \mu \neq 0, 1, 2, \ldots , k-1$
\item $\langle \  {\rm Fp}\, x_+^n \ , \ \frac{1}{x^k}\, \phi(x) \ \rangle = \langle \  {\rm Fp}\, x_+^{\mu-k} \ , \ \phi(x) \ \rangle \quad , \quad n = 0, 1, 2, \ldots , k-1$
\end{itemize}
\end{lemma}


\newpage
\section{Appendix 2: The distributions $T_\la$ and $U_\la$ in $\mR^m$}
\label{appendixTU}


Let $\la$ be a complex parameter, let $\mu = \la+m-1$, let Fp\,$r_+^\mu$ be the finite part distribution on the real $r$--axis, defined by the monomial pseudofunctions for $\mu = -\, n, n \in \mN$, and let $\varphi(\ux)$ be a scalar test function in $\mR^m$. Then we define the following two families of distributions in $\mR^m$.

\begin{definition}
One puts, for $\la \in \mC$,
$$
\langle \ T_\la \ , \ \varphi(\ux) \ \rangle = a_m\, \langle \ {\rm Fp}\, r_+^\mu \ , \ \Sigma^{(0)}\, [\varphi](r) \ \rangle
$$
and
$$
\langle \ U_\la \ , \ \varphi(\ux) \ \rangle = a_m\, \langle \ {\rm Fp}\, r_+^\mu \ , \ \Sigma^{(1)}\, [\varphi](r) \ \rangle
$$
\end{definition}

\noindent
The scalar distributions $T_\la$ are classical distributions in harmonic analysis. The $1$--vector distributions $U_\la$ are the Clifford analysis counterparts of vector valued distributions, see the history comment in \cite{BDSch}. Both families of distributions were thoroughly studied in the framework of Clifford analysis, see the references mentioned in the introductory Section \ref{intro}.\\

Seen the holomorphy properties of the one--dimensional distribution Fp\, $r_+^\mu$, we expect $T_\la$ to be {\em a priori} a meromorphic function of the complex paramter $\la$ with simple poles at $\mu = -1, -2, \ldots$, or $\la = -m, -m-1, -m-2, \ldots$ with residues given by
\begin{align*}
{\rm Res }_{\, \la=-m-2\ell+1} \,  \langle \ T_\la \ , \ \varphi(\ux) \ \rangle  &=  a_m\, \langle \ -\, \frac{1}{(2\ell-1)!}\, \delta^{(2\ell-1)}(r) \ , \ \Sigma^{(0)}\, [\varphi](r) \ \rangle\\[3mm]
{\rm Res }_{\, \la=-m-2\ell} \,  \langle \ T_\la \ , \ \varphi(\ux) \ \rangle  &=  a_m\, \langle \ \frac{1}{(2\ell)!}\, \delta^{(2\ell)}(r) \ , \ \Sigma^{(0)}\, [\varphi](r) \ \rangle
\end{align*}
As the spherical mean $\Sigma^{(0)}\, [\varphi](r)$ has vanishing odd order derivatives at the origin $r=0$, it follows that the singularities of $T_\la$ at the points $\la= -m-2\ell+1, \ell=1,2,\ldots$ are removable. This should be confirmed by coinciding expressions for $T_\la$ in the vertical strips $-m-2\ell-1 < \Re\, \la < -m-2\ell$ and $-m-2\ell-2 < \Re\, \la < -m-2\ell-1$, and a matching value at the removable singularity $\la=-m-2\ell-1$. We indeed have, using the short--hand notation $\Phi$ for the spherical mean $\Sigma^{(0)}[\varphi]$,

\begin{itemize}
\item for $-m-2\ell-1 < \Re\, \la < -m-2\ell$ or $-2\ell - 2 < \Re\, \mu < -2\ell-1$
\begin{align*}
\frac{1}{a_m}\, \langle \ T_\la \ , \ \varphi \ \rangle &= \lim_{\epsilon \rightarrow +0}\, \int_\epsilon^{+\infty}\, r^\mu\, \Phi(r)\, dr + \Phi(0)\, \frac{\epsilon^{\mu+1}}{\mu+1} + \frac{\Phi'(0)}{1!}\, \frac{\epsilon^{\mu+2}}{\mu+2} + \cdots + \frac{\Phi^{(2\ell-1)}(0)}{(2\ell-1)!}\, \frac{\epsilon^{\mu+2\ell}}{\mu+2\ell} +  \frac{\Phi^{(2\ell)}(0)}{(2\ell)!}\, \frac{\epsilon^{\mu+2\ell+1}}{\mu+2\ell+1}\\[2mm]
&= \lim_{\epsilon \rightarrow +0}\, \int_\epsilon^{+\infty}\, r^\mu\, \Phi(r)\, dr + \Phi(0)\, \frac{\epsilon^{\mu+1}}{\mu+1} + \frac{\Phi''(0)}{2!}\, \frac{\epsilon^{\mu+3}}{\mu+3} + \cdots +  \frac{\Phi^{(2\ell)}(0)}{(2\ell)!}\, \frac{\epsilon^{\mu+2\ell+1}}{\mu+2\ell+1}
\end{align*}

\item for $-m-2\ell-2 < \Re\, \la < -m-2\ell-1$ or $-2\ell - 3 < \Re\, \mu < -2\ell-2$
\begin{align*}
\frac{1}{a_m}\, \langle \ T_\la \ , \ \varphi \ \rangle &= \lim_{\epsilon \rightarrow +0}\, \int_\epsilon^{+\infty}\, r^\mu\, \Phi(r)\, dr + \Phi(0)\, \frac{\epsilon^{\mu+1}}{\mu+1} + \frac{\Phi'(0)}{1!}\, \frac{\epsilon^{\mu+2}}{\mu+2} + \cdots + \frac{\Phi^{(2\ell)}(0)}{(2\ell)!}\, \frac{\epsilon^{\mu+2\ell+1}}{\mu+2\ell+1} +  \frac{\Phi^{(2\ell+1)}(0)}{(2\ell+1)!}\, \frac{\epsilon^{\mu+2\ell+2}}{\mu+2\ell+2}\\[2mm]
&= \lim_{\epsilon \rightarrow +0}\, \int_\epsilon^{+\infty}\, r^\mu\, \Phi(r)\, dr + \Phi(0)\, \frac{\epsilon^{\mu+1}}{\mu+1} + \frac{\Phi''(0)}{2!}\, \frac{\epsilon^{\mu+3}}{\mu+3} + \cdots +  \frac{\Phi^{(2\ell)}(0)}{(2\ell)!}\, \frac{\epsilon^{\mu+2\ell+1}}{\mu+2\ell+1}
\end{align*}

\item for $\la=-m-2\ell-1$ or $\mu = -2\ell-2$
\begin{align*}
\frac{1}{a_m}\, \langle \ T_{-m-2\ell-1} \ , \ \varphi \ \rangle &= \lim_{\epsilon \rightarrow +0}\, \int_\epsilon^{+\infty}\, r^{-2\ell-2}\, \Phi(r)\, dr + \Phi(0)\, \frac{\epsilon^{-2\ell-1}}{-2\ell-1} + \frac{\Phi'(0)}{1!}\, \frac{\epsilon^{-2\ell}}{-2\ell} + \cdots + \frac{\Phi^{(2\ell)}(0)}{(2\ell)!}\, \frac{\epsilon^{-1}}{-1} +  \frac{\Phi^{(2\ell+1)}(0)}{(2\ell+1)!}\, \ln{\epsilon}\\[2mm]
&= \lim_{\epsilon \rightarrow +0}\, \int_\epsilon^{+\infty}\, r^{-2\ell-2}\, \Phi(r)\, dr + \Phi(0)\, \frac{\epsilon^{-2\ell-1}}{-2\ell-1} +  \cdots + \frac{\Phi^{(2\ell)}(0)}{(2\ell)!}\, \frac{\epsilon^{-1}}{-1} 
\end{align*}
\end{itemize}

\noindent
We may conclude that the distribution $T_\la$ is a meromorphic function of the complex parameter $\la$ showing simple poles at $\la = -m-2\ell, \ell=0,1,2,\ldots$ where there is an {\em ad hoc} defintion through the monomial pseudofunctions. Let us compute the residue at these singularities:
\begin{align*}
{\rm Res}_{\, \la=-m-2\ell}\, \langle \ T_\la \ , \ \varphi(\ux) \ \rangle &=    {\rm Res}_{\, \mu=-2\ell-1} \, a_m\, \langle \    {\rm Fp}\, r_+^\mu \ , \ \Sigma^{(0)} [\varphi](r) \ \rangle\\[2mm]
&= \frac{1}{(2\ell)!}\, a_m\, \langle \   \delta^{(2\ell)}(r)  \ , \ \Sigma^{(0)} [\varphi](r) \ \rangle\\[2mm]
&= (-1)^\ell\, \frac{m+2\ell}{C(m,\ell)}\, a_m\, \{ \dirac^{2\ell}\, \varphi(\ux) \}|_{\ux=0}
\end{align*}
whence
$$
{\rm Res}_{\, \la=-m-2\ell}\,  T_\la = (-1)^\ell\, \frac{m+2\ell}{C(m,\ell)}\, a_m\,  \dirac^{2\ell}\, \delta(\ux)
$$\hfill\\

In a similar way we  eventually find that the distribution $U_\la$ is a meromorphic function of the complex parameter $\la$ showing simple poles at $\la = -m-2\ell-1, \ell=0,1,2,\ldots$ where there is an {\em ad hoc} definition through the monomial pseudofunctions. 
For the residue at these singularities we have
\begin{align*}
{\rm Res}_{\, \la=-m-2\ell-1}\, \langle \ U_\la \ , \ \varphi(\ux) \ \rangle &=    {\rm Res}_{\, \mu=-2\ell-2} \, a_m\, \langle \    {\rm Fp}\, r_+^\mu \ , \ \Sigma^{(1)} [\varphi](r) \ \rangle\\[2mm]
&= -\, \frac{1}{(2\ell+1)!}\, a_m\, \langle \   \delta^{(2\ell+1)}(r)  \ , \ \Sigma^{(1)} [\varphi](r) \ \rangle\\[2mm]
&= (-1)^{\ell+1}\, \frac{1}{C(m,\ell)}\, a_m\, \{ \dirac^{2\ell+1}\, \varphi(\ux) \}|_{\ux=0}
\end{align*}
whence
$$
{\rm Res}_{\, \la=-m-2\ell-1}\,  U_\la = (-1)^{\ell+1}\, \frac{1}{C(m,\ell)}\, a_m\,  \dirac^{2\ell+1}\, \delta(\ux)
$$


\newpage
\section{Appendix 3: Some specific distributions}
\label{appendixspecific}


In this appendix we give an overview of the properties of some specific distributions of the $T_\la$ and $U_\la$ families.\\

\vspace{1cm}

$$
\boxed{  \ \ T_0 \ \ }
$$\hfil\\

\noindent
One has
\begin{align*}
\langle \ T_0 \ , \ \varphi(\ux) \ \rangle &= a_m\, \langle \ {\rm Fp}\, r_+^{m-1} \ , \ \Sigma^{(0)}[\varphi](r) \ \rangle = a_m\, \int_{0}^{+\infty}\, r^{m-1}\, dr\, \frac{1}{a_m}\, \int_{S^{m-1}}\, \varphi(r\, \uom)\, dS_{\uom} = \int_{\mR^m}\, \varphi(\ux) d\ux \\[3mm] 
&= \langle \ 1 \, , \ \varphi(\ux) \ \rangle
\end{align*}

\noindent
whence 
$T_0 = 1$ is the regular distribution associated to the constant function $1 = \ux^0$. Its properties are trivial.

\newpage

$$
\boxed{  \ \ U_0 \ \ }
$$\hfil\\

\noindent
One has
\begin{align*}
\langle \ U_0 \ , \ \varphi(\ux) \ \rangle &= a_m\, \langle \ {\rm Fp}\, r_+^{m-1} \ ,  \ \Sigma^{(1)}[\varphi](r) \ \rangle = a_m\, \int_{0}^{+\infty}\, r^{m-1}\, dr\, \frac{1}{a_m}\, \int_{S^{m-1}}\, \uom\, \varphi(r\, \uom)\,  dS_{\uom} = \int_{\mR^m}\, \uom\, \varphi(\ux) d\ux \\[3mm] 
&= \langle \ \uom \, , \ \varphi(\ux) \ \rangle
\end{align*}

\noindent
whence 
$U_0 = \uom$ is the regular distribution associated to the function $\uom = \frac{\ux}{|\ux|}$, which is the higher dimensional counterpart to the {\em signum distribution} on the real line. Its properties are:
\begin{itemize}
\item $\ux\, U_0 = -\, T_1 = -\, r$
\item $\dirac\, U_0 = -\, (m-1)\, T_{-1}$ \ \  \\[2mm]which splits up into
$$
(\uom\, \p_r)\, U_0 = 0 \quad {\rm and} \quad (\invr\, \p_{\uom})\, U_0 = -\, (m-1)\, T_{-1}
$$
\item $\Delta\, U_0 = -\, (m-1)\, U_{-2}$
\item $\mE\, U_0 = 0$
\item $\Gamma\, U_0 = (m-1)\, U_0$ \ \ and \ \ $\p_{\uom}\, U_0 = -\, (m-1)$
\item $\Delta^*\, U_0 = -\, (m-1)\, \uom$ \ \ and \ \ $\invrsq\, \Delta^*\, U_0 = -\, (m-1)\, U_{-2}$
\item $\p_r^2\, U_0 = 0$ \ \ and \ \ $\invr\, \p_r\, U_0 = 0$
\item $\uD\, U_0 = 0$
\item $\bZ\, U_0 = 0$
\item $\invrsq\, \bZ^*\, U_0 = 0$

\end{itemize}

\newpage

$$
\boxed{  \ \ T_{-m+2} \ \ }
$$\hfil\\

\noindent
One has
\begin{align*}
\langle \ T_{-m+2} \ , \ \varphi(\ux) \ \rangle &= a_m\, \langle \ {\rm Fp}\, r_+ \  , \  \Sigma^{(0)}[\varphi](r) \ \rangle = a_m\, \int_{0}^{+\infty}\, r\, dr\, \frac{1}{a_m}\, \int_{S^{m-1}}\,  \varphi(r\, \uom)\,  dS_{\uom} = \int_{\mR^m}\, \frac{1}{r^{m-2}}\, \varphi(\ux) d\ux \\[3mm] 
&= \langle \  \frac{1}{r^{m-2}} \ , \ \varphi(\ux) \ \rangle
\end{align*}

\noindent
whence 
$T_{-m+2}$ is the regular distribution associated to the locally integrable function $\frac{1}{r^{m-2}}$ in $\mR^m$. Its properties are:
\begin{itemize}
\item $\ux\, T_{-m+2} = U_{-m+3}$
\item $\dirac\, T_{-m+2} = -\, (m-2)\, U_{-m+1}$ \ \  \\[2mm]which splits up into
$$
(\uom\, \p_r)\, T_{-m+2} = -\, (m-2)\, U_{-m+1} \quad {\rm and} \quad (\invr\, \p_{\uom})\, T_{-m+2} = 0
$$
\item $\Delta\, T_{-m+2}= -\, a_m\, (m-2)\, \delta(\ux)$\\[2mm]
which expresses the fact that $-\, \frac{1}{a_m}\, \frac{1}{m-2}\, \frac{1}{r^{m-2}} $ is the fundamental solution of the Laplace operator
\item $\mE\, T_{-m+2} = -\, (m-2)\, T_{-m+2}$
\item $\Gamma\, T_{-m+2}= 0$
\item $\Delta^*\, T_{-m+2} = 0$
\item $\p_r^2\, T_{-m+2} = (m-2)(m-1)\, T_{-m} - a_m\, (m-2)\, \delta(\ux)$ 
\item $\invr\, \p_r\,  T_{-m+2} = -\, (m-2)\, T_{-m}$
\item $\uD\, T_{-m+2} = U_{-m+1}$
\item $\bZ\, T_{-m+2} = -\, (m-1)\, T_{-m} + a_m\, \delta(\ux)$
\item $\invrsq\, \bZ^*\, T_{-m+2} = -\, (m-1)\, T_{-m} + (m-1)\, a_m\, \delta(\ux)$
\end{itemize}

\newpage

$$
\boxed{  \ \ T_{-m+1} \ \ }
$$\hfil\\

\noindent
One has
\begin{align*}
\langle \ T_{-m+1} \ , \ \varphi(\ux) \ \rangle &= a_m\, \langle \ {\rm Fp}\, r_+^0 \  , \  \Sigma^{(0)}[\varphi](r) \ \rangle = a_m\, \int_{0}^{+\infty}\,  dr\, \frac{1}{a_m}\, \int_{S^{m-1}}\,  \varphi(r\, \uom)\,  dS_{\uom} = \int_{\mR^m}\, \frac{1}{r^{m-1}}\, \varphi(\ux) d\ux \\[3mm] 
&= \langle \  \frac{1}{r^{m-1}} \ , \ \varphi(\ux) \ \rangle
\end{align*}

\noindent
whence 
$T_{-m+1}$ is the regular distribution associated to the locally integrable function $\frac{1}{r^{m-1}}$ in $\mR^m$. Its properties are:

\begin{itemize}
\item $\ux\, T_{-m+1} = U_{-m+2}$
\item $\dirac\, T_{-m+1} = -\, (m-1)\, U_{-m}$ \ \  \\[2mm]which splits up into
$$
(\uom\, \p_r)\, T_{-m+1} = -\, (m-1)\, U_{-m} \quad {\rm and} \quad (\invr\, \p_{\uom})\, T_{-m+1} = 0
$$
\item $\Delta\, T_{-m+1}= (m-1)\, T_{-m-1}$
\item $\mE\, T_{-m+1} = -\, (m-1)\, T_{-m+1}$
\item $\Gamma\, T_{-m+1}= 0$
\item $\Delta^*\, T_{-m+1} = 0$
\item $\p_r^2\, T_{-m+1} = (m-1)(m)\, T_{-m-1}$ 
\item $\invr\, \p_r\,  T_{-m+1} = -\, (m-1)\, T_{-m-1}$
\item $\uD\, T_{-m+1} = 0$
\item $\bZ\, T_{-m+1} = 0$
\item $\invrsq\, \bZ^*\, T_{-m+1} = -\, (m-1)\, T_{-m-1}$
\end{itemize}

\newpage

$$
\boxed{  \ \ U_{-m+1} \ \ }
$$\hfil\\

\noindent
One has
\begin{align*}
\langle \ U_{-m+1} \ , \ \varphi(\ux) \ \rangle &= a_m\, \langle \ {\rm Fp}\, r_+^0 \  , \  \Sigma^{(1)}[\varphi](r) \ \rangle = a_m\, \int_{0}^{+\infty}\,  dr\, \frac{1}{a_m}\, \int_{S^{m-1}}\,  \uom\, \varphi(r\, \uom)\,  dS_{\uom} = \int_{\mR^m}\, \frac{1}{r^{m-1}}\, \uom\, \varphi(\ux) d\ux \\[3mm] 
&= \langle \  \frac{\uom}{r^{m-1}} \ , \ \varphi(\ux) \ \rangle
\end{align*}

\noindent
whence 
$U_{-m+1}$ is the regular distribution associated to the locally integrable function $\frac{1}{r^{m-1}}\, \uom$ in $\mR^m$. Its properties are:

\begin{itemize}
\item $\ux\, U_{-m+1} = -\, T_{-m+2}$
\item $\dirac\, U_{-m+1} = -\, a_m\, \delta(\ux)$ \ \  \\[2mm]which splits up into
$$
(\uom\, \p_r)\, U_{-m+1} = (m-1)\, T_{-m} - a_m\, \delta(\ux) \quad {\rm and} \quad (\invr\, \p_{\uom})\, U_{-m+1} = -\, (m-1)\, T_{-m}
$$
and expresses the fact that $-\, \frac{1}{a_m}\,  \frac{1}{r^{m-1}}\, \uom = -\, \frac{1}{a_m}\,  \frac{1}{r^{m}}\, \ux $ is the fundamental solution of the Dirac operator $\dirac$
\item $\Delta\, U_{-m+1}= a_m\, \dirac \delta(\ux)$
\item $\mE\, U_{-m+1} = -\, (m-1)\, U_{-m+1}$
\item $\Gamma\, U_{-m+1}= (m-1)\, U_{-m+1}$
\item $\Delta^*\, U_{-m+1} = -\, (m-1)\, U_{-m+1}$ \ \ and \ \ $\invrsq\, \Delta^*\, U_{-m+1} = -\, (m-1)\, U_{-m-1}$
\item $\p_r^2\, U_{-m+1} = (m-1)(m)\, U_{-m-1} + \frac{2m-1}{m}\, a_m\, \dirac \delta(\ux)$ 
\item $\invr\, \p_r\,  U_{-m+1} = -\, (m-1)\, U_{-m-1} - \frac{1}{m}\, a_m\, \dirac \delta(\ux)$
\item $\uD\, U_{-m+1} = (m-1)\, T_{-m} - a_m\, \delta(\ux)$
\item $\bZ\, U_{-m+1} =  (m-1)\, U_{-m-1} + a_m\, \dirac \delta(\ux)$
\item $\invrsq\, \bZ^*\, U_{-m+1} = 0$
\end{itemize}

\newpage

$$
\boxed{  \ \ T_{-m} \ \ }
$$\hfil\\

\noindent
One has
\begin{align*}
\langle \ T_{-m} \ , \ \varphi(\ux) \ \rangle &= a_m\, \langle \ {\rm Fp}\, r_+^{-1} \  , \  \Sigma^{(0)}[\varphi](r) \ \rangle =  {\rm Fp}\, \int_{0}^{+\infty}\,  \invr\, dr\,  \int_{S^{m-1}}\,  \varphi(r\, \uom)\,  dS_{\uom}  \\[3mm] 
&= a_m\, \lim_{\epsilon \rightarrow 0+}\, \int_\epsilon^{+\infty}\, \invr\, \Sigma^{(0)}[\varphi](r) + \varphi(0)\, \ln{\epsilon}
\end{align*}
\noindent
Properties are:

\begin{itemize}
\item $\ux\, T_{-m} = U_{-m+1}$
\item $\dirac\, T_{-m} = -\, m\, U_{-m-1} -\, a_m\, \frac{1}{m}\, \dirac \delta(\ux)$ \ \  \\[2mm]which splits up into
$$
(\uom\, \p_r)\, T_{-m} = -\, m\, U_{-m-1} -\, a_m\, \frac{1}{m}\, \dirac \delta(\ux)\quad {\rm and} \quad (\invr\, \p_{\uom})\, T_{-m} = 0
$$
\item $\Delta\, T_{-m}= 2\, m\, T_{-m-2} + \frac{m+2}{2m}\, a_m\, \dirac^2\, \delta(\ux)$
\item $\mE\, T_{-m}= -\, m\, T_{-m} + a_m\, \delta(\ux)$
\item $\Gamma\, T_{-m}= 0$
\item $\Delta^*\, T_{-m} = 0$
\item $\p_r^2\, T_{-m} = m(m+1)\, T_{-m-2}  + \frac{2m+1}{2m}\, a_m\, \dirac^2\,  \delta(\ux)$ 
\item $\invr\, \p_r\,  T_{-m}= -\, m\, T_{-m-2} - \frac{1}{2m}\, a_m\, \dirac^2\,  \delta(\ux)$
\item $\uD\, T_{-m} = -\, U_{-m-1} - \frac{1}{m}\, a_m\, \dirac \delta(\ux)$
\item $\bZ\, T_{-m} = (m+1)\, T_{-m-2} + \frac{m+2}{2m}\, a_m\, \dirac^2 \delta(\ux)$
\item $\invrsq\, \bZ^*\, T_{-m} = -\, (m-1)\, T_{-m-2}$

\end{itemize}

\newpage

$$
\boxed{  \ \ U_{-m} \ \ }
$$\hfil\\

\noindent
One has
\begin{align*}
\langle \ U_{-m} \ , \ \varphi(\ux) \ \rangle &= a_m\, \langle \ {\rm Fp}\, r_+^{-1} \  , \  \Sigma^{(1)}[\varphi](r) \ \rangle  = a_m\, \lim_{\epsilon \rightarrow 0+}\, \int_\epsilon^{+\infty}\, \invr\, \Sigma^{(1)}[\varphi](r) + \Sigma^{(1)}[\varphi](0)\, \ln{\epsilon}\\[3mm]
&=  \lim_{\epsilon \rightarrow 0+}\, \int_\epsilon^{+\infty}\, \invr\, dr\, \int_{S^{m-1}}\, \uom\, \varphi(r\, \uom)\, dS_{\uom} = \lim_{\epsilon \rightarrow 0+}\, \int_{\mR^m \backslash B(0,\epsilon)}\, \frac{\uom}{r^m}\, \varphi(\ux)\, d\ux \\[3mm]
&= {\rm Pv}\, \int_{\mR^m }\, \frac{\uom}{r^m}\, \varphi(\ux)\, d\ux = \langle \ {\rm Pv}\, \frac{\uom}{r^m}  \ , \ \varphi(\ux) \ \rangle
\end{align*}

\noindent
whence 
$$
U_{-m} = {\rm Pv}\, \frac{\uom}{r^m}\, 
$$
which is the higher dimensional analogue of the {\em principal value} distribution Pv $\frac{1}{t}$ on the real line. It is, up to a constant, the convolution kernel $\mcH$ of the Hilbert transform, given, for a suitable function or distribution $f$, by
$$
H[f] = \mcH * f =  - \frac{2}{a_{m+1}}\, {\rm Pv}\, \frac{\uom}{r^m}\, * f
$$

\noindent
Properties are

\begin{itemize}
\item $\ux\, U_{-m} = -\, T_{-m+1}$ \ or \  $(r\, \uom)\, {\rm Pv}\, \frac{\uom}{r^m} = -\, \frac{1}{r^{m-1}}$
\item $\dirac\, U_{-m} = T_{-m-1}$ \quad
which splits up into \quad
$
(\uom\, \p_r)\, U_{-m} =    m\, T_{-m-1}  \quad {\rm and} \quad (\invr\, \p_{\uom})\, U_{-m} = -\, (m-1)\, T_{-m-1}
$
\item $\Delta\, U_{-m}= (m+1)\, U_{-m-2}$
\item $\mE\, U_{-m}= -\, m\, U_{-m} $
\item $\Gamma\, U_{-m}= (m-1)\, U_{-m}$
\item $\Delta^*\, U_{-m} = -\, (m-1)\, U_{-m} \quad {\rm and}   \quad \invrsq\,  \Delta^*\, U_{-m} = -\, (m-1)\, U_{-m-2}$
\item $\p_r^2\, U_{-m} = m(m+1)\, U_{-m-2}  $ 
\item $\invr\, \p_r\,  U_{-m} = -\, m\, U_{-m-2} $
\item $\uD\, U_{-m}= m\, T_{-m-1}$
\item $\bZ\, U_{-m} = 2\, m\, \, U_{-m-2}$
\item $\invrsq\, \bZ^*\, U_{-m} = 0$
\end{itemize}

\newpage

$$
\boxed{  \ \ T_{-m-1} \ \ }
$$\hfil\\

\noindent
One has
\begin{align*}
\langle \ T_{-m-1} \ , \ \varphi(\ux) \ \rangle &= a_m\, \langle \ {\rm Fp}\, r_+^{-2} \  , \  \Sigma^{(0)}[\varphi](r) \ \rangle  \\[3mm] 
&=  {\rm Fp}\, \int_{\mR^m}\,  \frac{1}{r^{m+1}}\,  \varphi(\ux)\,  d\ux  =  \langle \ {\rm Fp}\, \frac{1}{r^{m+1}}  \ , \ \varphi(\ux) \ \rangle
\end{align*}
The scalar distribution $T_{-m-1} = {\rm Fp}\, \frac{1}{r^{m+1}}$ is, up to a constant, the convolution kernel for the so--called {\em square root of the negative Laplacian}, which, for an appropriate function  or distribution $f$, is defined by
$$
\left(  -\, \Delta \right)^{\onehalf}\, [f] = -\, \frac{2}{a_{m+1}}\, {\rm Fp}\, \frac{1}{r^{m+1}} * f
$$

\noindent
Properties are:

\begin{itemize}
\item $\ux\, T_{-m-1} = U_{-m}$
\item $\dirac\, T_{-m-1} = -\, (m+1)\, U_{-m-2} $
\item $\Delta\, T_{-m-1}= 3\, (m+1)\, T_{-m-3} $
\item $\mE\, T_{-m-1}= -\, (m+1)\, T_{-m-1} $
\item $\Gamma\, T_{-m-1}= 0$
\item $\Delta^*\, T_{-m-1} = 0$
\item $\p_r^2\, T_{-m-1} = (m+1)(m+2)\, T_{-m-3}  $ 
\item $\invr\, \p_r\,  T_{-m-1}= -\, (m+1)\, T_{-m-3} $
\item $\uD\, T_{-m-1} = -2\, U_{-m-2}$
\item $\bZ\, T_{-m-1} = 2\, (m+2)\, \, T_{-m-3}$
\item $\invrsq\, \bZ^*\, T_{-m-1}= -\, (m-1)\, T_{-m-3}$
\end{itemize}\hfill\\

\noindent
Note that in \cite{convolution} the so--called Hilbert--Dirac operator $HD$ was introduced as the convolution operator with convolution kernel $\dirac\, \mcH = \mcH\, \dirac$. As it holds that $\dirac\, U_{-m} = T_{-m-1}$ or
$$
\dirac\, \left( - \frac{2}{a_{m+1}}\, {\rm Pv}\, \frac{\uom}{r^m} \right) = -\, \frac{2}{a_{m+1}}\, {\rm Fp}\, \frac{1}{r^{m+1}}
$$
it becomes clear that the convolution kernel of the operator $\left(  -\, \Delta \right)^{\onehalf}$ is precisely the convolution kernel $\dirac\, \mcH = \mcH\, \dirac$ of the Hilbert--Dirac operator $HD$, whence
$$
HD[f] = \dirac\, \mcH\, * f = \mcH\, \dirac\, * f = \left( -\, \frac{2}{a_{m+1}}\, {\rm Fp}\, \frac{1}{r^{m+1}} \right)\, * f = \left(  -\, \Delta \right)^{\onehalf}\, [f]
$$



\end{document}